\input amstex
%

\def\next{AMS-SEKR}\ifx\styname\next \endinput\fi
\catcode`\@=11
\def\styname{AMS-SEKR}
\def\styversion{2.0}
{\W@{}\W@{\styname.STY - Version \styversion}\W@{}}
\hyphenation{acad-e-my acad-e-mies af-ter-thought anom-aly anom-alies
an-ti-deriv-a-tive an-tin-o-my an-tin-o-mies apoth-e-o-ses apoth-e-o-sis
ap-pen-dix ar-che-typ-al as-sign-a-ble as-sist-ant-ship as-ymp-tot-ic
asyn-chro-nous at-trib-uted at-trib-ut-able bank-rupt bank-rupt-cy
bi-dif-fer-en-tial blue-print busier busiest cat-a-stroph-ic
cat-a-stroph-i-cally con-gress cross-hatched data-base de-fin-i-tive
de-riv-a-tive dis-trib-ute dri-ver dri-vers eco-nom-ics econ-o-mist
elit-ist equi-vari-ant ex-quis-ite ex-tra-or-di-nary flow-chart
for-mi-da-ble forth-right friv-o-lous ge-o-des-ic ge-o-det-ic geo-met-ric
griev-ance griev-ous griev-ous-ly hexa-dec-i-mal ho-lo-no-my ho-mo-thetic
ideals idio-syn-crasy in-fin-ite-ly in-fin-i-tes-i-mal ir-rev-o-ca-ble
key-stroke lam-en-ta-ble light-weight mal-a-prop-ism man-u-script
mar-gin-al meta-bol-ic me-tab-o-lism meta-lan-guage me-trop-o-lis
met-ro-pol-i-tan mi-nut-est mol-e-cule mono-chrome mono-pole mo-nop-oly
mono-spline mo-not-o-nous mul-ti-fac-eted mul-ti-plic-able non-euclid-ean
non-iso-mor-phic non-smooth par-a-digm par-a-bol-ic pa-rab-o-loid
pa-ram-e-trize para-mount pen-ta-gon phe-nom-e-non post-script pre-am-ble
pro-ce-dur-al pro-hib-i-tive pro-hib-i-tive-ly pseu-do-dif-fer-en-tial
pseu-do-fi-nite pseu-do-nym qua-drat-ics quad-ra-ture qua-si-smooth
qua-si-sta-tion-ary qua-si-tri-an-gu-lar quin-tes-sence quin-tes-sen-tial
re-arrange-ment rec-tan-gle ret-ri-bu-tion retro-fit retro-fit-ted
right-eous right-eous-ness ro-bot ro-bot-ics sched-ul-ing se-mes-ter
semi-def-i-nite semi-ho-mo-thet-ic set-up se-vere-ly side-step sov-er-eign
spe-cious spher-oid spher-oid-al star-tling star-tling-ly
sta-tis-tics sto-chas-tic straight-est strange-ness strat-a-gem strong-hold
sum-ma-ble symp-to-matic syn-chro-nous topo-graph-i-cal tra-vers-a-ble
tra-ver-sal tra-ver-sals treach-ery turn-around un-at-tached un-err-ing-ly
white-space wide-spread wing-spread wretch-ed wretch-ed-ly Brown-ian
Eng-lish Euler-ian Feb-ru-ary Gauss-ian Grothen-dieck Hamil-ton-ian
Her-mit-ian Jan-u-ary Japan-ese Kor-te-weg Le-gendre Lip-schitz
Lip-schitz-ian Mar-kov-ian Noe-ther-ian No-vem-ber Rie-mann-ian
Schwarz-schild Sep-tem-ber
form per-iods Uni-ver-si-ty cri-ti-sism for-ma-lism}
\Invalid@\nofrills
\Invalid@\usualspace
\newif\ifnofrills@
\def\nofrills@#1#2{\relaxnext@
  \DN@{\ifx\next\nofrills
    \nofrills@true\let#2\relax\DN@\nofrills{\nextii@}%
  \else
    \nofrills@false\def#2{#1}\let\next@\nextii@\fi
\next@}}
\def\usualspace@#1{\ifnofrills@\def\usualspace{#1}\fi}
\def\addto#1#2{\csname \expandafter\eat@\string#1@\endcsname
  \expandafter{\the\csname \expandafter\eat@\string#1@\endcsname#2}}
\newdimen\bigsize@
\def\big@#1#2{{\hbox{$\left#2\vcenter to#1\bigsize@{}%
  \right.\nulldelimiterspace\z@\m@th$}}}
\def\big{\big@\@ne}
\def\Big{\big@{1.5}}
\def\bigg{\big@\tw@}
\def\Bigg{\big@{2.5}}
\def\raggedcenter@{\leftskip\z@ plus.4\hsize \rightskip\leftskip
 \parfillskip\z@ \parindent\z@ \spaceskip.3333em \xspaceskip.5em
 \pretolerance9999\tolerance9999 \exhyphenpenalty\@M
 \hyphenpenalty\@M \let\\\linebreak}
\def\upperspecialchars{\def\ss{SS}\let\i=I\let\j=J\let\ae\AE\let\oe\OE
  \let\o\O\let\aa\AA\let\l\L}
\def\uppercasetext@#1{%
  {\spaceskip1.2\fontdimen2\the\font plus1.2\fontdimen3\the\font
   \upperspecialchars\uctext@#1$\m@th\aftergroup\eat@$}}
\def\uctext@#1$#2${\endash@#1-\endash@$#2$\uctext@}
\def\endash@#1-#2\endash@{\uppercase{#1}\if\notempty{#2}--\endash@#2\endash@\fi}
\def\runaway@#1{\DN@{#1}\ifx\envir@\next@
  \Err@{You seem to have a missing or misspelled \string\end#1 ...}%
  \let\envir@\empty\fi}
\newif\iftemp@
\def\notempty#1{TT\fi\def\test@{#1}\ifx\test@\empty\temp@false
  \else\temp@true\fi \iftemp@}
\font@\tensmc=cmcsc10
\font@\sevenex=cmex7
\font@\sevenit=cmti7
\font@\eightrm=cmr8 
\font@\sixrm=cmr6 
\font@\eighti=cmmi8     \skewchar\eighti='177 
\font@\sixi=cmmi6       \skewchar\sixi='177   
\font@\eightsy=cmsy8    \skewchar\eightsy='60 
\font@\sixsy=cmsy6      \skewchar\sixsy='60   
\font@\eightex=cmex8
\font@\eightbf=cmbx8 
\font@\sixbf=cmbx6   
\font@\eightit=cmti8 
\font@\eightsl=cmsl8 
\font@\eightsmc=cmcsc8
\font@\eighttt=cmtt8 


\loadmsam
\loadmsbm
\loadeufm
\UseAMSsymbols
\newtoks\tenpoint@
\def\tenpoint{\normalbaselineskip12\p@
 \abovedisplayskip12\p@ plus3\p@ minus9\p@
 \belowdisplayskip\abovedisplayskip
 \abovedisplayshortskip\z@ plus3\p@
 \belowdisplayshortskip7\p@ plus3\p@ minus4\p@
 \textonlyfont@\rm\tenrm \textonlyfont@\it\tenit
 \textonlyfont@\sl\tensl \textonlyfont@\bf\tenbf
 \textonlyfont@\smc\tensmc \textonlyfont@\tt\tentt
 \textonlyfont@\bsmc\tenbsmc
 \ifsyntax@ \def\big##1{{\hbox{$\left##1\right.$}}}%
  \let\Big\big \let\bigg\big \let\Bigg\big
 \else
  \textfont\z@=\tenrm  \scriptfont\z@=\sevenrm  \scriptscriptfont\z@=\fiverm
  \textfont\@ne=\teni  \scriptfont\@ne=\seveni  \scriptscriptfont\@ne=\fivei
  \textfont\tw@=\tensy \scriptfont\tw@=\sevensy \scriptscriptfont\tw@=\fivesy
  \textfont\thr@@=\tenex \scriptfont\thr@@=\sevenex
        \scriptscriptfont\thr@@=\sevenex
  \textfont\itfam=\tenit \scriptfont\itfam=\sevenit
        \scriptscriptfont\itfam=\sevenit
  \textfont\bffam=\tenbf \scriptfont\bffam=\sevenbf
        \scriptscriptfont\bffam=\fivebf
  \setbox\strutbox\hbox{\vrule height8.5\p@ depth3.5\p@ width\z@}%
  \setbox\strutbox@\hbox{\lower.5\normallineskiplimit\vbox{%
        \kern-\normallineskiplimit\copy\strutbox}}%
 \setbox\z@\vbox{\hbox{$($}\kern\z@}\bigsize@=1.2\ht\z@
 \fi
 \normalbaselines\rm\ex@.2326ex\jot3\ex@\the\tenpoint@}
\newtoks\eightpoint@
\def\eightpoint{\normalbaselineskip10\p@
 \abovedisplayskip10\p@ plus2.4\p@ minus7.2\p@
 \belowdisplayskip\abovedisplayskip
 \abovedisplayshortskip\z@ plus2.4\p@
 \belowdisplayshortskip5.6\p@ plus2.4\p@ minus3.2\p@
 \textonlyfont@\rm\eightrm \textonlyfont@\it\eightit
 \textonlyfont@\sl\eightsl \textonlyfont@\bf\eightbf
 \textonlyfont@\smc\eightsmc \textonlyfont@\tt\eighttt
 \textonlyfont@\bsmc\eightbsmc
 \ifsyntax@\def\big##1{{\hbox{$\left##1\right.$}}}%
  \let\Big\big \let\bigg\big \let\Bigg\big
 \else
  \textfont\z@=\eightrm \scriptfont\z@=\sixrm \scriptscriptfont\z@=\fiverm
  \textfont\@ne=\eighti \scriptfont\@ne=\sixi \scriptscriptfont\@ne=\fivei
  \textfont\tw@=\eightsy \scriptfont\tw@=\sixsy \scriptscriptfont\tw@=\fivesy
  \textfont\thr@@=\eightex \scriptfont\thr@@=\sevenex
   \scriptscriptfont\thr@@=\sevenex
  \textfont\itfam=\eightit \scriptfont\itfam=\sevenit
   \scriptscriptfont\itfam=\sevenit
  \textfont\bffam=\eightbf \scriptfont\bffam=\sixbf
   \scriptscriptfont\bffam=\fivebf
 \setbox\strutbox\hbox{\vrule height7\p@ depth3\p@ width\z@}%
 \setbox\strutbox@\hbox{\raise.5\normallineskiplimit\vbox{%
   \kern-\normallineskiplimit\copy\strutbox}}%
 \setbox\z@\vbox{\hbox{$($}\kern\z@}\bigsize@=1.2\ht\z@
 \fi
 \normalbaselines\eightrm\ex@.2326ex\jot3\ex@\the\eightpoint@}

\font@\twelverm=cmr10 scaled\magstep1
\font@\twelveit=cmti10 scaled\magstep1
\font@\twelvesl=cmsl10 scaled\magstep1
\font@\twelvesmc=cmcsc10 scaled\magstep1
\font@\twelvett=cmtt10 scaled\magstep1
\font@\twelvebf=cmbx10 scaled\magstep1
\font@\twelvei=cmmi10 scaled\magstep1
\font@\twelvesy=cmsy10 scaled\magstep1
\font@\twelveex=cmex10 scaled\magstep1
\font@\twelvemsa=msam10 scaled\magstep1
\font@\twelveeufm=eufm10 scaled\magstep1
\font@\twelvemsb=msbm10 scaled\magstep1
\newtoks\twelvepoint@
\def\twelvepoint{\normalbaselineskip15\p@
 \abovedisplayskip15\p@ plus3.6\p@ minus10.8\p@
 \belowdisplayskip\abovedisplayskip
 \abovedisplayshortskip\z@ plus3.6\p@
 \belowdisplayshortskip8.4\p@ plus3.6\p@ minus4.8\p@
 \textonlyfont@\rm\twelverm \textonlyfont@\it\twelveit
 \textonlyfont@\sl\twelvesl \textonlyfont@\bf\twelvebf
 \textonlyfont@\smc\twelvesmc \textonlyfont@\tt\twelvett
 \textonlyfont@\bsmc\twelvebsmc
 \ifsyntax@ \def\big##1{{\hbox{$\left##1\right.$}}}%
  \let\Big\big \let\bigg\big \let\Bigg\big
 \else
  \textfont\z@=\twelverm  \scriptfont\z@=\tenrm  \scriptscriptfont\z@=\sevenrm
  \textfont\@ne=\twelvei  \scriptfont\@ne=\teni  \scriptscriptfont\@ne=\seveni
  \textfont\tw@=\twelvesy \scriptfont\tw@=\tensy \scriptscriptfont\tw@=\sevensy
  \textfont\thr@@=\twelveex \scriptfont\thr@@=\tenex
        \scriptscriptfont\thr@@=\tenex
  \textfont\itfam=\twelveit \scriptfont\itfam=\tenit
        \scriptscriptfont\itfam=\tenit
  \textfont\bffam=\twelvebf \scriptfont\bffam=\tenbf
        \scriptscriptfont\bffam=\sevenbf
  \setbox\strutbox\hbox{\vrule height10.2\p@ depth4.2\p@ width\z@}%
  \setbox\strutbox@\hbox{\lower.6\normallineskiplimit\vbox{%
        \kern-\normallineskiplimit\copy\strutbox}}%
 \setbox\z@\vbox{\hbox{$($}\kern\z@}\bigsize@=1.4\ht\z@
 \fi
 \normalbaselines\rm\ex@.2326ex\jot3.6\ex@\the\twelvepoint@}

\def\headfonts{\twelvepoint\bf}

\font@\fourteenrm=cmr10 scaled\magstep2
\font@\fourteenit=cmti10 scaled\magstep2
\font@\fourteensl=cmsl10 scaled\magstep2
\font@\fourteensmc=cmcsc10 scaled\magstep2
\font@\fourteentt=cmtt10 scaled\magstep2
\font@\fourteenbf=cmbx10 scaled\magstep2
\font@\fourteeni=cmmi10 scaled\magstep2
\font@\fourteensy=cmsy10 scaled\magstep2
\font@\fourteenex=cmex10 scaled\magstep2
\font@\fourteenmsa=msam10 scaled\magstep2
\font@\fourteeneufm=eufm10 scaled\magstep2
\font@\fourteenmsb=msbm10 scaled\magstep2
\newtoks\fourteenpoint@
\def\fourteenpoint{\normalbaselineskip15\p@
 \abovedisplayskip18\p@ plus4.3\p@ minus12.9\p@
 \belowdisplayskip\abovedisplayskip
 \abovedisplayshortskip\z@ plus4.3\p@
 \belowdisplayshortskip10.1\p@ plus4.3\p@ minus5.8\p@
 \textonlyfont@\rm\fourteenrm \textonlyfont@\it\fourteenit
 \textonlyfont@\sl\fourteensl \textonlyfont@\bf\fourteenbf
 \textonlyfont@\smc\fourteensmc \textonlyfont@\tt\fourteentt
 \textonlyfont@\bsmc\fourteenbsmc
 \ifsyntax@ \def\big##1{{\hbox{$\left##1\right.$}}}%
  \let\Big\big \let\bigg\big \let\Bigg\big
 \else
  \textfont\z@=\fourteenrm  \scriptfont\z@=\twelverm  \scriptscriptfont\z@=\tenrm
  \textfont\@ne=\fourteeni  \scriptfont\@ne=\twelvei  \scriptscriptfont\@ne=\teni
  \textfont\tw@=\fourteensy \scriptfont\tw@=\twelvesy \scriptscriptfont\tw@=\tensy
  \textfont\thr@@=\fourteenex \scriptfont\thr@@=\twelveex
        \scriptscriptfont\thr@@=\twelveex
  \textfont\itfam=\fourteenit \scriptfont\itfam=\twelveit
        \scriptscriptfont\itfam=\twelveit
  \textfont\bffam=\fourteenbf \scriptfont\bffam=\twelvebf
        \scriptscriptfont\bffam=\tenbf
  \setbox\strutbox\hbox{\vrule height12.2\p@ depth5\p@ width\z@}%
  \setbox\strutbox@\hbox{\lower.72\normallineskiplimit\vbox{%
        \kern-\normallineskiplimit\copy\strutbox}}%
 \setbox\z@\vbox{\hbox{$($}\kern\z@}\bigsize@=1.7\ht\z@
 \fi
 \normalbaselines\rm\ex@.2326ex\jot4.3\ex@\the\fourteenpoint@}

\def\chapheadfonts{\fourteenpoint\bf}

\font@\seventeenrm=cmr10 scaled\magstep3
\font@\seventeenit=cmti10 scaled\magstep3
\font@\seventeensl=cmsl10 scaled\magstep3
\font@\seventeensmc=cmcsc10 scaled\magstep3
\font@\seventeentt=cmtt10 scaled\magstep3
\font@\seventeenbf=cmbx10 scaled\magstep3
\font@\seventeeni=cmmi10 scaled\magstep3
\font@\seventeensy=cmsy10 scaled\magstep3
\font@\seventeenex=cmex10 scaled\magstep3
\font@\seventeenmsa=msam10 scaled\magstep3
\font@\seventeeneufm=eufm10 scaled\magstep3
\font@\seventeenmsb=msbm10 scaled\magstep3
\newtoks\seventeenpoint@
\def\seventeenpoint{\normalbaselineskip18\p@
 \abovedisplayskip21.6\p@ plus5.2\p@ minus15.4\p@
 \belowdisplayskip\abovedisplayskip
 \abovedisplayshortskip\z@ plus5.2\p@
 \belowdisplayshortskip12.1\p@ plus5.2\p@ minus7\p@
 \textonlyfont@\rm\seventeenrm \textonlyfont@\it\seventeenit
 \textonlyfont@\sl\seventeensl \textonlyfont@\bf\seventeenbf
 \textonlyfont@\smc\seventeensmc \textonlyfont@\tt\seventeentt
 \textonlyfont@\bsmc\seventeenbsmc
 \ifsyntax@ \def\big##1{{\hbox{$\left##1\right.$}}}%
  \let\Big\big \let\bigg\big \let\Bigg\big
 \else
  \textfont\z@=\seventeenrm  \scriptfont\z@=\fourteenrm  \scriptscriptfont\z@=\twelverm
  \textfont\@ne=\seventeeni  \scriptfont\@ne=\fourteeni  \scriptscriptfont\@ne=\twelvei
  \textfont\tw@=\seventeensy \scriptfont\tw@=\fourteensy \scriptscriptfont\tw@=\twelvesy
  \textfont\thr@@=\seventeenex \scriptfont\thr@@=\fourteenex
        \scriptscriptfont\thr@@=\fourteenex
  \textfont\itfam=\seventeenit \scriptfont\itfam=\fourteenit
        \scriptscriptfont\itfam=\fourteenit
  \textfont\bffam=\seventeenbf \scriptfont\bffam=\fourteenbf
        \scriptscriptfont\bffam=\twelvebf
  \setbox\strutbox\hbox{\vrule height14.6\p@ depth6\p@ width\z@}%
  \setbox\strutbox@\hbox{\lower.86\normallineskiplimit\vbox{%
        \kern-\normallineskiplimit\copy\strutbox}}%
 \setbox\z@\vbox{\hbox{$($}\kern\z@}\bigsize@=2\ht\z@
 \fi
 \normalbaselines\rm\ex@.2326ex\jot5.2\ex@\the\seventeenpoint@}

\font@\rrrrrm=cmr10 scaled\magstep4
\font@\bigtitlefont=cmbx10 scaled\magstep4

\parindent1pc
\normallineskiplimit\p@
\newdimen\indenti \indenti=2pc
\def\pageheight#1{\vsize#1}
\def\pagewidth#1{\hsize#1%
   \captionwidth@\hsize \advance\captionwidth@-2\indenti}
\pagewidth{30pc} \pageheight{47pc}
\def\topmatter{%
 \ifx\undefined\msafam
 \else\font@\eightmsa=msam8 \font@\sixmsa=msam6
   \ifsyntax@\else \addto\tenpoint{\textfont\msafam=\tenmsa
              \scriptfont\msafam=\sevenmsa \scriptscriptfont\msafam=\fivemsa}%
     \addto\eightpoint{\textfont\msafam=\eightmsa \scriptfont\msafam=\sixmsa
              \scriptscriptfont\msafam=\fivemsa}%
   \fi
 \fi
 \ifx\undefined\msbfam
 \else\font@\eightmsb=msbm8 \font@\sixmsb=msbm6
   \ifsyntax@\else \addto\tenpoint{\textfont\msbfam=\tenmsb
         \scriptfont\msbfam=\sevenmsb \scriptscriptfont\msbfam=\fivemsb}%
     \addto\eightpoint{\textfont\msbfam=\eightmsb \scriptfont\msbfam=\sixmsb
         \scriptscriptfont\msbfam=\fivemsb}%
   \fi
 \fi
 \ifx\undefined\eufmfam
 \else \font@\eighteufm=eufm8 \font@\sixeufm=eufm6
   \ifsyntax@\else \addto\tenpoint{\textfont\eufmfam=\teneufm
       \scriptfont\eufmfam=\seveneufm \scriptscriptfont\eufmfam=\fiveeufm}%
     \addto\eightpoint{\textfont\eufmfam=\eighteufm
       \scriptfont\eufmfam=\sixeufm \scriptscriptfont\eufmfam=\fiveeufm}%
   \fi
 \fi
 \ifx\undefined\eufbfam
 \else \font@\eighteufb=eufb8 \font@\sixeufb=eufb6
   \ifsyntax@\else \addto\tenpoint{\textfont\eufbfam=\teneufb
      \scriptfont\eufbfam=\seveneufb \scriptscriptfont\eufbfam=\fiveeufb}%
    \addto\eightpoint{\textfont\eufbfam=\eighteufb
      \scriptfont\eufbfam=\sixeufb \scriptscriptfont\eufbfam=\fiveeufb}%
   \fi
 \fi
 \ifx\undefined\eusmfam
 \else \font@\eighteusm=eusm8 \font@\sixeusm=eusm6
   \ifsyntax@\else \addto\tenpoint{\textfont\eusmfam=\teneusm
       \scriptfont\eusmfam=\seveneusm \scriptscriptfont\eusmfam=\fiveeusm}%
     \addto\eightpoint{\textfont\eusmfam=\eighteusm
       \scriptfont\eusmfam=\sixeusm \scriptscriptfont\eusmfam=\fiveeusm}%
   \fi
 \fi
 \ifx\undefined\eusbfam
 \else \font@\eighteusb=eusb8 \font@\sixeusb=eusb6
   \ifsyntax@\else \addto\tenpoint{\textfont\eusbfam=\teneusb
       \scriptfont\eusbfam=\seveneusb \scriptscriptfont\eusbfam=\fiveeusb}%
     \addto\eightpoint{\textfont\eusbfam=\eighteusb
       \scriptfont\eusbfam=\sixeusb \scriptscriptfont\eusbfam=\fiveeusb}%
   \fi
 \fi
 \ifx\undefined\eurmfam
 \else \font@\eighteurm=eurm8 \font@\sixeurm=eurm6
   \ifsyntax@\else \addto\tenpoint{\textfont\eurmfam=\teneurm
       \scriptfont\eurmfam=\seveneurm \scriptscriptfont\eurmfam=\fiveeurm}%
     \addto\eightpoint{\textfont\eurmfam=\eighteurm
       \scriptfont\eurmfam=\sixeurm \scriptscriptfont\eurmfam=\fiveeurm}%
   \fi
 \fi
 \ifx\undefined\eurbfam
 \else \font@\eighteurb=eurb8 \font@\sixeurb=eurb6
   \ifsyntax@\else \addto\tenpoint{\textfont\eurbfam=\teneurb
       \scriptfont\eurbfam=\seveneurb \scriptscriptfont\eurbfam=\fiveeurb}%
    \addto\eightpoint{\textfont\eurbfam=\eighteurb
       \scriptfont\eurbfam=\sixeurb \scriptscriptfont\eurbfam=\fiveeurb}%
   \fi
 \fi
 \ifx\undefined\cmmibfam
 \else \font@\eightcmmib=cmmib8 \font@\sixcmmib=cmmib6
   \ifsyntax@\else \addto\tenpoint{\textfont\cmmibfam=\tencmmib
       \scriptfont\cmmibfam=\sevencmmib \scriptscriptfont\cmmibfam=\fivecmmib}%
    \addto\eightpoint{\textfont\cmmibfam=\eightcmmib
       \scriptfont\cmmibfam=\sixcmmib \scriptscriptfont\cmmibfam=\fivecmmib}%
   \fi
 \fi
 \ifx\undefined\cmbsyfam
 \else \font@\eightcmbsy=cmbsy8 \font@\sixcmbsy=cmbsy6
   \ifsyntax@\else \addto\tenpoint{\textfont\cmbsyfam=\tencmbsy
      \scriptfont\cmbsyfam=\sevencmbsy \scriptscriptfont\cmbsyfam=\fivecmbsy}%
    \addto\eightpoint{\textfont\cmbsyfam=\eightcmbsy
      \scriptfont\cmbsyfam=\sixcmbsy \scriptscriptfont\cmbsyfam=\fivecmbsy}%
   \fi
 \fi
 \let\topmatter\relax}
\def\chapterno@{\uppercase\expandafter{\romannumeral\chaptercount@}}
\newcount\chaptercount@
\def\chapter{\nofrills@{\afterassignment\chapterno@
                        CHAPTER \global\chaptercount@=}\chapter@
 \DNii@##1{\leavevmode\hskip-\leftskip
   \rlap{\vbox to\z@{\vss\centerline{\eightpoint
   \chapter@##1\unskip}\baselineskip2pc\null}}\hskip\leftskip
   \nofrills@false}%
 \FN@\next@}
\newbox\titlebox@

\def\title{\nofrills@{\relax}\title@%
 \DNii@##1\endtitle{\global\setbox\titlebox@\vtop{\tenpoint\bf
 \raggedcenter@\ignorespaces
 \baselineskip1.3\baselineskip\title@{##1}\endgraf}%
 \ifmonograph@ \edef\next{\the\leftheadtoks}\ifx\next\empty
    \leftheadtext{##1}\fi
 \fi
 \edef\next{\the\rightheadtoks}\ifx\next\empty \rightheadtext{##1}\fi
 }\FN@\next@}
\newbox\authorbox@
\def\author#1\endauthor{\global\setbox\authorbox@
 \vbox{\tenpoint\smc\raggedcenter@\ignorespaces
 #1\endgraf}\relaxnext@ \edef\next{\the\leftheadtoks}%
 \ifx\next\empty\leftheadtext{#1}\fi}
\newbox\affilbox@
\def\affil#1\endaffil{\global\setbox\affilbox@
 \vbox{\tenpoint\raggedcenter@\ignorespaces#1\endgraf}}
\newcount\addresscount@
\addresscount@\z@
\def\address#1\endaddress{\global\advance\addresscount@\@ne
  \expandafter\gdef\csname address\number\addresscount@\endcsname
  {\vskip12\p@ minus6\p@\noindent\eightpoint\smc\ignorespaces#1\par}}
\def\email{\nofrills@{\eightpoint{\it E-mail\/}:\enspace}\email@
  \DNii@##1\endemail{%
  \expandafter\gdef\csname email\number\addresscount@\endcsname
  {\def\usualspace{{\it\enspace}}\smallskip\noindent\eightpoint\email@
  \ignorespaces##1\par}}%
 \FN@\next@}
\def\thedate@{}
\def\date#1\enddate{\gdef\thedate@{\tenpoint\ignorespaces#1\unskip}}
\def\thethanks@{}
\def\thanks#1\endthanks{\gdef\thethanks@{\eightpoint\ignorespaces#1.\unskip}}
\def\thekeywords@{}
\def\keywords{\nofrills@{{\it Key words and phrases.\enspace}}\keywords@
 \DNii@##1\endkeywords{\def\thekeywords@{\def\usualspace{{\it\enspace}}%
 \eightpoint\keywords@\ignorespaces##1\unskip.}}%
 \FN@\next@}
\def\thesubjclass@{}
\def\subjclass{\nofrills@{{\rm2000 {\it Mathematics Subject
   Classification\/}.\enspace}}\subjclass@
 \DNii@##1\endsubjclass{\def\thesubjclass@{\def\usualspace
  {{\rm\enspace}}\eightpoint\subjclass@\ignorespaces##1\unskip.}}%
 \FN@\next@}
\newbox\abstractbox@
\def\abstract{\nofrills@{{\smc Abstract.\enspace}}\abstract@
 \DNii@{\setbox\abstractbox@\vbox\bgroup\noindent$$\vbox\bgroup
  \def\envir@{abstract}\advance\hsize-2\indenti
  \usualspace@{{\enspace}}\eightpoint \noindent\abstract@\ignorespaces}%
 \FN@\next@}
\def\endabstract{\par\unskip\egroup$$\egroup}
\def\widestnumber#1#2{\begingroup\let\head\null\let\subhead\empty
   \let\subsubhead\subhead
   \ifx#1\head\global\setbox\tocheadbox@\hbox{#2.\enspace}%
   \else\ifx#1\subhead\global\setbox\tocsubheadbox@\hbox{#2.\enspace}%
   \else\ifx#1\key\bgroup\let\endrefitem@\egroup
        \key#2\endrefitem@\global\refindentwd\wd\keybox@
   \else\ifx#1\no\bgroup\let\endrefitem@\egroup
        \no#2\endrefitem@\global\refindentwd\wd\nobox@
   \else\ifx#1\page\global\setbox\pagesbox@\hbox{\quad\bf#2}%
   \else\ifx#1\item\setboxz@h{#2}\global\rosteritemwd\wdz@
        \global\advance\rosteritemwd by.5\parindent
   \else\message{\string\widestnumber is not defined for this option
   (\string#1)}%
\fi\fi\fi\fi\fi\fi\endgroup}
\newif\ifmonograph@
\def\Monograph{\monograph@true \let\headmark\rightheadtext
  \let\varindent@\indent \def\headfont@{\bf}\def\proclaimheadfont@{\smc}%
  \def\demofont@{\smc}}
\let\varindent@\indent

\newbox\tocheadbox@    \newbox\tocsubheadbox@
\newbox\tocbox@
\def\toc{\toc@{Contents}}
\def\newtocdefs{%
   \def \title##1\endtitle
       {\penaltyandskip@\z@\smallskipamount
        \hangindent\wd\tocheadbox@\noindent{\bf##1}}%
   \def \chapter##1{%
        Chapter \uppercase\expandafter{\romannumeral##1.\unskip}\enspace}%
   \def \specialhead##1\endspecialhead
       {\par\hangindent\wd\tocheadbox@ \noindent##1\par}%
   \def \head##1 ##2\endhead
       {\par\hangindent\wd\tocheadbox@ \noindent
        \if\notempty{##1}\hbox to\wd\tocheadbox@{\hfil##1\enspace}\fi
        ##2\par}%
   \def \subhead##1 ##2\endsubhead
       {\par\vskip-\parskip {\normalbaselines
        \advance\leftskip\wd\tocheadbox@
        \hangindent\wd\tocsubheadbox@ \noindent
        \if\notempty{##1}\hbox to\wd\tocsubheadbox@{##1\unskip\hfil}\fi
         ##2\par}}%
   \def \subsubhead##1 ##2\endsubsubhead
       {\par\vskip-\parskip {\normalbaselines
        \advance\leftskip\wd\tocheadbox@
        \hangindent\wd\tocsubheadbox@ \noindent
        \if\notempty{##1}\hbox to\wd\tocsubheadbox@{##1\unskip\hfil}\fi
        ##2\par}}}
\def\toc@#1{\relaxnext@
   \def\page##1%
       {\unskip\penalty0\null\hfil
        \rlap{\hbox to\wd\pagesbox@{\quad\hfil##1}}\hfilneg\penalty\@M}%
 \DN@{\ifx\next\nofrills\DN@\nofrills{\nextii@}%
      \else\DN@{\nextii@{{#1}}}\fi
      \next@}%
 \DNii@##1{%
\ifmonograph@\bgroup\else\setbox\tocbox@\vbox\bgroup
   \centerline{\headfont@\ignorespaces##1\unskip}\nobreak
   \vskip\belowheadskip \fi
   \setbox\tocheadbox@\hbox{0.\enspace}%
   \setbox\tocsubheadbox@\hbox{0.0.\enspace}%
   \leftskip\indenti \rightskip\leftskip
   \setbox\pagesbox@\hbox{\bf\quad000}\advance\rightskip\wd\pagesbox@
   \newtocdefs
 }%
 \FN@\next@}
\def\endtoc{\par\egroup}
\let\pretitle\relax
\let\preauthor\relax
\let\preaffil\relax
\let\predate\relax
\let\preabstract\relax
\let\prepaper\relax
\def\dedicatory #1\enddedicatory{\def\preabstract{{\medskip
  \eightpoint\it \raggedcenter@#1\endgraf}}}
\def\thetranslator@{}
\def\translator#1\endtranslator{\def\thetranslator@{\nobreak\medskip
 \line{\eightpoint\hfil Translated by \uppercase{#1}\qquad\qquad}\nobreak}}
\outer\def\endtopmatter{\runaway@{abstract}%
 \edef\next{\the\leftheadtoks}\ifx\next\empty
  \expandafter\leftheadtext\expandafter{\the\rightheadtoks}\fi
 \ifmonograph@\else
   \ifx\thesubjclass@\empty\else \makefootnote@{}{\thesubjclass@}\fi
   \ifx\thekeywords@\empty\else \makefootnote@{}{\thekeywords@}\fi
   \ifx\thethanks@\empty\else \makefootnote@{}{\thethanks@}\fi
 \fi
  \pretitle
  \ifmonograph@ \topskip7pc \else \topskip4pc \fi
  \box\titlebox@
  \topskip10pt
  \preauthor
  \ifvoid\authorbox@\else \vskip2.5pc plus1pc \unvbox\authorbox@\fi
  \preaffil
  \ifvoid\affilbox@\else \vskip1pc plus.5pc \unvbox\affilbox@\fi
  \predate
  \ifx\thedate@\empty\else \vskip1pc plus.5pc \line{\hfil\thedate@\hfil}\fi
  \preabstract
  \ifvoid\abstractbox@\else \vskip1.5pc plus.5pc \unvbox\abstractbox@ \fi
  \ifvoid\tocbox@\else\vskip1.5pc plus.5pc \unvbox\tocbox@\fi
  \prepaper
  \vskip2pc plus1pc
}
\def\document{\let\fontlist@\relax\let\alloclist@\relax
  \tenpoint}

\newskip\aboveheadskip       \aboveheadskip1.8\bigskipamount
\newdimen\belowheadskip      \belowheadskip1.8\medskipamount

\def\headfont@{\smc}
\def\penaltyandskip@#1#2{\relax\ifdim\lastskip<#2\relax\removelastskip
      \ifnum#1=\z@\else\penalty@#1\relax\fi\vskip#2%
  \else\ifnum#1=\z@\else\penalty@#1\relax\fi\fi}
\def\nobreak{\penalty\@M
  \ifvmode\def\penalty@{\let\penalty@\penalty\count@@@}%
  \everypar{\let\penalty@\penalty\everypar{}}\fi}
\let\penalty@\penalty
\def\heading#1\endheading{\head#1\endhead}

\def\specialheadfont@{\bf}
\outer\def\specialhead{\par\penaltyandskip@{-200}\aboveheadskip
  \begingroup\interlinepenalty\@M\rightskip\z@ plus\hsize \let\\\linebreak
  \specialheadfont@\noindent\ignorespaces}
\def\endspecialhead{\par\endgroup\nobreak\vskip\belowheadskip}
\let\headmark\eat@
\newskip\subheadskip       \subheadskip\medskipamount
\def\subheadfont@{\bf}
\outer\def\subhead{\nofrills@{.\enspace}\subhead@
 \DNii@##1\endsubhead{\par\penaltyandskip@{-100}\subheadskip
  \varindent@{\usualspace@{{\subheadfont@\enspace}}%
 \subheadfont@\ignorespaces##1\unskip\subhead@}\ignorespaces}%
 \FN@\next@}
\outer\def\subsubhead{\nofrills@{.\enspace}\subsubhead@
 \DNii@##1\endsubsubhead{\par\penaltyandskip@{-50}\medskipamount
      {\usualspace@{{\it\enspace}}%
  \it\ignorespaces##1\unskip\subsubhead@}\ignorespaces}%
 \FN@\next@}
\def\proclaimheadfont@{\bf}
\outer\def\proclaim{\runaway@{proclaim}\def\envir@{proclaim}%
  \nofrills@{.\enspace}\proclaim@
 \DNii@##1{\penaltyandskip@{-100}\medskipamount\varindent@
   \usualspace@{{\proclaimheadfont@\enspace}}\proclaimheadfont@
   \ignorespaces##1\unskip\proclaim@
  \sl\ignorespaces}%
 \FN@\next@}
\outer\def\endproclaim{\let\envir@\relax\par\rm
  \penaltyandskip@{55}\medskipamount}
\def\demoheadfont@{\it}
\def\demo{\runaway@{proclaim}\nofrills@{.\enspace}\demo@
     \DNii@##1{\par\penaltyandskip@\z@\medskipamount
  {\usualspace@{{\demoheadfont@\enspace}}%
  \varindent@\demoheadfont@\ignorespaces##1\unskip\demo@}\rm
  \ignorespaces}\FN@\next@}
\def\enddemo{\par\medskip}
\def\qed{\ifhmode\unskip\nobreak\fi\quad\ifmmode\square\else$\m@th\square$\fi}
\let\remark\demo
\let\endremark\enddemo
\def\definition{\runaway@{proclaim}%
  \nofrills@{.\demoheadfont@\enspace}\definition@
        \DNii@##1{\penaltyandskip@{-100}\medskipamount
        {\usualspace@{{\demoheadfont@\enspace}}%
        \varindent@\demoheadfont@\ignorespaces##1\unskip\definition@}%
        \rm \ignorespaces}\FN@\next@}


\newdimen\rosteritemwd
\newcount\rostercount@
\newif\iffirstitem@
\let\plainitem@\item
\newtoks\everypartoks@
\def\par@{\everypartoks@\expandafter{\the\everypar}\everypar{}}
\def\roster{\edef\leftskip@{\leftskip\the\leftskip}%
 \relaxnext@
 \rostercount@\z@  
 \def\item{\FN@\rosteritem@}%
 \DN@{\ifx\next\runinitem\let\next@\nextii@\else
  \let\next@\nextiii@\fi\next@}%
 \DNii@\runinitem  
  {\unskip  
   \DN@{\ifx\next[\let\next@\nextii@\else
    \ifx\next"\let\next@\nextiii@\else\let\next@\nextiv@\fi\fi\next@}%
   \DNii@[####1]{\rostercount@####1\relax
    \enspace{\rm(\number\rostercount@)}~\ignorespaces}%
   \def\nextiii@"####1"{\enspace{\rm####1}~\ignorespaces}%
   \def\nextiv@{\enspace{\rm(1)}\rostercount@\@ne~}%
   \par@\firstitem@false  
   \FN@\next@}%
 \def\nextiii@{\par\par@  
  \penalty\@m\smallskip\vskip-\parskip
  \firstitem@true}%
 \FN@\next@}
\def\rosteritem@{\iffirstitem@\firstitem@false\else\par\vskip-\parskip\fi
 \leftskip3\parindent\noindent  
 \DNii@[##1]{\rostercount@##1\relax
  \llap{\hbox to2.5\parindent{\hss\rm(\number\rostercount@)}%
   \hskip.5\parindent}\ignorespaces}%
 \def\nextiii@"##1"{%
  \llap{\hbox to2.5\parindent{\hss\rm##1}\hskip.5\parindent}\ignorespaces}%
 \def\nextiv@{\advance\rostercount@\@ne
  \llap{\hbox to2.5\parindent{\hss\rm(\number\rostercount@)}%
   \hskip.5\parindent}}%
 \ifx\next[\let\next@\nextii@\else\ifx\next"\let\next@\nextiii@\else
  \let\next@\nextiv@\fi\fi\next@}

\newif\ifnextRunin@
\def\endroster{\relaxnext@
 \par\leftskip@  
 \penalty-50 \vskip-\parskip\smallskip  
 \DN@{\ifx\next\Runinitem\let\next@\relax
  \else\nextRunin@false\let\item\plainitem@  
   \ifx\next\par 
    \DN@\par{\everypar\expandafter{\the\everypartoks@}}%
   \else  
    \DN@{\noindent\everypar\expandafter{\the\everypartoks@}}%
  \fi\fi\next@}%
 \FN@\next@}
\newcount\rosterhangafter@
\def\Runinitem#1\roster\runinitem{\relaxnext@
 \rostercount@\z@ 
 \def\item{\FN@\rosteritem@}%
 \def\runinitem@{#1}%
 \DN@{\ifx\next[\let\next\nextii@\else\ifx\next"\let\next\nextiii@
  \else\let\next\nextiv@\fi\fi\next}%
 \DNii@[##1]{\rostercount@##1\relax
  \def\item@{{\rm(\number\rostercount@)}}\nextv@}%
 \def\nextiii@"##1"{\def\item@{{\rm##1}}\nextv@}%
 \def\nextiv@{\advance\rostercount@\@ne
  \def\item@{{\rm(\number\rostercount@)}}\nextv@}%
 \def\nextv@{\setbox\z@\vbox  
  {\ifnextRunin@\noindent\fi  
  \runinitem@\unskip\enspace\item@~\par  
  \global\rosterhangafter@\prevgraf}%
  \firstitem@false  
  \ifnextRunin@\else\par\fi  
  \hangafter\rosterhangafter@\hangindent3\parindent
  \ifnextRunin@\noindent\fi  
  \runinitem@\unskip\enspace 
  \item@~\ifnextRunin@\else\par@\fi  
  \nextRunin@true\ignorespaces}%
 \FN@\next@}
\def\footmarkform@#1{$\m@th^{#1}$}
\let\thefootnotemark\footmarkform@
\def\makefootnote@#1#2{\insert\footins
 {\interlinepenalty\interfootnotelinepenalty
 \eightpoint\splittopskip\ht\strutbox\splitmaxdepth\dp\strutbox
 \floatingpenalty\@MM\leftskip\z@\rightskip\z@\spaceskip\z@\xspaceskip\z@
 \leavevmode{#1}\footstrut\ignorespaces#2\unskip\lower\dp\strutbox
 \vbox to\dp\strutbox{}}}
\newcount\footmarkcount@
\footmarkcount@\z@
\def\footnotemark{\let\@sf\empty\relaxnext@
 \ifhmode\edef\@sf{\spacefactor\the\spacefactor}\/\fi
 \DN@{\ifx[\next\let\next@\nextii@\else
  \ifx"\next\let\next@\nextiii@\else
  \let\next@\nextiv@\fi\fi\next@}%
 \DNii@[##1]{\footmarkform@{##1}\@sf}%
 \def\nextiii@"##1"{{##1}\@sf}%
 \def\nextiv@{\iffirstchoice@\global\advance\footmarkcount@\@ne\fi
  \footmarkform@{\number\footmarkcount@}\@sf}%
 \FN@\next@}
\def\footnotetext{\relaxnext@
 \DN@{\ifx[\next\let\next@\nextii@\else
  \ifx"\next\let\next@\nextiii@\else
  \let\next@\nextiv@\fi\fi\next@}%
 \DNii@[##1]##2{\makefootnote@{\footmarkform@{##1}}{##2}}%
 \def\nextiii@"##1"##2{\makefootnote@{##1}{##2}}%
 \def\nextiv@##1{\makefootnote@{\footmarkform@{\number\footmarkcount@}}{##1}}%
 \FN@\next@}
\def\footnote{\let\@sf\empty\relaxnext@
 \ifhmode\edef\@sf{\spacefactor\the\spacefactor}\/\fi
 \DN@{\ifx[\next\let\next@\nextii@\else
  \ifx"\next\let\next@\nextiii@\else
  \let\next@\nextiv@\fi\fi\next@}%
 \DNii@[##1]##2{\footnotemark[##1]\footnotetext[##1]{##2}}%
 \def\nextiii@"##1"##2{\footnotemark"##1"\footnotetext"##1"{##2}}%
 \def\nextiv@##1{\footnotemark\footnotetext{##1}}%
 \FN@\next@}
\def\adjustfootnotemark#1{\advance\footmarkcount@#1\relax}
\def\footnoterule{\kern-3\p@
  \hrule width 5pc\kern 2.6\p@} 
\def\captionfont@{\smc}
\def\topcaption#1#2\endcaption{%
  {\dimen@\hsize \advance\dimen@-\captionwidth@
   \rm\raggedcenter@ \advance\leftskip.5\dimen@ \rightskip\leftskip
  {\captionfont@#1}%
  \if\notempty{#2}.\enspace\ignorespaces#2\fi
  \endgraf}\nobreak\bigskip}
\def\botcaption#1#2\endcaption{%
  \nobreak\bigskip
  \setboxz@h{\captionfont@#1\if\notempty{#2}.\enspace\rm#2\fi}%
  {\dimen@\hsize \advance\dimen@-\captionwidth@
   \leftskip.5\dimen@ \rightskip\leftskip
   \noindent \ifdim\wdz@>\captionwidth@ 
   \else\hfil\fi 
  {\captionfont@#1}\if\notempty{#2}.\enspace\rm#2\fi\endgraf}}
\def\@ins{\par\begingroup\def\vspace##1{\vskip##1\relax}%
  \def\captionwidth##1{\captionwidth@##1\relax}%
  \setbox\z@\vbox\bgroup} 
\def\block{\RIfMIfI@\nondmatherr@\block\fi
       \else\ifvmode\vskip\abovedisplayskip\noindent\fi
        $$\def\endblock{\par\egroup$$}\fi
  \vbox\bgroup\advance\hsize-2\indenti\noindent}
\def\endblock{\par\egroup}
\def\cite#1{{\rm[{\citefont@\m@th#1}]}}
\def\citefont@{\rm}
\def\refsfont@{\eightpoint}
\outer\def\Refs{\runaway@{proclaim}%
 \relaxnext@ \DN@{\ifx\next\nofrills\DN@\nofrills{\nextii@}\else
  \DN@{\nextii@{References}}\fi\next@}%
 \DNii@##1{\penaltyandskip@{-200}\aboveheadskip
  \line{\hfil\headfont@\ignorespaces##1\unskip\hfil}\nobreak
  \vskip\belowheadskip
  \begingroup\refsfont@\sfcode`.=\@m}%
 \FN@\next@}
\def\endRefs{\par\endgroup}
\newbox\nobox@            \newbox\keybox@           \newbox\bybox@
\newbox\paperbox@         \newbox\paperinfobox@     \newbox\jourbox@
\newbox\volbox@           \newbox\issuebox@         \newbox\yrbox@
\newbox\pagesbox@         \newbox\bookbox@          \newbox\bookinfobox@
\newbox\publbox@          \newbox\publaddrbox@      \newbox\finalinfobox@
\newbox\edsbox@           \newbox\langbox@
\newif\iffirstref@        \newif\iflastref@
\newif\ifprevjour@        \newif\ifbook@            \newif\ifprevinbook@
\newif\ifquotes@          \newif\ifbookquotes@      \newif\ifpaperquotes@
\newdimen\bysamerulewd@
\setboxz@h{\refsfont@\kern3em}
\bysamerulewd@\wdz@
\newdimen\refindentwd
\setboxz@h{\refsfont@ 00. }
\refindentwd\wdz@
\outer\def\ref{\begingroup \noindent\hangindent\refindentwd
 \firstref@true \def\nofrills{\def\refkern@{\kern3sp}}%
 \ref@}
\def\ref@{\book@false \bgroup\let\endrefitem@\egroup \ignorespaces}
\def\moreref{\endrefitem@\endref@\firstref@false\ref@}%
\def\transl{\endrefitem@\endref@\firstref@false
  \book@false
  \prepunct@
  \setboxz@h\bgroup \aftergroup\unhbox\aftergroup\z@
    \def\endrefitem@{\unskip\refkern@\egroup}\ignorespaces}%
\def\emptyifempty@{\dimen@\wd\currbox@
  \advance\dimen@-\wd\z@ \advance\dimen@-.1\p@
  \ifdim\dimen@<\z@ \setbox\currbox@\copy\voidb@x \fi}
\let\refkern@\relax
\def\endrefitem@{\unskip\refkern@\egroup
  \setboxz@h{\refkern@}\emptyifempty@}\ignorespaces
\def\refdef@#1#2#3{\edef\next@{\noexpand\endrefitem@
  \let\noexpand\currbox@\csname\expandafter\eat@\string#1box@\endcsname
    \noexpand\setbox\noexpand\currbox@\hbox\bgroup}%
  \toks@\expandafter{\next@}%
  \if\notempty{#2#3}\toks@\expandafter{\the\toks@
  \def\endrefitem@{\unskip#3\refkern@\egroup
  \setboxz@h{#2#3\refkern@}\emptyifempty@}#2}\fi
  \toks@\expandafter{\the\toks@\ignorespaces}%
  \edef#1{\the\toks@}}
\refdef@\no{}{. }
\refdef@\key{[\m@th}{] }
\refdef@\by{}{}
\def\bysame{\by\hbox to\bysamerulewd@{\hrulefill}\thinspace
   \kern0sp}
\def\manyby{\message{\string\manyby is no longer necessary; \string\by
  can be used instead, starting with version 2.0 of \styname.STY}\by}
\refdef@\paper{\ifpaperquotes@``\fi\it}{}
\refdef@\paperinfo{}{}
\def\jour{\endrefitem@\let\currbox@\jourbox@
  \setbox\currbox@\hbox\bgroup
  \def\endrefitem@{\unskip\refkern@\egroup
    \setboxz@h{\refkern@}\emptyifempty@
    \ifvoid\jourbox@\else\prevjour@true\fi}%
\ignorespaces}
\refdef@\vol{\ifbook@\else\bf\fi}{}
\refdef@\issue{no. }{}
\refdef@\yr{}{}
\refdef@\pages{}{}
\def\page{\endrefitem@\def\pp@{\def\pp@{pp.~}p.~}\let\currbox@\pagesbox@
  \setbox\currbox@\hbox\bgroup\ignorespaces}
\def\pp@{pp.~}
\def\book{\endrefitem@ \let\currbox@\bookbox@
 \setbox\currbox@\hbox\bgroup\def\endrefitem@{\unskip\refkern@\egroup
  \setboxz@h{\ifbookquotes@``\fi}\emptyifempty@
  \ifvoid\bookbox@\else\book@true\fi}%
  \ifbookquotes@``\fi\it\ignorespaces}
\def\inbook{\endrefitem@
  \let\currbox@\bookbox@\setbox\currbox@\hbox\bgroup
  \def\endrefitem@{\unskip\refkern@\egroup
  \setboxz@h{\ifbookquotes@``\fi}\emptyifempty@
  \ifvoid\bookbox@\else\book@true\previnbook@true\fi}%
  \ifbookquotes@``\fi\ignorespaces}
\refdef@\eds{(}{, eds.)}
\def\ed{\endrefitem@\let\currbox@\edsbox@
 \setbox\currbox@\hbox\bgroup
 \def\endrefitem@{\unskip, ed.)\refkern@\egroup
  \setboxz@h{(, ed.)}\emptyifempty@}(\ignorespaces}
\refdef@\bookinfo{}{}
\refdef@\publ{}{}
\refdef@\publaddr{}{}
\refdef@\finalinfo{}{}
\refdef@\lang{(}{)}

\let\refdef@\relax 
\def\ppunbox@#1{\ifvoid#1\else\prepunct@\unhbox#1\fi}
\def\nocomma@#1{\ifvoid#1\else\changepunct@3\prepunct@\unhbox#1\fi}
\def\changepunct@#1{\ifnum\lastkern<3 \unkern\kern#1sp\fi}
\def\prepunct@{\count@\lastkern\unkern
  \ifnum\lastpenalty=0
    \let\penalty@\relax
  \else
    \edef\penalty@{\penalty\the\lastpenalty\relax}%
  \fi
  \unpenalty
  \let\refspace@\ \ifcase\count@,
\or;\or.\or 
  \or\let\refspace@\relax
  \else,\fi
  \ifquotes@''\quotes@false\fi \penalty@ \refspace@
}
\def\transferpenalty@#1{\dimen@\lastkern\unkern
  \ifnum\lastpenalty=0\unpenalty\let\penalty@\relax
  \else\edef\penalty@{\penalty\the\lastpenalty\relax}\unpenalty\fi
  #1\penalty@\kern\dimen@}
\def\endref{\endrefitem@\lastref@true\endref@
  \par\endgroup \prevjour@false \previnbook@false }
\def\endref@{%
\iffirstref@
  \ifvoid\nobox@\ifvoid\keybox@\indent\fi
  \else\hbox to\refindentwd{\hss\unhbox\nobox@}\fi
  \ifvoid\keybox@
  \else\ifdim\wd\keybox@>\refindentwd
         \box\keybox@
       \else\hbox to\refindentwd{\unhbox\keybox@\hfil}\fi\fi
  \kern4sp\ppunbox@\bybox@
\fi 
  \ifvoid\paperbox@
  \else\prepunct@\unhbox\paperbox@
    \ifpaperquotes@\quotes@true\fi\fi
  \ppunbox@\paperinfobox@
  \ifvoid\jourbox@
    \ifprevjour@ \nocomma@\volbox@
      \nocomma@\issuebox@
      \ifvoid\yrbox@\else\changepunct@3\prepunct@(\unhbox\yrbox@
        \transferpenalty@)\fi
      \ppunbox@\pagesbox@
    \fi 
  \else \prepunct@\unhbox\jourbox@
    \nocomma@\volbox@
    \nocomma@\issuebox@
    \ifvoid\yrbox@\else\changepunct@3\prepunct@(\unhbox\yrbox@
      \transferpenalty@)\fi
    \ppunbox@\pagesbox@
  \fi 
  \ifbook@\prepunct@\unhbox\bookbox@ \ifbookquotes@\quotes@true\fi \fi
  \nocomma@\edsbox@
  \ppunbox@\bookinfobox@
  \ifbook@\ifvoid\volbox@\else\prepunct@ vol.~\unhbox\volbox@
  \fi\fi
  \ppunbox@\publbox@ \ppunbox@\publaddrbox@
  \ifbook@ \ppunbox@\yrbox@
    \ifvoid\pagesbox@
    \else\prepunct@\pp@\unhbox\pagesbox@\fi
  \else
    \ifprevinbook@ \ppunbox@\yrbox@
      \ifvoid\pagesbox@\else\prepunct@\pp@\unhbox\pagesbox@\fi
    \fi \fi
  \ppunbox@\finalinfobox@
  \iflastref@
    \ifvoid\langbox@.\ifquotes@''\fi
    \else\changepunct@2\prepunct@\unhbox\langbox@\fi
  \else
    \ifvoid\langbox@\changepunct@1%
    \else\changepunct@3\prepunct@\unhbox\langbox@
      \changepunct@1\fi
  \fi
}
\outer\def
\Refs



\ref\no \AignAB\by M.    Aigner \yr 2001 \paper Catalan and other
numbers: a recurrent theme\inbook Algebraic Combinatorics and Computer
Science\eds H.~Crapo, D.~Senato\publ Springer--Verlag\publaddr
Berlin\pages 347--390\endref 

\ref\no \AmZeAA\by T. Amdeberhan and D. Zeilberger\paper
Determinants through the looking glass\jour
\hbox{\hskip-13pt}\linebreak Adv\. Appl\. Math\.\vol 27\yr 2001\pages 225--230\finalinfo
{\sl Maple} package DODGSON available at\linebreak
{\tt http://www.math.rutgers.edu/\~{}zeilberg/tokhniot/DODGSON}
\endref





\ref\no \BresAO\by D. M. Bressoud \yr 1999 \book Proofs and confirmations 
--- The story of the alternating sign matrix conjecture\publ 
Cambridge University Press\publaddr Cambridge\endref





\ref\no \ComtAA \by L. Comtet\book Advanced Combinatorics\publ D.~Reidel\publaddr
Dordrecht, Holland\yr 1974\endref





\ref\no \EgRRAB\by \"O.  E\u gecio\u glu, T. Redmond and C. Ryavec 
\yr 2008 \paper Almost product evaluation of Hankel determinants\jour 
Electron\. J. Combin\. \vol 15 \rm(1)\pages Article~\#R6, 58~pp
\endref

\ref\no \EgRRAD\by \"O.  E\u gecio\u glu, T. Redmond and C. Ryavec 
\yr 2010 \paper A multilinear operator for 
almost product evaluation of Hankel determinants\jour 
J. Combin\. Theory Ser.~A\vol 117\pages 77--103\endref



\ref\no \FuKrAA\by M.    Fulmek and C. Krattenthaler \yr 1997 \paper
Lattice path proofs for determinant formulas for symplectic and
orthogonal characters\jour J. Combin\. Theory Ser.~A\vol 77\pages
3--50\endref 

\ref\no \FuHaAA\by W.    Fulton and J. Harris \yr 1991 
\book Representation Theory
\publ Springer--Verlag
\publaddr New York\endref

\ref\no \GeViAA\by I. M. Gessel and X. Viennot \yr 1985 \paper Binomial 
determinants, paths, and hook length formulae\jour Adv\. in Math\. 
\vol 58\pages 300--321\endref

\ref\no \GeViAB\by I. M. Gessel and X. Viennot \yr 1989 \paper
Determinants, paths, and plane partitions \paperinfo preprint,
1989\finalinfo available at {\tt
http://www.cs.brandeis.edu/\~{}ira}\endref 


\ref\no \GhKrAA\by S. R. Ghorpade and C. Krattenthaler \yr 2004 \paper
The Hilbert series of Pfaffian rings\inbook Algebra, Arithmetic and
Geometry with Applications\eds C.~Christensen, G.~Sundaram, A.~Sathaye
and C.~Bajaj\publ Springer-Ver\-lag\publaddr New York\pages
337--356\endref 






\ref\no \KratBL\by C.    Krattenthaler \yr 1997 \paper The enumeration
of lattice paths with respect to their number of turns\inbook Advances
in Combinatorial Methods and Applications to Probability and
Statistics\ed N.~Balakrishnan\publ Birkh\"auser\publaddr Boston\pages
29--58\endref 

\ref\no \KratBN\by C. Krattenthaler
\paper Advanced determinant calculus\jour S\'eminaire Lotharingien 
Combin\.\vol 42 \rm(``The Andrews Festschrift")\yr 1999\pages 
Article~B42q, 67~pp\endref

\ref\no \KratBZ\by C. Krattenthaler
\paper Advanced determinant calculus: a complement\yr 2005\jour 
Linear Algebra Appl\.\vol 411\pages 64--166\endref

\ref\no \KratBV\by C.    Krattenthaler \yr 2005 \paper On
multiplicities of points on Schubert varieties in Gra\ss mannians
II\jour J. Algebraic Combin\.\vol 22\pages 273--288\endref 

\ref\no \KratBW\by C.    Krattenthaler \yr 2006 \paper Watermelon 
configurations with wall interaction: exact and asymptotic results\jour 
J. Physics: Conf\. Series\vol 42\pages 179--212\endref



\ref\no \LindAA\by B.    Lindstr\"om \yr 1973 \paper On the vector 
representations of induced matroids\jour 
Bull\. London Math\. Soc\.\vol 5\pages 85--90\endref









\ref\no \StanBI\by R. P. Stanley \yr 1999 \book Enumerative
Combinatorics\bookinfo Vol.~2\publ Cambridge University Press\publaddr
Cambridge\endref 

\ref\no \StemAE\by J. R. Stembridge \yr 1990 \paper Nonintersecting paths, pfaffians 
and plane partitions\jour Adv\. in Math\.\vol 83\pages 96--131\endref


\ref\no \VienAE\by X.    Viennot \yr 1983 \book Une th\'eorie combinatoire 
des polyn\^omes orthogonaux g\'en\'eraux\publ UQAM\publaddr Montr\'eal, 
Qu\'e\-bec\endref


\ref\no \ZeilZZ\by D. Zeilberger\paper
The holonomic Ansatz II. Automatic discovery\ (!)\ and proof\ (!!)\ of 
holonomic\linebreak 
determinant evaluations\jour Ann\. Combin\.\vol 11\yr 2007\pages 241--247
\finalinfo {\sl Maple} package DET available at\linebreak 
{\tt http://www.math.rutgers.edu/\~{}zeilberg/tokhniot/DET}
\endref

\endRefs
\enddocument{%
 \runaway@{proclaim}%
\ifmonograph@ 
\else
 \nobreak
 \thetranslator@
 \count@\z@ \loop\ifnum\count@<\addresscount@\advance\count@\@ne
 \csname address\number\count@\endcsname
 \csname email\number\count@\endcsname
 \repeat
\fi
 \vfill\supereject\end}

\def\headfont@{\headfonts}
\def\proclaimheadfont@{\bf}
\def\specialheadfont@{\bf}
\def\subheadfont@{\bf}
\def\demoheadfont@{\smc}

\newif\ifThisToToc \ThisToTocfalse
\newif\iftocloaded \tocloadedfalse

\def\C@L{\noexpand\Cal}\def\B@B{\noexpand\Bbb}\def\fR@K{\noexpand\frak}
\def\S@{\noexpand\S}\def\P@P{\noexpand\"}
\def\xpar{\\}

\def\writetoc#1{\iftocloaded\ifThisToToc\begingroup\def\totoc{}
  \def\Cal{\noexpand\C@L}\def\Bbb{\noexpand\B@B}
  \def\frak{\noexpand\fR@K}\def\goth{\frak}\def\S{\noexpand\S@}
  \def\"{\noexpand\P@P}
  \def\xpar{\par\penalty100000 }\def\idx##1{##1}\def\\{\xpar}
  \edef\next@{\write\toc{\noindent#1\leaderfill\noexpand\folio\par}}%
  \next@\endgroup\global\ThisToTocfalse\fi\fi}
\def\leaderfill{\leaders\hbox to 1em{\hss.\hss}\hfill}

\newif\ifindexloaded \indexloadedfalse
\def\idx#1{\ifindexloaded\begingroup\def\ign{}\def\it{}\def\/{}%
 \def\smc{}\def\bf{}\def\tt{}%
 \def\Cal{\noexpand\C@L}\def\Bbb{\noexpand\B@B}%
 \def\frak{\noexpand\fR@K}\def\goth{\frak}\def\S{\noexpand\S@}%
  \def\"{\noexpand\P@P}%
 {\edef\next@{\write\index{#1, \noexpand\folio}}\next@}%
 \endgroup\fi{#1}}
\def\ign#1{}

\def\totoc{\global\ThisToToctrue}

\outer\def\head#1\endhead{\par\penaltyandskip@{-200}\aboveheadskip
 {\headfont@\raggedcenter@\interlinepenalty\@M
 \ignorespaces#1\endgraf}\nobreak
 \vskip\belowheadskip
 \headmark{#1}\writetoc{#1}}

\outer\def\chaphead#1\endchaphead{\par\penaltyandskip@{-200}\aboveheadskip
 {\chapheadfonts\raggedcenter@\interlinepenalty\@M
 \ignorespaces#1\endgraf}\nobreak
 \vskip3\belowheadskip
 \headmark{#1}\writetoc{#1}}

\def\folio{{\foliofont@\ifnum\pageno<\z@ \romannumeral-\pageno
 \else\number\pageno \fi}}
\newtoks\leftheadtoks
\newtoks\rightheadtoks

\def\leftheadtext{\nofrills@{\relax}\lht@
  \DNii@##1{\leftheadtoks\expandafter{\lht@{##1}}%
    \mark{\the\leftheadtoks\noexpand\else\the\rightheadtoks}
    \ifsyntax@\setboxz@h{\def\\{\unskip\space\ignorespaces}%
        \headlinefont@##1}\fi}%
  \FN@\next@}
\def\rightheadtext{\nofrills@{\relax}\rht@
  \DNii@##1{\rightheadtoks\expandafter{\rht@{##1}}%
    \mark{\the\leftheadtoks\noexpand\else\the\rightheadtoks}%
    \ifsyntax@\setboxz@h{\def\\{\unskip\space\ignorespaces}%
        \headlinefont@##1}\fi}%
  \FN@\next@}
\def\NoRunningHeads{\global\runheads@false\global\let\headmark\eat@}
\def\NoPageNumbers{\gdef\folio{}}
\newif\iffirstpage@     \firstpage@true
\newif\ifrunheads@      \runheads@true

\newdimen\fullhsize \fullhsize=\hsize
\newdimen\fullvsize \fullvsize=\vsize
\def\fullline{\hbox to\fullhsize}

\def\pagenumbers{\gdef\folio{\folio@}}

\let\norunningheads\NoRunningHeads
\def\userunningheads{\global\runheads@true}
\norunningheads

\headline={\def\chapter#1{\chapterno@. }%
  \def\\{\unskip\space\ignorespaces}\ifrunheads@\headlinefont@
    \ifodd\pageno\rightheadline \else\leftheadline\fi
   \else\hfil\fi\ifNoRunHeadline\global\NoRunHeadlinefalse\fi}
\let\folio@\folio
\def\foliofont@{\foliofont}
\def\foliofont{\eightrm}
\def\headlinefont@{\headlinefont}
\def\headlinefont{\eightpoint\smc}
\def\leftheadline{\rlap{\folio}\hfill
   \ifNoRunHeadline\else\iftrue\topmark\fi\fi \hfill}
\def\rightheadline{\hfill\ifNoRunHeadline
   \else \expandafter\fi
  \hfill \llap{\folio}}
\footline={{\eightpoint\bottremark}%
   \ifrunheads@\else\hfil{\let\foliofont\tenrm\folio}\fi\hfil}
\def\bottremark{}
 
\newif\ifNoRunHeadline      
\def\norunninghead{\global\NoRunHeadlinetrue}
\norunninghead

\output={\output@}
%
\newif\ifoffset\offsetfalse
\output={\output@}
\def\output@{%
 \ifoffset 
  \ifodd\count0\advance\hoffset by0.5truecm
   \else\advance\hoffset by-0.5truecm\fi\fi
 \shipout\vbox{%
  \makeheadline \pagebody \makefootline }%
 \advancepageno \ifnum\outputpenalty>-\@MM\else\dosupereject\fi}

\def\indexoutput#1{%
  \ifoffset 
   \ifodd\count0\advance\hoffset by0.5truecm
    \else\advance\hoffset by-0.5truecm\fi\fi
  \shipout\vbox{\makeheadline
  \vbox to\fullvsize{\boxmaxdepth\maxdepth%
  \ifvoid\topins\else\unvbox\topins\fi%
  #1 %
  \ifvoid\footins\else 
    \vskip\skip\footins
    \footnoterule
    \unvbox\footins\fi
  \ifr@ggedbottom \kern-\dimen@ \vfil \fi}%
  \baselineskip2pc
  \makefootline}%
 \global\advance\pageno\@ne
 \ifnum\outputpenalty>-\@MM\else\dosupereject\fi}
 
 \newbox\partialpage \newdimen\halfsize \halfsize=0.5\fullhsize
 \advance\halfsize by-0.5em

 \def\begindoublecolumns{\output={\indexoutput{\unvbox255}}%
   \begingroup \def\line{\fullline}
   \output={\global\setbox\partialpage=\vbox{\unvbox255\bigskip}}\eject
   \output={\doublecolumnout}\hsize=\halfsize \vsize=2\fullvsize}
 \def\enddoublecolumns{\output={\balancecolumns}\eject
  \endgroup \pagegoal=\fullvsize%
  \output={\output@}}
\def\doublecolumnout{\splittopskip=\topskip \splitmaxdepth=\maxdepth
  \dimen@=\fullvsize \advance\dimen@ by-\ht\partialpage
  \setbox0=\vsplit255 to \dimen@ \setbox2=\vsplit255 to \dimen@
  \indexoutput{\pagesofar} \unvbox255 \penalty\outputpenalty}
\def\pagesofar{\unvbox\partialpage
  \wd0=\hsize \wd2=\hsize \hbox to\fullhsize{\box0\hfil\box2}}
\def\balancecolumns{\setbox0=\vbox{\unvbox255} \dimen@=\ht0
  \advance\dimen@ by\topskip \advance\dimen@ by-\baselineskip
  \divide\dimen@ by2 \splittopskip=\topskip
  {\vbadness=10000 \loop \global\setbox3=\copy0
    \global\setbox1=\vsplit3 to\dimen@
    \ifdim\ht3>\dimen@ \global\advance\dimen@ by1pt \repeat}
  \setbox0=\vbox to\dimen@{\unvbox1} \setbox2=\vbox to\dimen@{\unvbox3}
  \pagesofar}

\def\indexitems{\catcode`\^^M=\active \def\^^M={\par}%
\everypar{\parindent=0pt\hangindent=1em\hangafter=1}%
\noindent}

\tenpoint
\catcode`\@=\active

\def\smallheadings{\let\chapheadfonts\tenpoint\let\headfonts\tenpoint}

\tenpoint
\catcode`\@=\active

\def\LL{\leavevmode\setbox0=\hbox{L}\hbox to\wd0{\hss\char'40L}}

\def\la{\lambda}

\def\si{\sigma}

\def\om{\omega}

\def\P{{\Bbb P}}

\def\today{\ifcase\month\or
 January\or February\or March\or April\or May\or June\or
 July\or August\or September\or October\or November\or December\fi
 \space\number\day, \number\year}

\def\({\left(}
\def\){\right)}
\def\[{\left[}
\def\]{\right]}

\def\sgn{\operatorname{sgn}}

\def\3{\ss}
\catcode`\@=11
\def\dddot#1{\vbox{\ialign{##\crcr
      .\hskip-.5pt.\hskip-.5pt.\crcr\noalign{\kern1.5\p@\nointerlineskip}
      $\hfil\displaystyle{#1}\hfil$\crcr}}}

\newif\iftab@\tab@false
\newif\ifvtab@\vtab@false
\def\tab{\bgroup\tab@true\vtab@false\vst@bfalse\Strich@false%
   \def\\{\global\hline@@false%
     \ifhline@\global\hline@false\global\hline@@true\fi\cr}
   \edef\l@{\the\leftskip}\ialign\bgroup\hskip\l@##\hfil&&##\hfil\cr}
\def\endtab{\cr\egroup\egroup}
\def\vtab{\vtop\bgroup\vst@bfalse\vtab@true\tab@true\Strich@false%
   \bgroup\def\\{\cr}\ialign\bgroup&##\hfil\cr}
\def\endvtab{\cr\egroup\egroup\egroup}
\def\stab{\D@cke0.5pt\null 
 \bgroup\tab@true\vtab@false\vst@bfalse\Strich@true\Let@@\vspace@
 \normalbaselines\offinterlineskip
  \openup\spreadmlines@
 \edef\l@{\the\leftskip}\ialign
 \bgroup\hskip\l@##\hfil&&##\hfil\crcr}
\def\endstab{\crcr\egroup
 \egroup}
\newif\ifvst@b\vst@bfalse
\def\vstab{\D@cke0.5pt\null
 \vtop\bgroup\tab@true\vtab@false\vst@btrue\Strich@true\bgroup\Let@@\vspace@
 \normalbaselines\offinterlineskip
  \openup\spreadmlines@\bgroup}
\def\endvstab{\crcr\egroup\egroup
 \egroup\tab@false\Strich@false}

\newdimen\htstrut@
\htstrut@8.5\p@
\newdimen\htStrut@
\htStrut@12\p@
\newdimen\dpstrut@
\dpstrut@3.5\p@
\newdimen\dpStrut@
\dpStrut@3.5\p@
\def\openup{\afterassignment\@penup\dimen@=}
\def\@penup{\advance\lineskip\dimen@
  \advance\baselineskip\dimen@
  \advance\lineskiplimit\dimen@
  \divide\dimen@ by2
  \advance\htstrut@\dimen@
  \advance\htStrut@\dimen@
  \advance\dpstrut@\dimen@
  \advance\dpStrut@\dimen@}
\def\Let@@{\relax%
    \def\\{\global\hline@@false%
     \ifhline@\global\hline@false\global\hline@@true\fi\cr}%
    \iffalse}\fi}
\def\matrix{\null\,\vcenter\bgroup
 \tab@false\vtab@false\vst@bfalse\Strich@false\Let@@\vspace@
 \normalbaselines\openup\spreadmlines@\ialign
 \bgroup\hfil$\m@th##$\hfil&&\quad\hfil$\m@th##$\hfil\crcr
 \Mathstrut@\crcr\noalign{\kern-\baselineskip}}
\def\endmatrix{\crcr\Mathstrut@\crcr\noalign{\kern-\baselineskip}\egroup
 \egroup\,}
\def\smatrix{\D@cke0.5pt\null\,
 \vcenter\bgroup\tab@false\vtab@false\vst@bfalse\Strich@true\Let@@\vspace@
 \normalbaselines\offinterlineskip
  \openup\spreadmlines@\ialign
 \bgroup\hfil$\m@th##$\hfil&&\quad\hfil$\m@th##$\hfil\crcr}
\def\endsmatrix{\crcr\egroup
 \egroup\,\Strich@false}
\newdimen\D@cke
\def\Dicke#1{\global\D@cke#1}
\newtoks\tabs@\tabs@{&}
\newif\ifStrich@\Strich@false
\newif\iff@rst

\def\Stricherr@{\iftab@\ifvtab@\errmessage{\noexpand\s not allowed
     here. Use \noexpand\vstab!}%
  \else\errmessage{\noexpand\s not allowed here. Use \noexpand\stab!}%
  \fi\else\errmessage{\noexpand\s not allowed
     here. Use \noexpand\smatrix!}\fi}
\def\format{\ifvst@b\else\crcr\fi\egroup\iffalse{\fi\ifnum`}=0 \fi\format@}
\def\format@#1\\{\def\preamble@{#1}%
 \def\Str@chfehlt##1{\ifx##1\s\Stricherr@\fi\ifx##1\\\let\Next\relax%
   \else\let\Next\Str@chfehlt\fi\Next}%
 \def\c{\hfil\noexpand\ifhline@@\hbox{\vrule height\htStrut@%
   depth\dpstrut@ width\z@}\noexpand\fi%
   \ifStrich@\hbox{\vrule height\htstrut@ depth\dpstrut@ width\z@}%
   \fi\iftab@\else$\m@th\fi\the\hashtoks@\iftab@\else$\fi\hfil}%
 \def\r{\hfil\noexpand\ifhline@@\hbox{\vrule height\htStrut@%
   depth\dpstrut@ width\z@}\noexpand\fi%
   \ifStrich@\hbox{\vrule height\htstrut@ depth\dpstrut@ width\z@}%
   \fi\iftab@\else$\m@th\fi\the\hashtoks@\iftab@\else$\fi}%
 \def\l{\noexpand\ifhline@@\hbox{\vrule height\htStrut@%
   depth\dpstrut@ width\z@}\noexpand\fi%
   \ifStrich@\hbox{\vrule height\htstrut@ depth\dpstrut@ width\z@}%
   \fi\iftab@\else$\m@th\fi\the\hashtoks@\iftab@\else$\fi\hfil}%
 \def\s{\ifStrich@\ \the\tabs@\vrule width\D@cke\the\hashtoks@%
          \fi\the\tabs@\ }%
 \def\sa{\ifStrich@\vrule width\D@cke\the\hashtoks@%
            \the\tabs@\ %
            \fi}%
 \def\se{\ifStrich@\ \the\tabs@\vrule width\D@cke\the\hashtoks@\fi}%
 \def\cd{\hfil\noexpand\ifhline@@\hbox{\vrule height\htStrut@%
   depth\dpstrut@ width\z@}\noexpand\fi%
   \ifStrich@\hbox{\vrule height\htstrut@ depth\dpstrut@ width\z@}%
   \fi$\dsize\m@th\the\hashtoks@$\hfil}%
 \def\rd{\hfil\noexpand\ifhline@@\hbox{\vrule height\htStrut@%
   depth\dpstrut@ width\z@}\noexpand\fi%
   \ifStrich@\hbox{\vrule height\htstrut@ depth\dpstrut@ width\z@}%
   \fi$\dsize\m@th\the\hashtoks@$}%
 \def\ld{\noexpand\ifhline@@\hbox{\vrule height\htStrut@%
   depth\dpstrut@ width\z@}\noexpand\fi%
   \ifStrich@\hbox{\vrule height\htstrut@ depth\dpstrut@ width\z@}%
   \fi$\dsize\m@th\the\hashtoks@$\hfil}%
 \ifStrich@\else\Str@chfehlt#1\\\fi%
 \setbox\z@\hbox{\xdef\Preamble@{\preamble@}}\ifnum`{=0 \fi\iffalse}\fi
 \ialign\bgroup\span\Preamble@\crcr}
\newif\ifhline@\hline@false
\newif\ifhline@@\hline@@false
\def\hlinefor#1{\multispan@{\strip@#1 }\leaders\hrule height\D@cke\hfill%
    \global\hline@true\ignorespaces}
\def\Item "#1"{\par\noindent\hangindent2\parindent%
  \hangafter1\setbox0\hbox{\rm#1\enspace}\ifdim\wd0>2\parindent%
  \box0\else\hbox to 2\parindent{\rm#1\hfil}\fi\ignorespaces}
\def\ITEM #1"#2"{\par\noindent\hangafter1\hangindent#1%
  \setbox0\hbox{\rm#2\enspace}\ifdim\wd0>#1%
  \box0\else\hbox to 0pt{\rm#2\hss}\hskip#1\fi\ignorespaces}
\def\item"#1"{\par\noindent\hang%
  \setbox0=\hbox{\rm#1\enspace}\ifdim\wd0>\the\parindent%
  \box0\else\hbox to \parindent{\rm#1\hfil}\enspace\fi\ignorespaces}
\let\plainitem@\item
\catcode`\@=13

\catcode`\@=11
\font\tenln    = line10
\font\tenlnw   = linew10

\newskip\Einheit \Einheit=0.5cm
\newcount\xcoord \newcount\ycoord
\newdimen\xdim \newdimen\ydim \newdimen\PfadD@cke \newdimen\Pfadd@cke

\newcount\@tempcnta
\newcount\@tempcntb

\newdimen\@tempdima
\newdimen\@tempdimb

\newdimen\@wholewidth
\newdimen\@halfwidth

\newcount\@xarg
\newcount\@yarg
\newcount\@yyarg
\newbox\@linechar
\newbox\@tempboxa
\newdimen\@linelen
\newdimen\@clnwd
\newdimen\@clnht

\newif\if@negarg

\def\@whilenoop#1{}
\def\@whiledim#1\do #2{\ifdim #1\relax#2\@iwhiledim{#1\relax#2}\fi}
\def\@iwhiledim#1{\ifdim #1\let\@nextwhile=\@iwhiledim
        \else\let\@nextwhile=\@whilenoop\fi\@nextwhile{#1}}

\def\@whileswnoop#1\fi{}
\def\@whilesw#1\fi#2{#1#2\@iwhilesw{#1#2}\fi\fi}
\def\@iwhilesw#1\fi{#1\let\@nextwhile=\@iwhilesw
         \else\let\@nextwhile=\@whileswnoop\fi\@nextwhile{#1}\fi}

\def\thinlines{\let\@linefnt\tenln \let\@circlefnt\tencirc
  \@wholewidth\fontdimen8\tenln \@halfwidth .5\@wholewidth}
\def\thicklines{\let\@linefnt\tenlnw \let\@circlefnt\tencircw
  \@wholewidth\fontdimen8\tenlnw \@halfwidth .5\@wholewidth}
\thinlines

\PfadD@cke1pt \Pfadd@cke0.5pt
\def\PfadDicke#1{\PfadD@cke#1 \divide\PfadD@cke by2 \Pfadd@cke\PfadD@cke \multiply\PfadD@cke by2}
\long\def\LOOP#1\REPEAT{\def\BODY{#1}\ITERATE}
\def\ITERATE{\BODY \let\next\ITERATE \else\let\next\relax\fi \next}
\let\REPEAT=\fi
\def\Punkt{\hbox{\raise-2pt\hbox to0pt{\hss$\ssize\bullet$\hss}}}
\def\DuennPunkt(#1,#2){\unskip
  \raise#2 \Einheit\hbox to0pt{\hskip#1 \Einheit
          \raise-2.5pt\hbox to0pt{\hss$\bullet$\hss}\hss}}
\def\NormalPunkt(#1,#2){\unskip
  \raise#2 \Einheit\hbox to0pt{\hskip#1 \Einheit
          \raise-3pt\hbox to0pt{\hss\twelvepoint$\bullet$\hss}\hss}}
\def\DickPunkt(#1,#2){\unskip
  \raise#2 \Einheit\hbox to0pt{\hskip#1 \Einheit
          \raise-4pt\hbox to0pt{\hss\fourteenpoint$\bullet$\hss}\hss}}
\def\Kreis(#1,#2){\unskip
  \raise#2 \Einheit\hbox to0pt{\hskip#1 \Einheit
          \raise-4pt\hbox to0pt{\hss\fourteenpoint$\circ$\hss}\hss}}

\def\Line@(#1,#2)#3{\@xarg #1\relax \@yarg #2\relax
\@linelen=#3\Einheit
\ifnum\@xarg =0 \@vline
  \else \ifnum\@yarg =0 \@hline \else \@sline\fi
\fi}

\def\@sline{\ifnum\@xarg< 0 \@negargtrue \@xarg -\@xarg \@yyarg -\@yarg
  \else \@negargfalse \@yyarg \@yarg \fi
\ifnum \@yyarg >0 \@tempcnta\@yyarg \else \@tempcnta -\@yyarg \fi
\ifnum\@tempcnta>6 \@badlinearg\@tempcnta0 \fi
\ifnum\@xarg>6 \@badlinearg\@xarg 1 \fi
\setbox\@linechar\hbox{\@linefnt\@getlinechar(\@xarg,\@yyarg)}%
\ifnum \@yarg >0 \let\@upordown\raise \@clnht\z@
   \else\let\@upordown\lower \@clnht \ht\@linechar\fi
\@clnwd=\wd\@linechar
\if@negarg \hskip -\wd\@linechar \def\@tempa{\hskip -2\wd\@linechar}\else
     \let\@tempa\relax \fi
\@whiledim \@clnwd <\@linelen \do
  {\@upordown\@clnht\copy\@linechar
   \@tempa
   \advance\@clnht \ht\@linechar
   \advance\@clnwd \wd\@linechar}%
\advance\@clnht -\ht\@linechar
\advance\@clnwd -\wd\@linechar
\@tempdima\@linelen\advance\@tempdima -\@clnwd
\@tempdimb\@tempdima\advance\@tempdimb -\wd\@linechar
\if@negarg \hskip -\@tempdimb \else \hskip \@tempdimb \fi
\multiply\@tempdima \@m
\@tempcnta \@tempdima \@tempdima \wd\@linechar \divide\@tempcnta \@tempdima
\@tempdima \ht\@linechar \multiply\@tempdima \@tempcnta
\divide\@tempdima \@m
\advance\@clnht \@tempdima
\ifdim \@linelen <\wd\@linechar
   \hskip \wd\@linechar
  \else\@upordown\@clnht\copy\@linechar\fi}

\def\@hline{\ifnum \@xarg <0 \hskip -\@linelen \fi
\vrule height\Pfadd@cke width \@linelen depth\Pfadd@cke
\ifnum \@xarg <0 \hskip -\@linelen \fi}

\def\@getlinechar(#1,#2){\@tempcnta#1\relax\multiply\@tempcnta 8
\advance\@tempcnta -9 \ifnum #2>0 \advance\@tempcnta #2\relax\else
\advance\@tempcnta -#2\relax\advance\@tempcnta 64 \fi
\char\@tempcnta}

\def\Vektor(#1,#2)#3(#4,#5){\unskip\leavevmode
  \xcoord#4\relax \ycoord#5\relax
      \raise\ycoord \Einheit\hbox to0pt{\hskip\xcoord \Einheit
         \Vector@(#1,#2){#3}\hss}}

\def\Vector@(#1,#2)#3{\@xarg #1\relax \@yarg #2\relax
\@tempcnta \ifnum\@xarg<0 -\@xarg\else\@xarg\fi
\ifnum\@tempcnta<5\relax
\@linelen=#3\Einheit
\ifnum\@xarg =0 \@vvector
  \else \ifnum\@yarg =0 \@hvector \else \@svector\fi
\fi
\else\@badlinearg\fi}

\def\@hvector{\@hline\hbox to 0pt{\@linefnt
\ifnum \@xarg <0 \@getlarrow(1,0)\hss\else
    \hss\@getrarrow(1,0)\fi}}

\def\@vvector{\ifnum \@yarg <0 \@downvector \else \@upvector \fi}

\def\@svector{\@sline
\@tempcnta\@yarg \ifnum\@tempcnta <0 \@tempcnta=-\@tempcnta\fi
\ifnum\@tempcnta <5
  \hskip -\wd\@linechar
  \@upordown\@clnht \hbox{\@linefnt  \if@negarg
  \@getlarrow(\@xarg,\@yyarg) \else \@getrarrow(\@xarg,\@yyarg) \fi}%
\else\@badlinearg\fi}

\def\@upline{\hbox to \z@{\hskip -.5\Pfadd@cke \vrule width \Pfadd@cke
   height \@linelen depth \z@\hss}}

\def\@downline{\hbox to \z@{\hskip -.5\Pfadd@cke \vrule width \Pfadd@cke
   height \z@ depth \@linelen \hss}}

\def\@upvector{\@upline\setbox\@tempboxa\hbox{\@linefnt\char'66}\raise
     \@linelen \hbox to\z@{\lower \ht\@tempboxa\box\@tempboxa\hss}}

\def\@downvector{\@downline\lower \@linelen
      \hbox to \z@{\@linefnt\char'77\hss}}

\def\@getlarrow(#1,#2){\ifnum #2 =\z@ \@tempcnta='33\else
\@tempcnta=#1\relax\multiply\@tempcnta \sixt@@n \advance\@tempcnta
-9 \@tempcntb=#2\relax\multiply\@tempcntb \tw@
\ifnum \@tempcntb >0 \advance\@tempcnta \@tempcntb\relax
\else\advance\@tempcnta -\@tempcntb\advance\@tempcnta 64
\fi\fi\char\@tempcnta}

\def\@getrarrow(#1,#2){\@tempcntb=#2\relax
\ifnum\@tempcntb < 0 \@tempcntb=-\@tempcntb\relax\fi
\ifcase \@tempcntb\relax \@tempcnta='55 \or
\ifnum #1<3 \@tempcnta=#1\relax\multiply\@tempcnta
24 \advance\@tempcnta -6 \else \ifnum #1=3 \@tempcnta=49
\else\@tempcnta=58 \fi\fi\or
\ifnum #1<3 \@tempcnta=#1\relax\multiply\@tempcnta
24 \advance\@tempcnta -3 \else \@tempcnta=51\fi\or
\@tempcnta=#1\relax\multiply\@tempcnta
\sixt@@n \advance\@tempcnta -\tw@ \else
\@tempcnta=#1\relax\multiply\@tempcnta
\sixt@@n \advance\@tempcnta 7 \fi\ifnum #2<0 \advance\@tempcnta 64 \fi
\char\@tempcnta}

\def\Diagonale(#1,#2)#3{\unskip\leavevmode
  \xcoord#1\relax \ycoord#2\relax
      \raise\ycoord \Einheit\hbox to0pt{\hskip\xcoord \Einheit
         \Line@(1,1){#3}\hss}}
\def\AntiDiagonale(#1,#2)#3{\unskip\leavevmode
  \xcoord#1\relax \ycoord#2\relax 
      \raise\ycoord \Einheit\hbox to0pt{\hskip\xcoord \Einheit
         \Line@(1,-1){#3}\hss}}
\def\Pfad(#1,#2),#3\endPfad{\unskip\leavevmode
  \xcoord#1 \ycoord#2 \thicklines\ZeichnePfad#3\endPfad\thinlines}
\def\ZeichnePfad#1{\ifx#1\endPfad\let\next\relax
  \else\let\next\ZeichnePfad
    \ifnum#1=1
      \raise\ycoord \Einheit\hbox to0pt{\hskip\xcoord \Einheit
         \vrule height\Pfadd@cke width1 \Einheit depth\Pfadd@cke\hss}%
      \advance\xcoord by 1
    \else\ifnum#1=2
      \raise\ycoord \Einheit\hbox to0pt{\hskip\xcoord \Einheit
        \hbox{\hskip-\PfadD@cke\vrule height1 \Einheit width\PfadD@cke depth0pt}\hss}%
      \advance\ycoord by 1
    \else\ifnum#1=3
      \raise\ycoord \Einheit\hbox to0pt{\hskip\xcoord \Einheit
         \Line@(1,1){1}\hss}
      \advance\xcoord by 1
      \advance\ycoord by 1
    \else\ifnum#1=4
      \raise\ycoord \Einheit\hbox to0pt{\hskip\xcoord \Einheit
         \Line@(1,-1){1}\hss}
      \advance\xcoord by 1
      \advance\ycoord by -1
    \else\ifnum#1=5
      \advance\xcoord by -1
      \raise\ycoord \Einheit\hbox to0pt{\hskip\xcoord \Einheit
         \vrule height\Pfadd@cke width1 \Einheit depth\Pfadd@cke\hss}%
    \else\ifnum#1=6
      \advance\ycoord by -1
      \raise\ycoord \Einheit\hbox to0pt{\hskip\xcoord \Einheit
        \hbox{\hskip-\PfadD@cke\vrule height1 \Einheit width\PfadD@cke depth0pt}\hss}%
    \else\ifnum#1=7
      \advance\xcoord by -1
      \advance\ycoord by -1
      \raise\ycoord \Einheit\hbox to0pt{\hskip\xcoord \Einheit
         \Line@(1,1){1}\hss}
    \else\ifnum#1=8
      \advance\xcoord by -1
      \advance\ycoord by +1
      \raise\ycoord \Einheit\hbox to0pt{\hskip\xcoord \Einheit
         \Line@(1,-1){1}\hss}
    \fi\fi\fi\fi
    \fi\fi\fi\fi
  \fi\next}
\def\hSSchritt{\leavevmode\raise-.4pt\hbox to0pt{\hss.\hss}\hskip.2\Einheit
  \raise-.4pt\hbox to0pt{\hss.\hss}\hskip.2\Einheit
  \raise-.4pt\hbox to0pt{\hss.\hss}\hskip.2\Einheit
  \raise-.4pt\hbox to0pt{\hss.\hss}\hskip.2\Einheit
  \raise-.4pt\hbox to0pt{\hss.\hss}\hskip.2\Einheit}
\def\vSSchritt{\vbox{\baselineskip.2\Einheit\lineskiplimit0pt
\hbox{.}\hbox{.}\hbox{.}\hbox{.}\hbox{.}}}
\def\DSSchritt{\leavevmode\raise-.4pt\hbox to0pt{%
  \hbox to0pt{\hss.\hss}\hskip.2\Einheit
  \raise.2\Einheit\hbox to0pt{\hss.\hss}\hskip.2\Einheit
  \raise.4\Einheit\hbox to0pt{\hss.\hss}\hskip.2\Einheit
  \raise.6\Einheit\hbox to0pt{\hss.\hss}\hskip.2\Einheit
  \raise.8\Einheit\hbox to0pt{\hss.\hss}\hss}}
\def\dSSchritt{\leavevmode\raise-.4pt\hbox to0pt{%
  \hbox to0pt{\hss.\hss}\hskip.2\Einheit
  \raise-.2\Einheit\hbox to0pt{\hss.\hss}\hskip.2\Einheit
  \raise-.4\Einheit\hbox to0pt{\hss.\hss}\hskip.2\Einheit
  \raise-.6\Einheit\hbox to0pt{\hss.\hss}\hskip.2\Einheit
  \raise-.8\Einheit\hbox to0pt{\hss.\hss}\hss}}
\def\SPfad(#1,#2),#3\endSPfad{\unskip\leavevmode
  \xcoord#1 \ycoord#2 \ZeichneSPfad#3\endSPfad}
\def\ZeichneSPfad#1{\ifx#1\endSPfad\let\next\relax
  \else\let\next\ZeichneSPfad
    \ifnum#1=1
      \raise\ycoord \Einheit\hbox to0pt{\hskip\xcoord \Einheit
         \hSSchritt\hss}%
      \advance\xcoord by 1
    \else\ifnum#1=2
      \raise\ycoord \Einheit\hbox to0pt{\hskip\xcoord \Einheit
        \hbox{\hskip-2pt \vSSchritt}\hss}%
      \advance\ycoord by 1
    \else\ifnum#1=3
      \raise\ycoord \Einheit\hbox to0pt{\hskip\xcoord \Einheit
         \DSSchritt\hss}
      \advance\xcoord by 1
      \advance\ycoord by 1
    \else\ifnum#1=4
      \raise\ycoord \Einheit\hbox to0pt{\hskip\xcoord \Einheit
         \dSSchritt\hss}
      \advance\xcoord by 1
      \advance\ycoord by -1
    \else\ifnum#1=5
      \advance\xcoord by -1
      \raise\ycoord \Einheit\hbox to0pt{\hskip\xcoord \Einheit
         \hSSchritt\hss}%
    \else\ifnum#1=6
      \advance\ycoord by -1
      \raise\ycoord \Einheit\hbox to0pt{\hskip\xcoord \Einheit
        \hbox{\hskip-2pt \vSSchritt}\hss}%
    \else\ifnum#1=7
      \advance\xcoord by -1
      \advance\ycoord by -1
      \raise\ycoord \Einheit\hbox to0pt{\hskip\xcoord \Einheit
         \DSSchritt\hss}
    \else\ifnum#1=8
      \advance\xcoord by -1
      \advance\ycoord by 1
      \raise\ycoord \Einheit\hbox to0pt{\hskip\xcoord \Einheit
         \dSSchritt\hss}
    \fi\fi\fi\fi
    \fi\fi\fi\fi
  \fi\next}
\def\Koordinatenachsen(#1,#2){\unskip
 \hbox to0pt{\hskip-.5pt\vrule height#2 \Einheit width.5pt depth1 \Einheit}%
 \hbox to0pt{\hskip-1 \Einheit \xcoord#1 \advance\xcoord by1
    \vrule height0.25pt width\xcoord \Einheit depth0.25pt\hss}}
\def\Koordinatenachsen(#1,#2)(#3,#4){\unskip
 \hbox to0pt{\hskip-.5pt \ycoord-#4 \advance\ycoord by1
    \vrule height#2 \Einheit width.5pt depth\ycoord \Einheit}%
 \hbox to0pt{\hskip-1 \Einheit \hskip#3\Einheit 
    \xcoord#1 \advance\xcoord by1 \advance\xcoord by-#3 
    \vrule height0.25pt width\xcoord \Einheit depth0.25pt\hss}}
\def\Gitter(#1,#2){\unskip \xcoord0 \ycoord0 \leavevmode
  \LOOP\ifnum\ycoord<#2
    \loop\ifnum\xcoord<#1
      \raise\ycoord \Einheit\hbox to0pt{\hskip\xcoord \Einheit\Punkt\hss}%
      \advance\xcoord by1
    \repeat
    \xcoord0
    \advance\ycoord by1
  \REPEAT}
\def\Gitter(#1,#2)(#3,#4){\unskip \xcoord#3 \ycoord#4 \leavevmode
  \LOOP\ifnum\ycoord<#2
    \loop\ifnum\xcoord<#1
      \raise\ycoord \Einheit\hbox to0pt{\hskip\xcoord \Einheit\Punkt\hss}%
      \advance\xcoord by1
    \repeat
    \xcoord#3
    \advance\ycoord by1
  \REPEAT}
\def\Label#1#2(#3,#4){\unskip \xdim#3 \Einheit \ydim#4 \Einheit
  \def\lo{\advance\xdim by-.5 \Einheit \advance\ydim by.5 \Einheit}%
  \def\llo{\advance\xdim by-.25cm \advance\ydim by.5 \Einheit}%
  \def\loo{\advance\xdim by-.5 \Einheit \advance\ydim by.25cm}%
  \def\o{\advance\ydim by.25cm}%
  \def\ro{\advance\xdim by.5 \Einheit \advance\ydim by.5 \Einheit}%
  \def\rro{\advance\xdim by.25cm \advance\ydim by.5 \Einheit}%
  \def\roo{\advance\xdim by.5 \Einheit \advance\ydim by.25cm}%
  \def\l{\advance\xdim by-.30cm}%
  \def\r{\advance\xdim by.30cm}%
  \def\lu{\advance\xdim by-.5 \Einheit \advance\ydim by-.6 \Einheit}%
  \def\llu{\advance\xdim by-.25cm \advance\ydim by-.6 \Einheit}%
  \def\luu{\advance\xdim by-.5 \Einheit \advance\ydim by-.30cm}%
  \def\u{\advance\ydim by-.30cm}%
  \def\ru{\advance\xdim by.5 \Einheit \advance\ydim by-.6 \Einheit}%
  \def\rru{\advance\xdim by.25cm \advance\ydim by-.6 \Einheit}%
  \def\ruu{\advance\xdim by.5 \Einheit \advance\ydim by-.30cm}%
  #1\raise\ydim\hbox to0pt{\hskip\xdim
     \vbox to0pt{\vss\hbox to0pt{\hss$#2$\hss}\vss}\hss}%
}
\catcode`\@=13

\hsize13cm
\vsize19cm
\newdimen\fullhsize
\newdimen\fullvsize
\newdimen\halfsize
\fullhsize13cm
\fullvsize19cm
\halfsize=0.5\fullhsize
\advance\halfsize by-0.5em

\magnification1200

\TagsOnRight

\def\AignAB{1}
\def\AmZeAA{2}
\def\BresAO{3}
\def\ComtAA{4}
\def\EgRRAB{5}
\def\EgRRAD{6}
\def\FuKrAA{7}
\def\FuHaAA{8}
\def\GeViAA{9}
\def\GeViAB{10}
\def\GhKrAA{11}
\def\KratBL{12}
\def\KratBN{13}
\def\KratBZ{14}
\def\KratBV{15}
\def\KratBW{16}
\def\LindAA{17}
\def\StanBI{18}
\def\StemAE{19}
\def\VienAE{20}
\def\ZeilZZ{21}

\def\AA{1.1}
\def\AB{1.2}
\def\ABa{1.3}
\def\AC{1.4}
\def\AD{1.5}
\def\ADa{1.6}
\def\ADb{1.7}
\def\ADc{1.8}
\def\ADcc{1.9}
\def\ADd{1.10}
\def\ADdd{1.11}
\def\ADddd{1.12}
\def\AE{2.1}
\def\AF{2.2}
\def\AG{2.3}
\def\AGa{2.4}
\def\SA{2.5}
\def\SB{2.6}
\def\SC{2.7}
\def\SD{2.8}
\def\SE{2.9}
\def\SF{2.10}
\def\SG{3.1}
\def\SH{3.2}
\def\SI{4.1}
\def\SJ{4.2}
\def\SK{4.3}
\def\SL{4.4}
\def\SM{4.5}
\def\AI{5.1}
\def\AIa{5.2}
\def\AIc{5.3}
\def\AId{5.4}
\def\AIb{5.5}
\def\AM{5.6}
\def\AN{5.7}
\def\ANa{5.8}
\def\AH{5.9}
\def\AJ{5.10}
\def\AJa{5.11}
\def\AJb{5.12}
\def\AJc{5.13}
\def\AK{5.14}
\def\AL{5.15}
\def\ALa{5.16}
\def\AO{5.17}
\def\AP{5.18}
\def\AQ{5.19}
\def\AR{5.20}
\def\ASAS{5.21}
\def\AS{5.21a}
\def\AT{5.21b}
\def\AU{5.21c}
\def\AV{5.22}
\def\AW{5.23}
\def\AX{5.24}
\def\AY{5.25}
\def\AZa{5.26}
\def\AZ{5.27}
\def\AZAZ{5.28}
\def\BA{6.1}
\def\BAa{6.2}
\def\BAb{6.3}
\def\BC{6.4}
\def\BD{6.5}
\def\BE{6.6}
\def\BF{6.7}
\def\BG{6.8}
\def\BB{6.9}
\def\BBa{6.10}
\def\BBb{6.11}
\def\BBc{6.12}
\def\BH{6.13}
\def\BHd{7.1}
\def\BHdd{7.2}
\def\BHa{7.3}
\def\BHe{7.4}
\def\BHb{7.5}
\def\BI{8.1}
\def\BJ{8.2}
\def\BK{8.3}
\def\BL{8.4}

\def\TA{1}
\def\TB{2}
\def\TAa{3}
\def\TBa{4}
\def\TC{5}
\def\TD{6}
\def\GV{7}
\def\TE{8}
\def\TF{9}
\def\TG{10}
\def\TH{11}
\def\TI{12}
\def\TJ{13}
\def\TK{14}
\def\TL{15}
\def\TM{16}
\def\TN{17}
\def\TO{18}
\def\TP{19}
\def\TQ{20}
\def\TQa{21}
\def\TR{22}
\def\TS{23}
\def\TSa{24}
\def\TT{25}
\def\TU{26}

\def\FA{1}
\def\FB{2}
\def\FC{3}

\def\P{{\Cal P}}
\def\PP{{\bold P}}
\def\GF{\operatorname{GF}}
\def\sp{{sp}}
\def\so{{so}}
\def\fl#1{\lfloor#1\rfloor}
\def\cl#1{\lceil#1\rceil}
\def\coef#1{\left\langle#1\right\rangle}

\topmatter 
\title Some determinants of path generating functions
\endtitle 
\author J. Cigler and C.~Krattenthaler$^{\dagger}$
\endauthor 
\affil 
Fakult\"at f\"ur Mathematik, Universit\"at Wien,\\
Nordbergstra{\ss}e~15, A-1090 Vienna, Austria.\\
WWW: {\tt http://homepage.univie.ac.at/johann.cigler}\\
WWW: \tt http://www.mat.univie.ac.at/\~{}kratt
\endaffil
\address Fakult\"at f\"ur Mathematik, Universit\"at Wien,
Nordbergstra{\ss}e~15, A-1090 Vienna, Austria.
WWW: \tt http://homepage.univie.ac.at/johann.cigler,
http://www.mat.univie.ac.at/\~{}kratt
\endaddress
\thanks $^\dagger$Research partially supported 
by the Austrian Science Foundation FWF, grants Z130-N13 and S9607-N13,
the latter in the framework of the National Research Network
``Analytic Combinatorics and Probabilistic Number Theory"%
\endthanks

\subjclass Primary 05A19;
 Secondary 05A10 05A15 11C20 15A15
\endsubjclass
\keywords Hankel determinants, Catalan numbers, ballot numbers, 
Motzkin numbers, Motzkin paths, non-intersecting lattice paths
\endkeywords
\abstract 
We evaluate four families of determinants of matrices, where the
entries are sums or differences of generating functions for paths
consisting of up-steps, down-steps and level steps. By specialisation,
these determinant evaluations have numerous corollaries. In
particular, they cover numerous determinant evaluations of
combinatorial numbers --- most notably of Catalan, ballot, and of
Motzkin numbers --- that appeared previously in the literature. 
\endabstract
\endtopmatter
\document

\subhead 1. Introduction\endsubhead
Determinants (and Hankel determinants in particular) of path counting numbers
(respectively, more generally, of path generating functions) appear frequently
in the literature. The reason of this ubiquity is two-fold:
first, via the theory of non-intersecting lattice paths
(cf\. \cite{\GeViAA, \GeViAB, \StemAE}), such determinants represent the solution
to counting problems of combinatorial, probabilistic, or algebraic origin (see e.g\.
\cite{\BresAO, \GhKrAA, \KratBL, \KratBV, \KratBW, \StemAE} and the
references contained therein). Second,
it turns out that such determinants can be often
evaluated into attractive, compact closed formulae. This latter theme will
be the underlying theme of the present paper.

(Hankel) Determinant evaluations such as 
$$\align 
\det_{0\le i,j\le n-1}(C_{i+j})&=1,
\tag\AA\\
\det_{0\le i,j\le n-1}(C_{i+j+1})&=1,
\tag\AB\\
\det_{0\le i,j\le n-1}(C_{i+j+2})&=n+1,
\tag\ABa\\
\endalign$$
where $C_n=\frac {1} {n+1}\binom {2n}n$ is the $n$-th {\it Catalan
number}, and 
$$\align 
\det_{0\le i,j\le n-1}(M_{i+j})&=1,
\tag\AC\\
\det_{0\le i,j\le n-1}(M_{i+j+1})&=\cases 
(-1)^{n/3}&\text{if }n\equiv 0\pmod 3,\\
(-1)^{(n-1)/3}&\text{if }n\equiv 1\pmod 3,\\
0&\text{if }n\equiv 2\pmod 3,\\
\endcases
\tag\AD
\endalign$$
where $M_n=\sum _{k=0} ^{\fl{n/2}}\binom {n} {2k} \frac {1} {k+1}\binom {2k}k$ is
the $n$-th {\it Motzkin number}, belong to the folklore of (orthogonal
polynomials) literature (cf\. e.g\. \cite{\VienAE, \AignAB}). In our
paper, we shall consider common weighted generalisations of these
determinant evaluations. 

It is well-known (cf\. \cite{\StanBI, Exercises~6.19 and 6.38}) that
$C_n$ counts the number of lattice paths from $(0,0)$ to $(2n,0)$
consisting of up-steps $(1,1)$ and
down-steps $(1,-1)$, which never run below the $x$-axis
(see Figure~\FA.a for an example with $n=4$), and that
$M_n$ counts the number of lattice paths from $(0,0)$ to $(n,0)$
consisting of up-steps $(1,1)$, level steps $(1,0)$, and
down-steps $(1,-1)$, which never run below the $x$-axis
(see Figure~\FA.b for an example with $n=11$).
Our weighted generalisations will feature different weights for the
three types of steps in such paths.

\midinsert
$$
\Gitter(9,5)(-1,0)
\Koordinatenachsen(9,5)(-1,0)
\Pfad(0,0),33344344\endPfad
\DickPunkt(0,0)
\DickPunkt(8,0)
\PfadDicke{.5pt}
\hbox{\hskip6.5cm}
\Gitter(12,5)(-1,0)
\Koordinatenachsen(12,5)(-1,0)
\Pfad(0,0),33141433144\endPfad
\DickPunkt(0,0)
\DickPunkt(11,0)
\PfadDicke{.5pt}
\hbox{\hskip5.5cm}
$$
\centerline{\eightpoint a. A Catalan path\quad \kern4cm
b. A Motzkin path}
\vskip7pt
\centerline{\eightpoint Figure \FA}
\endinsert

Let us define $\P_n(l,k)$ as the generating function $\sum _{P} ^{}w(P),$
where $P$ runs over all paths from $(0,l)$ to $(n,k)$ consisting of
steps from $\{(1,0), (1,1), (1,-1)\}$ (for the sake of simplicity,
such paths will in the sequel
be referred to as {\it three-step paths}), and where $w(P)$ is the
product of all weights of the steps of $P$, where the weights of the
steps are defined by $w((1,0))=x+y$, $w((1,1))=1$, and $w((1,-1))=xy$.
Furthermore, let $\P^+_n(l,k)$ be the analogous generating function $\sum _{P}
^{}w(P)$, where $P$ runs over the subset of the set of the above
three-step paths which never run below the $x$-axis.
It should be observed that our choice of edge weights essentially
amounts to giving independent weights to the three kinds of steps of
the paths. The somewhat unusual parametrisation that we have chosen
here will turn out to be useful in presenting our results in more
compact forms than would be possible when using a more straightforward
parametrisation.

Clearly, if we specialise $x=-y=\sqrt{-1}$ 
(in which case $x+y=0$ and $xy=1$, that is, paths consisting entirely
of up- and down-steps are weighted by $1$, while all other paths
acquire vanishing weight), then
$\P^+_{2n}(0,0)$ reduces to $C_n$.
More generally, for this specialisation of $x$ and $y$, the numbers
$\P^+_n(0,k)$ are known as {\it ballot numbers}.
On the other hand,
if we specialise $x=\frac {1} {2}(1+\sqrt{-3})$,
$y=\frac {1} {2}(1-\sqrt{-3})$ (in which case we have $x+y=xy=1$,
that is, all three kinds of steps are weighted by $1$), then
$\P^+_{n}(0,0)$ reduces to $M_n$. A third kind of
specialisation that we shall make use of, which is less intuitive, is
$x=y=1$. In this case, up- and down-steps are weighted by 1, while
level steps are weighted by $2$. It is not difficult to see
(by either using (\AGa) below, or by combinatorial reasoning: each
up-step and each down-step is doubled, while level steps are replaced
by either an up-step followed by a down-step or by a down-step
followed by an up-step) that, for
this specialisation of $x$ and $y$, we have $\P^+_{n}(0,0)=C_{n+1}$ and
$\P_{n}(l,k)=\binom {2n} {n+k-l}$.

We present our results generalising (\AA)--(\AD) in the two theorems
below. 

\proclaim{Theorem \TA}
For all positive integers $n$ and non-negative integers $k$, we have
$$\det_{0\le i,j\le n-1}\(\P^+_{i+j}(0,k)\)=\cases (-1)^{n_1\binom
{k+1}2}(xy)^{(k+1)^2\binom {n_1}2}&n=n_1(k+1),\\
0&n\not\equiv 0~(\text {\rm mod }k+1).
\endcases
\tag\ADa$$
\endproclaim

\proclaim{Theorem \TB}
For all positive integers $n$ and non-negative integers $k$, we have
$$\multline 
\det_{0\le i,j\le n-1}\(\P^+_{i+j+1}(0,k)\)\\
=\cases (-1)^{n_1\binom
{k+1}2}(xy)^{(k+1)^2\binom {n_1}2}
\frac {y^{(k+1)(n_1+1)}-x^{(k+1)(n_1+1)}} {y^{k+1}-x^{k+1}}&n=n_1(k+1),\\
(-1)^{n_1\binom
{k+1}2+\binom k2}(xy)^{(k+1)^2\binom {n_1}2+n_1k(k+1)}\\
\kern2cm\times\frac {y^{(k+1)(n_1+1)}-x^{(k+1)(n_1+1)}}
{y^{k+1}-x^{k+1}}&n=n_1(k+1)+k,\\
0&n\not\equiv 0,k~(\text {\rm mod }k+1).
\endcases
\endmultline
\tag\ADb$$
\endproclaim

\remark{Remark}
If $k=0$, the formulae in Theorems~\TA\ and \TB\ have to be read
according to the convention that only the first line on the right-hand
sides of (\ADa) and (\ADb) applies; that is,
$$\det_{0\le i,j\le n-1}\(\P^+_{i+j}(0,0)\)=
(xy)^{\binom {n}2}
$$
and
$$
\det_{0\le i,j\le n-1}\(\P^+_{i+j+1}(0,0)\)
=(xy)^{\binom {n}2}
\frac {y^{n+1}-x^{n+1}} {y-x}.
$$
\endremark

For $x=-y=\sqrt{-1}$ and $k=0$, by the factorisation of determinants of 
``checkerboard matrices" given in Lemma~\TC\ in Section~3, Theorem~\TA\
implies both (\AA) and (\AB). Moreover, if we set $x=y=1$,
then Theorems~\TA\ and \TB\ reduce to (\AB) and (\ABa), respectively. 
On the other hand, if we specialise $x=\frac {1} {2}(1+\sqrt{-3})$,
$y=\frac {1} {2}(1-\sqrt{-3})$, then Theorems~\TA\ and \TB\ reduce to
(\AC) and (\AD), respectively. We list further interesting special cases
of the above two theorems in Section~7. 

We show furthermore that the analogous Hankel determinants, where the ``restricted" 
path generating functions $\P^+_{i+j}(0,k)$, respectively
$\P^+_{i+j+1}(0,k)$, are replaced by their ``unrestricted" counterparts, have
as well compact evaluations.

\proclaim{Theorem \TAa}
For all positive integers $n$ and $k$, we have
$$\det_{0\le i,j\le n-1}\(\P_{i+j}(0,k)\)=
\cases (-1)^{kn_1+\binom k2}(xy)^{k(n_1-1)(2kn_1-k+1)}&n=2kn_1-k+1,\\
(-1)^{kn_1}(xy)^{kn_1(2kn_1-k-1)}&n=2kn_1,\\
0&n\not\equiv 0,k+1~(\text {\rm mod }2k),
\endcases
\tag\ADc$$
while for $k=0$ we have
$$\det_{0\le i,j\le n-1}\(\P_{i+j}(0,0)\)=
2^{n-1}(xy)^{\binom n2}.
\tag\ADcc$$
\endproclaim

\remark{Remarks}
(1)
If $k=1$, the first two cases on the right-hand side of (\ADc)
coincide, so that we have
$$
\det_{0\le i,j\le n-1}\(\P_{i+j}(0,1)\)=
\cases (-1)^{n_1}(xy)^{2n_1(n_1-1)}&n=2n_1,\\
0&n\text { odd}.
\endcases
$$

\smallskip
(2) By Formula~(\AF), the determinant evaluation in Theorem~\TAa\
also implies a formula for negative $k$. We omit its explicit
statement for the sake of brevity.
\endremark

\proclaim{Theorem \TBa}
For all positive integers $n$ and integers $k\ge2$, we have
$$\multline 
\det_{0\le i,j\le n-1}\(\P_{i+j+1}(0,k)\)\\
=
\cases 
(-1)^{k(n_1-1)-1}(xy)^{kn_1(2kn_1-k-3)+k}P_{n-k+2,k}(x,y)
\hskip-1cm\\
&n=2kn_1-1,\\
(-1)^{kn_1+\binom k2}(xy)^{k(n_1-1)(2kn_1-k+1)}P_{n,k}(x,y)
\hskip-4cm\\
&n=2kn_1-k+1,\\
(-1)^{kn_1+\binom {k+1}2}(xy)^{k(n_1-1)(2kn_1-k-1)}P_{n-k,k}(x,y)
\hskip-4cm\\
&n=2kn_1-k,\\
(-1)^{kn_1}(xy)^{kn_1(2kn_1-k-1)}P_{n,k}(x,y)&n=2kn_1,\\
0\hskip7cm 
n\not\equiv 0,k,k+1,2k-1~(\text {\rm mod }2k),
\hskip-4cm
\endcases\\
\endmultline
\tag\ADd$$
where
$$
P_{m,k}(x,y)
=\cases 
\frac {x^{m+k}+(-1)^{m/k}y^{m+k}} {x^k+y^k}
&\text{if }k\mid m,\\
\frac {\(x^{k\fl{m/k}+k}+(-1)^{\fl{m/k}}y^{k\fl{m/k}+k}\)
\(x^{m-k\fl{m/k}}+(-1)^{\fl{m/k}}y^{m-k\fl{m/k}}\)} {x^k+y^k}
&\text{if }k\nmid m,
\endcases
$$
while for $k=1$ we have
$$
\det_{0\le i,j\le n-1}\(\P_{i+j+1}(0,1)\)
=
\cases 
(-1)^{n_1}(xy)^{2n_1(n_1-1)}
\frac {x^{n+1}+(-1)^{n}y^{n+1}} {x+y}
&n=2n_1,\\
(-1)^{n_1+1}(xy)^{2(n_1-1)^2}
\frac {x^{n}+(-1)^{n-1}y^{n}} {x+y}
&n=2n_1-1,
\endcases
\tag\ADdd$$
and for $k=0$ we have
$$
\det_{0\le i,j\le n-1}\(\P_{i+j+1}(0,0)\)
=2^{n-1}(xy)^{\binom n2}(x^n+y^n).
\tag\ADddd$$
\endproclaim

\remark{Remarks}
(1) Again, by Formula~(\AF), the determinant evaluations in Theorem~\TBa\
also imply formulae for negative $k$. We omit their explicit
statement for the sake of brevity.

\smallskip
(2) Inspection of those values of $n$ in (\ADd) which lead to non-zero
determinants shows that it suffices to use the following, restricted,
definition for $P_{m,k}(x,y)$: 
$$\multline
P_{m,k}(x,y)\\
=\cases 
\frac {x^{m+k}+y^{m+k}} {x^k+y^k}
&\text{if }k\mid m\text{ and $m/k$ is even},\\
\frac {\(x^{k\fl{m/k}+k}-y^{k\fl{m/k}+k}\)
\(x^{m-k\fl{m/k}}-y^{m-k\fl{m/k}}\)} {x^k+y^k}
&\text{if }k\nmid m\text{ and $\fl{m/k}$ is odd}.
\endcases
\endmultline
$$

\smallskip
(3) For $\alpha>1$, computer calculations show that the evaluations
of the ``higher order" Hankel determinants
$$
\det_{0\le i,j\le n-1}\(\P^+_{i+j+\alpha}(0,k)\)
\quad \text{and}\quad 
\det_{0\le i,j\le n-1}\(\P_{i+j+\alpha}(0,k)\)
$$
become increasingly unwieldy. Presumably it would be still possible
to work out, and subsequently prove, the corresponding evaluations
for $\alpha=2$, say. However, we did not actually try this.
In any case, we doubt that there is a reasonable
formula for generic $\alpha$.

\smallskip
(4) While, usually, Hankel determinants are intimately related to
orthogonal polynomials (cf\. e.g\. \cite{\KratBN, Sec.~2.7} and
\cite{\KratBZ, Sec.~5.4}), it does not seem to be the case
here (except for $k=0$ and $k=1$, where the determinants in
Theorems~\TA--\TBa\ are related to Chebyshev polynomials), since
the results in Theorems~\TA--\TBa\ follow modular patterns with
a frequent appearance of zeroes, something which is not allowed
in the theory of orthogonal polynomials. 

\smallskip
(5) A similar remark applies to applicability of available computer
packages: the evaluations of the determinants that we consider 
in this paper are certainly not amenable to the condensation method
(cf\. \cite{\KratBN, Sec.~2.3}), and therefore
the package DODGSON by Amdeberhan
and Zeilberger \cite{\AmZeAA} will not be useful here.
Zeilberger's package DET \cite{\ZeilZZ} (which is based on recurrence
methods) can do (\BHdd), but it must necessarily fail if the results
follow modular patterns, which is the case for our determinants if $k$
is not $0$ or $1$. It is conceivable that Zeilberger's algorithmic 
approach in \cite{\ZeilZZ} can be extended, respectively adapted, to
cover also patterns modulo $m$, say, {\it for a fixed $m$}. However,
such an extension could still not treat any of the corollaries in
Section~7 for {\it generic} $k$, it could only give {\it hints} towards
a general proof. An additional new idea is required to be able to
attack the determinant identities of our paper in complete
generality by algorithmic methods.

It should also be pointed out that the determinants in Theorems~\TA--\TBa,
\TE--\TH\ cause yet another problem when attacked by computer packages:
these are determinants the entries of which are polynomials in two,
respectively in three variables. This slows down computations considerably,
up to the effect that it may be impossible to carry them out by
current computer technology.  
\endremark

The purpose of the next two sections is to collect preliminary results
on our three-step paths and on determinants of ``checkerboard
matrices," respectively.
We then show in Section~4, that, by the Lindstr\"om--Gessel--Viennot
theorem, the determinants in Theorems~\TA\ and \TB\
have natural combinatorial interpretations in terms of
{\it non-intersecting lattice paths}. In particular, using
non-intersecting lattice paths, we reduce the
determinants in (\ADa) and (\ADb) to determinants of a similar, but
different kind (see (\SJ) and (\SK)). 
These latter determinants turn out to be special cases
of a more general family of determinants which we evaluate in
Theorems~\TE\ and \TF\ in Section~5. In this sense, these two theorems
are the first two main results of our article. Likewise, we show in
Section~4 that the determinants in Theorems~\TAa\ and \TBa\ are equal to
determinants that are of a very similar form as those in (\SJ) and
(\SK) (see (\SL) and (\SM)). The
second set of main results then consists of Theorems~\TG\ and \TH\
in Section~6, in
which we evaluate two further families of determinants, which
generalise (\SL) and (\SM). Corollaries of our main
results are collected in Section~7. We conclude our article by some
comments and questions (see Section~8). 
The most intriguing perhaps is the speculative
question on a potential relation between our determinants in
Theorems~\TE--\TH\ and Jacobi--Trudi-type formulae for symplectic and
orthogonal characters.

\subhead 2. Some facts about three-step paths\endsubhead
In order to prepare for the proofs of our theorems,
we collect some standard facts about our three-step paths.

By definition of our path generating functions, we have
$$\P_n(l,k)=\P_n(0,k-l)\tag\AE$$
and
$$\P_n(0,k)=(xy)^{-k}\P_n(0,-k).
\tag\AF$$
We shall use simple facts such as $P_n(0,k)=0$ for $n<k$ and
$P_n(0,n)=1$ without further reference frequently in the article.

The reflection principle (see e.g\. \cite{\ComtAA, p.~22}) allows us to
express the generating functions $\P^+_n(l,k)$ for {\it
restricted\/} paths in terms of the generating functions $\P_n(l,k)$
for {\it unrestricted\/} paths, namely by
$$\P^+_n(l,k)=\P_n(l,k)-(xy)^{l+1}\P_n(-l-2,k).\tag\AG$$
By using elementary combinatorial reasoning, the path generating
functions $\P_n(0,k)$ can be expressed in the form
$$\align
\P_n(0,k)&=
\coef{z^k}\(z+(x+y)+\frac {xy} {z}\)^n,\\
&=
\coef{z^0}z^{n-k}\(1+\frac {x} {z}\)^n\(1+\frac {y} {z}\)^n,
\tag\AGa
\endalign$$
where $\coef{z^{m}}f(z)$
denotes the coefficient of $z^{m}$ in the formal
Laurent series $f(z)$. 
From (\AGa), it is easy to derive the explicit formulae
$$\align 
\P_n(0,k)&=\sum _{\ell\ge0} ^{}\binom n
{\ell,\ell+k}(x+y)^{n-2\ell-k}(xy)^{\ell}\\
&=\sum _{\ell\ge0} ^{}\binom n
{\ell}\binom n{n-k-\ell}x^{\ell}y^{n-k-\ell},
\endalign$$
where
$$\binom n{k_1,k_2}=\frac {n!} {k_1!\,k_2!\,(n-k_1-k_2)!}$$
is a trinomial coefficient. Via (\AE) and (\AG), they imply explicit
formulae for $\P_n(l,k)$ and $\P^+_n(l,k)$.

For later use, we record the specialisations that were essentially
already discussed in the Introduction: with $\om$ denoting a primitive
sixth root of unity, we have
$$\align 
\P_n(l,k)\Big\vert_{x=-y=\sqrt{-1}}&=\chi(n+l+k\text{ even})
\binom n{\frac {1} {2}(n+k-l)}
\tag\SA\\
\P^+_n(l,k)\Big\vert_{x=-y=\sqrt{-1}}&=\chi(n+l+k\text{ even})
\(\binom n{\frac {1} {2}(n+k-l)}-\binom n{\frac {1} {2}(n+k+l+2)}\)
\tag\SB\\ 
\P_n(l,k)\Big\vert_{x=y^{-1}=\om}&=\sum _{\ell\ge0} ^{}
\binom n {\ell,\ell+k-l}
\tag\SC\\ 
\P^+_n(l,k)\Big\vert_{x=y^{-1}=\om}&=\sum _{\ell\ge0} ^{}\(
\binom n {\ell,\ell+k-l}-\binom n {\ell,\ell+k+l+2}\)
\tag\SD\\ 
\P_n(l,k)\Big\vert_{x=y=1}&=
\binom {2n} {n+k-l}
\tag\SE\\ 
\P^+_n(l,k)\Big\vert_{x=y=1}&=
\binom {2n} {n+k-l}-\binom {2n} {n+k+l+2},
\tag\SF
\endalign$$
where $\chi(\Cal A)=1$ if $\Cal A$ is true and  $\chi(\Cal A)$=0
otherwise.

\subhead 3. Determinants of ``checkerboard" matrices\endsubhead
By (\SA) and (\SB), if we specialise $x=-y=\sqrt{-1}$ in the
determinants in Theorems~\TA\ or \TB, then we obtain matrices in which
every other entry vanishes; more precisely, either the entries for
which the sum of the row index and the columnn index is even vanish, 
or the entries for which the sum of the row index and the 
columnn index is odd vanish. The next two lemmas record the well-known
(and easy to prove)
factorisations of the determinants of such ``checkerboard" matrices.

\proclaim{Lemma \TC}
Let $M=(M_{i,j})_{0\le i,j\le n-1}$ be a matrix for which $M_{i,j}=0$
whenever $i+j$ is odd. Then
$$
\det_{0\le i,j\le n-1}\big(M_{i,j}\big)=
\det_{0\le i,j\le \fl{(n-1)/2}}\big(M_{2i,2j}\big)\cdot
\det_{0\le i,j\le \fl{(n-2)/2}}\big(M_{2i+1,2j+1}\big).
\tag\SG$$
\endproclaim

\proclaim{Lemma \TD}
Let $M=(M_{i,j})_{0\le i,j\le n-1}$ be a matrix for which $M_{i,j}=0$
whenever $i+j$ is even. Then
$$\multline
\det_{0\le i,j\le n-1}\big(M_{i,j}\big)\\
=\cases
\det_{0\le i,j\le (n-2)/2}\big(M_{2i+1,2j}\big)\cdot
\det_{0\le i,j\le (n-2)/2}\big(M_{2i,2j+1}\big)&\text{if $n$ is
even,}\\
0&\text{if $n$ is
odd.}
\endcases
\endmultline\tag\SH$$
\endproclaim
The point here is that, in particular, if we are in the case of
Lemma~\TC, the knowledge of the determinants on the left-hand side of
(\SG) suffices to recursively calculate the values of the determinants
on the right-hand side of (\SG). The same cannot be done in the situation of
Lemma~\TD\ since every other determinant vanishes. The only exception
occurs if $M$ is a symmetric matrix. In that case, the two
determinants on the right-hand side of (\SH) in the case where $n$ is
even are equal to each other. We may therefore solve for
them. However, then the question of what the correct sign is remains,
since we have to take a square root.
Nevertheless, 
if there should be a ``nice" formula for the determinant on the left-hand
side of (\SH), then one expects that there are also ``nice" formulae for
the determinants on the right-hand side of (\SH).

\subhead 4. Non-intersecting lattice paths\endsubhead
The purpose of this section is to explain how the determinants in
Theorems~\TA\ and \TB\ can be combinatorially interpreted in terms of
non-intersecting lattice paths, and to use this interpretation to
transform them into different determinants, generalisations thereof will
subsequently be evaluated in the next section.

First, we recall the Lindstr\"om--Gessel--Viennot theorem on
non-intersecting lattice\linebreak paths, specialised to our context of
three-step paths. A family $(P_0,P_1,\dots,P_{n-1})$ of three-step paths
$P_i$, $i=0,1,\dots,n-1$, is called {\it non-intersecting}, if no two
paths share a {\it lattice point}. The reader should be well aware at this
point that, in our context, this notion has to be taken with care
since this definition {\it does allow} that two paths cross each other in
{\it non-lattice points}. See Figure~\FB\ for examples. In the figure,
the left half shows a pair of non-intersecting paths, while the two
paths shown in the right half (regardless how we read them) share one (!)
vertex (marked by a circle in the figure), and hence they are not non-intersecting. 

\midinsert
$$
\Gitter(9,5)(-1,0)
\Koordinatenachsen(9,5)(-1,0)
\Pfad(0,0),33344344\endPfad
\Pfad(0,2),14141333\endPfad
\DickPunkt(0,0)
\DickPunkt(0,2)
\DickPunkt(8,3)
\DickPunkt(8,0)
\PfadDicke{.5pt}
\hbox{\hskip6.5cm}
\Gitter(9,5)(-1,0)
\Koordinatenachsen(9,5)(-1,0)
\Pfad(0,0),33141133\endPfad
\Pfad(0,3),443133\endPfad
\DickPunkt(0,0)
\DickPunkt(0,3)
\DickPunkt(6,4)
\DickPunkt(8,3)
\Kreis(3,2)
\PfadDicke{.5pt}
\hbox{\hskip4cm}
$$
\centerline{\eightpoint a. Non-intersecting (!) paths\quad \kern2.5cm
b. Two paths that intersect}
\vskip7pt
\centerline{\eightpoint Figure \FB}
\endinsert

Let us now fix a sublattice $L$ of the plane integer lattice $\Bbb
Z^2$. For our purposes, $L$ will be either all of $\Bbb Z^2$ or the
upper half-plane including the $x$-axis.
Given lattice points $A$ and $E$, we write $\PP(A\to E)$ for the set
of three-step paths from $A$ to $E$ that stay in $L$. More generally,
given $n$-tuples $\bold A=(A_0,A_1,\dots,A_{n-1})$ and
$\bold E=(E_0,E_1,\dots,E_{n-1})$ of lattice points, we write
$\PP(\bold A\to\bold E)$ for the set of families
$(P_0,P_1,\dots,P_{n-1})$ of three-step paths that stay in $L$, 
where path $P_i$ runs
from $A_i$ to $E_i$, $i=0,1,\dots,n-1$, and we write
$\PP^{\text{nonint}}(\bold A\to\bold E)$ for the subset of $\PP(\bold
A\to\bold E)$ of non-intersecting path families. 
In order to not overload notation, we do not make the dependence on
$L$ explicit in the symbols $\PP(A\to E)$, etc. 

We extend the path
weight $w(\,.\,)$ of the Introduction to path families by
$$w\big((P_0,P_1,\dots,P_{n-1})\big):=
\prod _{i=0} ^{n-1}w(P_i).$$
Finally, given a set $M$ with weight function $w$, we write 
$\GF(\Cal M;w)$ for the generating function $\sum
_{x\in\Cal M} ^{}w(x)$. 

We are now in the position to state the Lindstr\"om--Gessel--Viennot
theorem. There, the symbol $S_n$
denotes the group of permutations of $\{0,1,\dots,n-1\}$, 
and, given a permutation $\si\in
S_n$, we write $\bold E_\si$ for
$(E_{\si(1)},E_{\si(2)},\dots,E_{\si(n)})$. 

\proclaim{Theorem \GV\ (\cite{\LindAA, \smc Lemma~1},
\cite{\GeViAB, \smc Theorem~1})}
Let $L$ be a fixed sublattice of $\Bbb Z^2$.
For all positive integers $n$ and
$n$-tuples $\bold A=(A_0,A_1,\dots,A_{n-1})$,
$\bold E=(E_0,E_1,\dots,E_{n-1})$ of lattice points, we have 
$$\det_{0\le i,j\le n-1}\big(\GF(\PP(A_{j}\to E_i);w)\big)
=\sum _{\si\in S_n} ^{}(\sgn\si) \cdot\GF(\PP^{\text{\rm nonint}}
(\bold A\to \bold
E_\si);w).
\tag\SI$$
\endproclaim

If we specialise the above theorem to $L$ being the upper half-plane,
$A_i=(-i,0)$ and $E_i=(i,k)$, $i=0,1,\dots,n-1$, then we see that the
determinant in Theorem~\TA\ can be interpreted in terms of
non-intersecting lattice paths. More precisely, with the above choice
of $L$, of the $A_i$'s, of the $E_i$'s, and of the weight $w(\,.\,)$
introduced in the Introduction, we have
$$
\det_{0\le i,j\le n-1}\(\P^+_{i+j}(0,k)\)
=\sum _{\si\in S_n} ^{}(\sgn\si) \cdot\GF(\PP^{\text{\rm nonint}}
(\bold A\to \bold
E_\si);w).
$$

\midinsert
$$\hbox{\hskip2.5cm}
\Einheit.5cm
\Gitter(6,6)(-5,0)
\Koordinatenachsen(6,6)(-5,0)
\Pfad(0,0),333\endPfad
\Pfad(-1,0),34333\endPfad
\Pfad(-2,0),3313\endPfad
\Pfad(-3,0),333\endPfad
\Pfad(-4,0),33334\endPfad
\DickPunkt(0,0)
\DickPunkt(-1,0)
\DickPunkt(-2,0)
\DickPunkt(-3,0)
\DickPunkt(-4,0)
\DickPunkt(0,3)
\DickPunkt(1,3)
\DickPunkt(2,3)
\DickPunkt(3,3)
\DickPunkt(4,3)
\Label\u{A_0}(0,0)
\Label\u{A_1}(-1,0)
\Label\u{A_2}(-2,0)
\Label\u{A_3}(-3,0)
\Label\u{A_4}(-4,0)
\Label\o{E_0}(0,3)
\Label\o{E_1}(1,3)
\Label\o{E_2}(2,3)
\Label\o{E_3}(3,3)
\Label\o{E_4}(4,3)
\PfadDicke{.5pt}
\hbox{\hskip6.5cm}
\Gitter(6,6)(-5,0)
\Koordinatenachsen(6,6)(-5,0)
\Pfad(0,0),333\endPfad
\Pfad(0,1),4333\endPfad
\Pfad(0,2),13\endPfad
\Pfad(0,3),\endPfad
\Pfad(0,4),4\endPfad
\DickPunkt(0,0)
\DickPunkt(0,1)
\DickPunkt(0,2)
\DickPunkt(0,3)
\DickPunkt(0,4)
\DickPunkt(0,3)
\DickPunkt(1,3)
\DickPunkt(2,3)
\DickPunkt(3,3)
\DickPunkt(4,3)
\Label\l{A'_0}(0,0)
\Label\l{A'_1}(0,1)
\Label\l{A'_2}(0,2)
\Label\l{A'_3}(0,3)
\Label\l{A'_4}(0,4)
\Label\o{E_0}(0,3)
\Label\o{E_1}(1,3)
\Label\o{E_2}(2,3)
\Label\o{E_3}(3,3)
\Label\o{E_4}(4,3)
\hbox{\hskip2.5cm}
$$
\centerline{\eightpoint a. A family of non-intersecting paths\quad \kern2cm
b. Omitting the superfluous initial portions}
\vskip7pt
\centerline{\eightpoint Figure \FC}
\endinsert

Now, if we consider a family $(P_0,P_1,\dots,P_{n-1})$ 
of paths in $\PP^{\text{\rm nonint}}
(\bold A\to \bold E_\si)$ that occurs on the right-hand
side (see Figure~\FC.a for an example with $n=5$, $k=3$, and $\si=34201$),
then we see that the first $i$ steps of $P_i$, $i=0,1,\dots,n-1$,
must all be up-steps since the path family is
non-intersecting. Therefore we may equally well omit these steps. 
Thereby, we obtain again a family of non-intersecting paths, say 
$(P'_0,P'_1,\dots,P'_{n-1})$, where $P'_i$ runs from $A'_i=(0,i)$ to
$E_{\si(i)}$. (Figure~\FC.b shows the family of non-intersecting paths
that is obtained in this way from the path family in Figure~\FC.a). 
Reading Theorem~\GV\ in the other direction, this argument
implies the equality
$$\det_{0\le i,j\le n-1}\(\P^+_{i+j}(0,k)\)
=\det_{0\le i,j\le n-1}\(\P^+_{j}(i,k)\).
\tag\SJ$$
An analogous argument establishes the equality
$$\det_{0\le i,j\le n-1}\(\P^+_{i+j+1}(0,k)\)
=\det_{0\le i,j\le n-1}\(\P^+_{j+1}(i,k)\).
\tag\SK$$
The determinants on the right-hand sides of (\SJ) and (\SK) will be
evaluated in the next section in Theorems~\TE\ and \TF, respectively,
thereby establishing Theorems~\TA\ and \TB.

\medskip
As we announced in the Introduction, also the determinants in
Theorems~\TAa\ and \TBa\ can be shown to equal different determinants,
which are very close to the determinants in (\SJ) and (\SK). 
We start with the determinant in (\ADc).
By cutting paths after $i$ steps, it is easy to see that the equation
$$\P_{i+j}(0,k)=
\sum _{\ell=-i} ^{i}\P_i(0,\ell)\P_j(\ell,k)$$
holds. We substitute this in the determinant in (\ADc), to obtain
$$\align
\det_{0\le i,j\le n-1}&\(\sum _{\ell=-i}
^{i}\P_i(0,\ell)\P_j(\ell,k)\)\\
&=
\det_{0\le i,j\le n-1}\(\P_i(0,0)\P_j(0,k)+\sum _{\ell=1}
^{i}\P_i(0,\ell)\P_j(\ell,k)+\sum _{\ell=-i}
^{-1}\P_i(0,\ell)\P_j(\ell,k)\)\\
&=
\det_{0\le i,j\le n-1}\(\P_i(0,0)\P_j(0,k)+\sum _{\ell=1}
^{i}\P_i(0,\ell)\(\P_j(\ell,k)+(xy)^{\ell}\P_j(-\ell,k)\)\),
\endalign
$$
where we used (\AF) to arrive at the last line.
Here, empty sums must be understood as $0$, so that the entry in row
$0$ and column $j$ is equal to $\P_j(0,k)$. We now use row $0$ to
eliminate the term $\P_i(0,0)\P_j(0,k)$ in rows
$i=1,2,\dots,n-1$. Thereby, the entry in row $1$ and column $j$ becomes
$$\P_1(0,1)\big(\P_j(1,k)+xy\P_j(-1,k)\big)=
\P_j(1,k)+xy\P_j(-1,k).$$ 
Hence, row $1$ can be used to
eliminate the terms for $\ell=1$ in the sums over $\ell$ in rows
$1,2,\dots,n-1$. Etc. At the end, we obtain that
$$\det_{0\le i,j\le n-1}\(\P_{i+j}(0,k)\)
=\frac {1} {2}\det_{0\le i,j\le n-1}\(\P_{j}(i,k)+(xy)^i\P_{j}(-i,k)\).
\tag\SL$$
(The reader should note that the fraction $\frac {1} {2}$ comes from
the fact that, written in the above form, in the determinant on the
right-hand side the entry in row $0$ and column $j$ is $2\P_{j}(0,k)$
instead of $\P_{j}(0,k)$.)

An analogous argument establishes the equality
$$\det_{0\le i,j\le n-1}\(\P_{i+j+1}(0,k)\)
=\frac {1} {2}\det_{0\le i,j\le n-1}\(\P_{j+1}(i,k)+(xy)^i\P_{j+1}(-i,k)\).
\tag\SM$$
The determinants on the right-hand sides of (\SL) and (\SM) will be
evaluated in Section~6 in Theorems~\TG\ and \TH, respectively,
thereby establishing Theorems~\TAa\ and \TBa.

\subhead 5. Main theorems, I\endsubhead
In Theorems~\TE\ and \TF\ below, we evaluate two families of
determinants in which the entries are (essentially) differences of path
generating functions. By (\AG), (\SJ), and (\SK), 
the special case $t=1$ of these two
theorems implies Theorems~\TA\ and \TB, respectively.

\proclaim{Theorem \TE}
For all positive integers $n$ and non-negative integers $k$, we have
$$\multline
\det_{0\le i,j\le n-1}\(\P_{j}(i,k)-t(xy)^{i+1}\P_{j}(-i-2,k)\)\\=
\cases (-1)^{n_1\binom
{k+1}2}t^{k\fl{\frac {n_1}2}}(xy)^{(k+1)^2\binom {n_1}2}&n=n_1(k+1),\\
0&n\not\equiv 0~(\text {\rm mod }k+1),
\endcases
\endmultline\tag\AI$$
while for $k=-1$ we have
$$
\det_{0\le i,j\le n-1}\(\P_{j}(i,-1)-t(xy)^{i+1}\P_{j}(-i-2,-1)\)=0.
\tag\AIa$$
\endproclaim

\remark{Remarks}
(1)
If $k=0$, the formula in Theorem~\TE\ has to be read
according to the convention that only the first line on the right-hand
side of (\AI) applies; that is,
$$
\det_{0\le i,j\le n-1}\(\P_{j}(i,0)-t(xy)^{i+1}\P_{j}(-i-2,0)\)=
(xy)^{\binom {n}2}.
\tag\AIc$$

\smallskip
(2) By Formula~(\AF), the determinant evaluation in Theorem~\TE\
also implies a formula for negative $k<-1$. More precisely, 
using as well (\AE), we have
$$\align 
\P_{j}(i,-k)-t(xy)^{i+1}\P_{j}(-i-{}&2,-k)
=\P_{j}(0,-i-k)-t(xy)^{i+1}\P_{j}(0,i+2-k)\\
&=(xy)^{i+k}\P_{j}(0,i+k)-t(xy)^{k-1}\P_{j}(0,k-i-2)\\
&=-t(xy)^{k-1}\left(\P_{j}(i,k-2)
-t^{-1}(xy)^{i+1}\P_{j}(-i-2,k-2)\right).
\tag\AId
\endalign$$
Aside from some trivial factors, the expression in the last line is
again in the form as the expression for the matrix element of the
determinant on the left-hand side of (\AI).
We omit the explicit statement of the resulting formula.
\endremark

\demo{Proof}
If $k=-1$, then the matrix 
$$
\(\P_{j}(i,k)-t(xy)^{i+1}\P_{j}(-i-2,k)\)_{0\le i,j\le n-1}\ ,
\tag\AIb$$
of which we want to compute the
determinant, is upper triangular with zeroes on the main diagonal.
Hence, its determinant vanishes.

If $k=0$, then the matrix (\AIb) is upper triangular, and the entry on the
main diagonal in the $i$-th row is
$$
\P_{i}(i,0)=(xy)^i,
$$
$i=0,1,\dots,n-1$.
The assertion in this case, given explicitly in (\AIc), follows
immediately.

From now on let $k\ge1$.
In the matrix (\AIb), we replace row $(h(2k+2)+b)$ by
$$\multline 
\sum _{\ell=0} ^{h}t^{\ell-h}(xy)^{(h-\ell)(k+1)}\cdot
\big(\text {row $(\ell(2k+2)+b)$}\big)\\
-\sum _{\ell=1} ^{h}t^{\ell-h-1}(xy)^{(h-\ell)(k+1)+b+1}\cdot
\big(\text {row $(\ell(2k+2)-b-2)$}\big)
\endmultline\tag\AM$$
if $0\le b\le k-1$, and by
$$\multline 
\sum _{\ell=0} ^{h}t^{\ell-h}(xy)^{(h-\ell)(k+1)}\cdot
\big(\text {row $(\ell(2k+2)+b)$}\big)\\
-\sum _{\ell=1} ^{h+1}t^{\ell-h-1}(xy)^{(h-\ell)(k+1)+b+1}\cdot
\big(\text {row $(\ell(2k+2)-b-2)$}\big)
\endmultline\tag\AN$$
if $k+1\le b\le 2k$. It is easy to see that this can be achieved by
elementary row manipulations: one starts with the last row, and one
works one's way up. (The reader should keep in mind that the rows are
labelled by $0,1,2,\dots$.) 

Let first $0\le b\le k-1$.
In the new matrix, the entry in the $j$-th column of row $i=h(2k+2)+b$
is equal to
$$\align 
&\sum _{\ell=0} ^{h}t^{\ell-h}(xy)^{(h-\ell)(k+1)}\bigg(
\P_{j}(\ell(2k+2)+b,k)\\
&\kern5cm
-t(xy)^{\ell(2k+2)+b+1}\P_j(-\ell(2k+2)-b-2,k)\bigg)\\
&\quad 
-\sum _{\ell=1} ^{h}t^{\ell-h-1}(xy)^{(h-\ell)(k+1)+b+1}\bigg(
\P_{j}(\ell(2k+2)-b-2,k)\\
&\kern5cm
-t(xy)^{\ell(2k+2)-b-1}\P_j(-\ell(2k+2)+b,k)\bigg)\\
&=\sum _{\ell=0} ^{h}t^{\ell-h}(xy)^{(h-\ell)(k+1)}
\P_{j}(0,-\ell(2k+2)-b+k)\\
&\kern2cm
-\sum _{\ell=0} ^{h}t^{\ell-h+1}(xy)^{(h+\ell)(k+1)+b+1}\P_j(0,\ell(2k+2)+b+k+2)\\
&\kern2cm
-\sum _{\ell=1} ^{h}t^{\ell-h-1}(xy)^{(h-\ell)(k+1)+b+1}
\P_{j}(0,-\ell(2k+2)+b+k+2)\\
&\kern2cm
+\sum _{\ell=1} ^{h}t^{\ell-h}(xy)^{(h+\ell)(k+1)}\P_j(0,\ell(2k+2)-b+k),
\endalign$$
where we used (\AE) repeatedly. If we subsequently apply (\AF) to
further simplify this expression, then we obtain
$$\align
&\sum _{\ell=0} ^{h}t^{\ell-h}(xy)^{(h+\ell)(k+1)+b-k}
\P_{j}(0,\ell(2k+2)+b-k)\\
&\kern2cm
-\sum _{\ell=1} ^{h+1}t^{\ell-h}(xy)^{(h+\ell)(k+1)+b-k}\P_j(0,\ell(2k+2)+b-k)\\
&\kern2cm
-\sum _{\ell=1} ^{h}t^{\ell-h-1}(xy)^{(h+\ell)(k+1)-k-1}
\P_{j}(0,\ell(2k+2)-b-k-2)\\
&\kern2cm
+\sum _{\ell=2} ^{h+1}t^{\ell-h-1}(xy)^{(h+\ell-1)(k+1)}\P_j(0,\ell(2k+2)-b-k-2)\\
&=t^{-h}(xy)^{h(k+1)+b-k}
\P_{j}(0,b-k)
-t(xy)^{(2h+1)(k+1)+b-k}\P_j(0,(h+1)(2k+2)+b-k)\\
&\kern1cm
-t^{-h}(xy)^{h(k+1)}\P_{j}(0,-b+k)
+(xy)^{2h(k+1)}\P_j(0,(h+1)(2k+2)-b-k-2)\\
&=-t(xy)^{(2h+1)(k+1)+b-k}\P_j(0,(h+1)(2k+2)+b-k)\\
&\kern4cm
+(xy)^{2h(k+1)}\P_j(0,h(2k+2)-b+k)
\tag\ANa
\endalign$$
for the $(i,j)$-entry of the new matrix, with $i=h(2k+2)+b$, $0\le
b\le k-1$.
An analogous calculation yields for the case $k+1\le b\le 2k$ that the
entry in the $j$-th column of
row $i=h(2k+2)+b$ in the new matrix equals
$$\multline 
-t(xy)^{(2h+1)(k+1)+b-k}\P_j(0,(h+1)(2k+2)+b-k)\\
+t(xy)^{(2h+1)(k+1)}\P_j(0,(h+1)(2k+2)-b+k).
\endmultline$$
The reader should recall that we did {\it not\/} change the $(i,j)$-entry if
$i\equiv k$~$(\text{mod }k+1)$, say $i=H(k+1)+k$, so that these entries are still given
by 
$$\multline 
\P_j(i,k)-t(xy)^{i+1}\P_j(-i-2,k)\\
=
(xy)^{H(k+1)}\P_j(0,H(k+1))-t(xy)^{(H+1)(k+1)}\P_j(0,(H+2)(k+1)),
\endmultline$$
which fits nicely with (\ANa) if $H=2h$.

In particular, this means that the $(i,j)$-entry, with $i=h(2k+2)+b$,
vanishes in the case where $0\le b\le k$ whenever
$j<h(2k+2)-b+k$, and that it vanishes in the case where $k+1\le b\le 2k+1$
whenever $j<(h+1)(2k+2)-b+k$. 
Hence, if
$n=h(2k+2)+b$ with $1\le b\le k$, then row
$h(2k+2)$ consists entirely of zeroes since $n-1<h(2k+2)+k$. 
Similarly, if $n=h(2k+2)+b$ with $k+1\le b\le 2k+1$, then row
$h(2k+2)+k+1$ consists entirely of zeroes since $n-1<(h+1)(2k+2)-1$.
Consequently, the determinant equals zero in the case
where $n\hbox{${}\not\equiv{}0$}$~(mod~$k+1$).

If $n=n_1(k+1)$, then one can transform the matrix which we have
obtained by the above manipulations into an upper triangular matrix,
using the permutation of the rows given by
$$i=h(2k+2)+b\mapsto
\cases 
h(2k+2)-b+k&0\le b\le k,\\
(h+1)(2k+2)-b+k&k+1\le b\le 2k+1,
\endcases$$
$0\le i\le n-1$, or, in simpler terms,
$$i=H(k+1)+b\mapsto
H(k+1)-b+k,\quad \quad 0\le b\le k,
\tag\AH$$
$0\le i\le n-1$. Reading the entries 
along the main diagonal of this upper triangular matrix, we find
$$\align 
&1,1,\dots,1,\\
&(xy)^{k+1},t(xy)^{k+1}\dots,t(xy)^{k+1},\\
&(xy)^{2(k+1)},(xy)^{2(k+1)},\dots,(xy)^{2(k+1)},\\
&(xy)^{3(k+1)},t(xy)^{3(k+1)}\dots,t(xy)^{3(k+1)},\\
&(xy)^{4(k+1)},(xy)^{4(k+1)}\dots,(xy)^{4(k+1)},\\
&\hbox to 5cm{\leaders \hbox{.}\hfil}\\
&(xy)^{(n_1-1)(k+1)},t^{\chi(n_1\text{ even})}(xy)^{(n_1-1)(k+1)},
\dots,t^{\chi(n_1\text{ even})}(xy)^{(n_1-1)(k+1)},
\endalign$$
where, when arranged as above, there are exactly $k+1$ entries in each
line. The notation $\chi(\,.\,)$ that we used in the last line has the
same meaning as at the end of Section~2:
$\chi(\Cal A)=1$ if $\Cal A$ is true and  $\chi(\Cal A)$=0
otherwise. The product of these entries is
$t^{k\fl{n_1/2}}(xy)^{(k+1)^2\binom {n_1}2}$, which is in agreement 
with our claim.

In order to determine the correct sign in front of the expression on
the right-hand side of (\AI), we must compute the sign of the
permutation in (\AH). The number of inversions of this permutation is
 $n_1\binom {k+1}2.$ Hence, the sign to be determined is
$(-1)^{n_1\binom {k+1}2}$.\quad \quad
\qed
\enddemo

\proclaim{Theorem \TF}
For all positive integers $n$ and non-negative integers $k$, we have
$$\multline 
\det_{0\le i,j\le n-1}\(\P_{j+1}(i,k)-t(xy)^{i+1}\P_{j+1}(-i-2,k)\)\\=
\cases (-1)^{n_1\binom
{k+1}2}t^{k\fl{\frac {n_1}2}}(xy)^{(k+1)^2\binom {n_1}2}\\
\kern2cm\times
\sum _{s=0} ^{n_1}t^{\min\{s,n_1-s\}}x^{s(k+1)}y^{(n_1-s)(k+1)}
&n=n_1(k+1),\\
(-1)^{n_1\binom
{k+1}2+\binom k2}t^{k\cl{\frac {n_1}2}}(xy)^{(k+1)^2\binom {n_1}2+n_1k(k+1)}\\
\kern2cm\times\sum _{s=0}
^{n_1}t^{\min\{s,n_1-s\}}x^{s(k+1)}y^{(n_1-s)(k+1)}&n=n_1(k+1)+k,\\
0&n\not\equiv 0,k~(\text {\rm mod }k+1),
\endcases\\
\endmultline\tag\AJ$$
while for $k=-1$ we have
$$
\det_{0\le i,j\le n-1}\(\P_{j+1}(i,-1)-t(xy)^{i+1}\P_{j+1}(-i-2,-1)\)=
(1-t)^n(xy)^{\binom {n+1}2}.
\tag\AJa$$
\endproclaim

\remark{Remarks}
(1)
If $k=0$, the first formula in Theorem~\TF\ has to be read
according to the convention that only the first and second lines
(which coincide) on the right-hand
side of (\AJ) apply; that is,
$$
\det_{0\le i,j\le n-1}\(\P_{j+1}(i,0)-t(xy)^{i+1}\P_{j+1}(-i-2,0)\)=
(xy)^{\binom {n}2}
\sum _{s=0} ^{n}t^{\min\{s,n-s\}}x^{s}y^{n-s}.
\tag\AJb$$

\smallskip
(2) Also here, by Formula~(\AF), the determinant evaluation in
Theorem~\TF\
also implies a formula for negative $k<-1$. To see this, all one has
to do is to replace $j$ by $j+1$ in the calculation (\AId).
We omit the explicit statement of the resulting formula.
\endremark

\demo{Proof}
If $k=-1$, then the matrix 
$$
\(\P_{j+1}(i,k)-t(xy)^{i+1}\P_{j+1}(-i-2,k)\)_{0\le i,j\le n-1}\ ,
\tag\AJc$$
of which we want to compute the
determinant, is upper triangular and the entry on the
main diagonal in the $i$-th row is
$$
\P_{i+1}(i,-1)=(1-t)(xy)^{i+1},
$$
$i=0,1,\dots,n-1$.
The assertion (\AJa) follows immediately.

For the remaining cases,
we proceed in the same way as in the proof of Theorem~\TE.
Let for the moment $k\ge1$.
We apply the row operations described in (\AM) and (\AN).
In this way,
we obtain a new matrix, where the $(i,j)$-entry, with $i=h(2k+2)+b$,
of the new matrix is given by
$$\multline 
-t(xy)^{(2h+1)(k+1)+b-k}\P_{j+1}(0,(h+1)(2k+2)+b-k)\\
+(xy)^{2h(k+1)}\P_{j+1}(0,h(2k+2)-b+k)
\endmultline\tag\AK$$
if $0\le b\le k$, by
$$\multline 
-t(xy)^{(2h+1)(k+1)+b-k}\P_{j+1}(0,(h+1)(2k+2)+b-k)\\
+t(xy)^{(2h+1)(k+1)}\P_{j+1}(0,(h+1)(2k+2)-b+k)
\endmultline\tag\AL$$
if $k+1\le b\le 2k$, and by
$$
(xy)^{(2h+1)(k+1)}\P_{j+1}(0,(2h+1)(k+1))
-t(xy)^{(2h+2)(k+1)}\P_{j+1}(0,(2h+3)(k+1))
\tag\ALa$$
if $b=2k+1$.
It should be noted that, if $k=0$, the definitions (\AK)--(\ALa) of
the new matrix entries (the case (\AL) being empty) coincide with the
original matrix entries for $k=0$. This allows us to continue with
(\AK)--(\ALa), assuming $k\ge0$.

In particular, from our new matrix we can read off 
that the $(i,j)$-entry, with $i=h(2k+2)+b$,
vanishes in the case where $0\le b\le k$ whenever
$j<h(2k+2)-b+k-1$, and that it vanishes in the case where $k+1\le b\le 2k+1$
whenever $j<(h+1)(2k+2)-b+k-1$. 
Hence, if
$n=h(2k+2)+b$ with $1\le b\le k-1$, then row
$h(2k+2)$ consists entirely of zeroes since $n-1<h(2k+2)+k-1$. 
Similarly, if $n=h(2k+2)+b$ with $k+2\le b\le 2k$, then row
$h(2k+2)+k+1$ consists entirely of zeroes since $n-1<(h+1)(2k+2)-2$.
Consequently, the determinant equals zero in the case where
$n\hbox{${}\not\equiv{}0,k$}$~(mod~$k+1$).

\medskip
Let now $n=n_1(k+1)$. We rearrange the rows of the matrix we have
obtained after the above manipulations according to the permutation
(\AH). This time, we do not obtain an upper triangular matrix, but an
``almost" upper triangular matrix, by which we mean a matrix
$(M_{i,j})$ for which $M_{i,j}=0$ if $j<i-1$. We now factor
$t^{\chi(h\text{ odd})}(xy)^{h(k+1)}$ from all the entries in rows
$h(k+1)+1,\dots,h(k+1)+k$, and we factor $(xy)^{h(k+1)}$ from the
entries in the rows $h(k+1)$, $h=0,1,\dots,n_1-1$. This yields an
overall factor of
$$t^{\fl{\frac {n_1} {2}}}(xy)^{(k+1)^2\binom {n_1}2}
\tag\AO
$$ 
by which we have to multiply the determinant of the remaining matrix
in the end. We must as well multiply by the sign
$$(-1)^{n_1\binom {k+1}2}
\tag\AP$$ 
of the permutation (\AH). 

The remaining matrix is the following matrix: its $(i,j)$-entry,
with $i=h(k+1)+b$ and $0\le b\le k$, is given by
(see (\AK)--(\ALa))
$$
\P_{j+1}(0,h(k+1)+b)
-t^{\chi(h\text{ even or }b\equiv 0\pmod{k+1})}(xy)^{k-b+1}\P_{j+1}(0,(h+2)(k+1)-b).
\tag\AQ
$$
We should observe that, for $i\ge1$, the first
non-zero entry in row $i$ (which is to be found in column $i-1$)
equals $1$.

In this matrix, we replace the $0$-th row by
$$\sum _{h=0} ^{n_1-1}\sum _{b=0} ^{k} (-1)^{h(k+1)+b}
\sum _{s=0} ^{h}c(h,b,s)\,x^{s(k+1)}y^{(h-s)(k+1)}
\cdot
\big(\text {row $(h(k+1)+b)$}\big),\tag\AR$$
where the coefficients $c(h,b,s)$ are given by
$$
c(h,b,s)=\cases 
t^{\min\{s+\chi(h\text{ odd}),h-s\}}x^b+
t^{\min\{s,h-s+\chi(h\text{ odd})\}}y^b,
&\text{if }b\ne0,\\
t^{\min\{s,h-s\}},
&\text{if }b=0.
\endcases
$$
Since the coefficient of the $0$-th row in the linear combination
(\AR) is $1$, this does not change the value of the determinant.

The $(0,j)$-entry in the new matrix is then given by
$$\align 
&\sum _{h=0} ^{n_1-1}\sum _{b=0} ^{k} (-1)^{h(k+1)+b}
\sum _{s=0} ^{h}c(h,b,s)\,x^{s(k+1)}y^{(h-s)(k+1)}\\
&\ \cdot
\big(\P_{j+1}(0,h(k+1)+b)
-t^{\chi(h\text{ even or }b\equiv 0\pmod{k+1})}(xy)^{k-b+1}
\P_{j+1}(0,(h+2)(k+1)-b)\big)\\
&=\sum _{h=0} ^{n_1-1}\sum _{s=0} ^{h} (-1)^{h(k+1)}
t^{\min\{s,h-s\}}x^{s(k+1)}y^{(h-s)(k+1)}\\
&\kern4cm 
\cdot
\big(\P_{j+1}(0,h(k+1))
-t(xy)^{k+1}\P_{j+1}(0,(h+2)(k+1))\big)
\tag\AS\\
&\kern1cm
+\sum _{h=0} ^{n_1-1}\sum _{b=1} ^{k}\sum _{s=0} ^{h} 
(-1)^{h(k+1)+b}
x^{s(k+1)}y^{(h-s)(k+1)}\\
&\kern2.5cm
\cdot
\Big(t^{\min\{s+\chi(h\text{ odd}),h-s\}}x^b+
t^{\min\{s,h-s+\chi(h\text{ odd})\}}y^b\Big)
\P_{j+1}(0,h(k+1)+b)
\tag\AT\\
&\kern1cm
-\sum _{h=0} ^{n_1-1}\sum _{b=1} ^{k}\sum _{s=0} ^{h} 
(-1)^{h(k+1)+b}
x^{s(k+1)}y^{(h-s)(k+1)}\\
&\kern2cm 
\cdot
\Big(t^{\min\{s+1,h-s+\chi(h\text{ even})\}}x^{k+1}y^{k-b+1}+
t^{\min\{s+\chi(h\text{ even}),h-s+1\}}x^{k-b+1}y^{k+1}\Big)\\
&\kern2cm
\cdot
\P_{j+1}(0,(h+2)(k+1)-b).
\tag\AU
\endalign$$
For the double sum (\AS), we have
$$\align 
&\sum _{h=0} ^{n_1-1}\sum _{s=0} ^{h} (-1)^{h(k+1)}
t^{\min\{s,h-s\}}x^{s(k+1)}y^{(h-s)(k+1)}\\
&\kern4cm 
\cdot
\big(\P_{j+1}(0,h(k+1))
-t(xy)^{k+1}\P_{j+1}(0,(h+2)(k+1))\big)\\
&\quad 
=\sum _{h=-1} ^{n_1-2}\sum _{s=0} ^{h+1} (-1)^{(h+1)(k+1)}
t^{\min\{s,h-s+1\}}x^{s(k+1)}y^{(h-s+1)(k+1)}
\P_{j+1}(0,(h+1)(k+1))\\
&\quad \quad 
-\sum _{h=1} ^{n_1}\sum _{s=1} ^{h} (-1)^{(h-1)(k+1)}
t^{\min\{s,h-s+1\}}x^{s(k+1)}y^{(h-s+1)(k+1)}
\P_{j+1}(0,(h+1)(k+1)).
\endalign$$
In this difference of double sums, almost everything cancels, the
exceptions being the terms for $h=-1$ and for $h=0$ in the first
double sum, the terms for $h\ge1$ and $s=0$ respectively
$s=h+1$ in the first double sum, and the terms for $h=n_1-1$ and for
$h=n_1$ in the second double sum. Thus, we obtain the expression
$$\align 
&\P_{j+1}(0,0)
+(-1)^{k+1}\(x^{k+1}+y^{k+1}\)\P_{j+1}(0,k+1)\\
&\kern1cm
+\sum _{h=1} ^{n_1-2} (-1)^{(h+1)(k+1)}
\(x^{(h+1)(k+1)}+y^{(h+1)(k+1)}\)\P_{j+1}(0,(h+1)(k+1))\\
&\kern1cm
-\sum _{s=1} ^{n_1-1} (-1)^{(n_1-2)(k+1)}
t^{\min\{s,n_1-s\}}x^{s(k+1)}y^{(n_1-s)(k+1)}
\P_{j+1}(0,n_1(k+1))\\
&\kern1cm
-\sum _{s=1} ^{n_1} (-1)^{(n_1-1)(k+1)}
t^{\min\{s,n_1-s+1\}}x^{s(k+1)}y^{(n_1-s+1)(k+1)}
\P_{j+1}(0,(n_1+1)(k+1))\\
&=-\P_{j+1}(0,0)+\sum _{h=0} ^{n_1-1} (-1)^{h(k+1)}
\(x^{h(k+1)}+y^{h(k+1)}\)\P_{j+1}(0,h(k+1))\\
&\kern1cm
-(-1)^{n_1(k+1)}\sum _{s=1} ^{n_1-1} 
t^{\min\{s,n_1-s\}}x^{s(k+1)}y^{(n_1-s)(k+1)}
\delta_{j,n-1}
\tag\AV
\endalign$$
for the double sum in (\AS),
where we used the fact that $j<n=n_1(k+1)$ to see that
$\P_{j+1}(0,(n_1+1)(k+1))=0$ and
$\P_{j+1}(0,n_1(k+1))=\delta_{j,n_1(k+1)-1}=\delta_{j,n-1}$, where $\delta_{a,b}$
denotes the Kronecker delta.

On the other hand, by expanding and shifting indices in (\AT), we obtain
the expression
$$\align 
&\sum _{h=-1} ^{n_1-2}\sum _{s=0} ^{h+1} \sum _{b=1} ^{k}
(-1)^{(h+1)(k+1)+b}
x^{s(k+1)+b}y^{(h-s+1)(k+1)}\\
&\kern5cm
\cdot
t^{\min\{s+\chi(h\text{ even}),h-s+1\}}\P_{j+1}(0,(h+1)(k+1)+b)\\
&+\sum _{h=-1} ^{n_1-2}\sum _{s=-1} ^{h} \sum _{b=1} ^{k}
(-1)^{(h+1)(k+1)+b}
x^{(s+1)(k+1)}y^{(h-s)(k+1)+b}\\
&\kern5cm
\cdot
t^{\min\{s+1,h-s+\chi(h\text{ even})\}}
\P_{j+1}(0,(h+1)(k+1)+b)
\tag\AW
\endalign$$
for the triple sum (\AT), while, by replacing $b$ by $k+1-b$ 
in (\AU), we obtain the expression
$$\align 
&-\sum _{h=0} ^{n_1-1}\sum _{s=0} ^{h}\sum _{b=1} ^{k} 
(-1)^{(h+1)(k+1)-b}
x^{(s+1)(k+1)}y^{(h-s)(k+1)+b}\\
&\kern5cm 
\cdot
t^{\min\{s+1,h-s+\chi(h\text{ even})\}}
\P_{j+1}(0,(h+1)(k+1)+b)\\
&-\sum _{h=0} ^{n_1-1}\sum _{s=0} ^{h}\sum _{b=1} ^{k} 
(-1)^{(h+1)(k+1)-b}
x^{s(k+1)+b}y^{(h-s+1)(k+1)}\\
&\kern5cm 
\cdot
t^{\min\{s+\chi(h\text{ even}),h-s+1\}}
\P_{j+1}(0,(h+1)(k+1)+b)
\tag\AX
\endalign$$
for the triple sum (\AU).
If we add (\AW) and (\AX), then there is again a large amount of
cancellation, with only the terms for $h=-1$, for $s=-1$, and for
$s=h+1$ in (\AW), and the terms for $h=n_1-1$ in (\AX) surviving.
However, the terms for $h=n_1-1$ in (\AX) involve
$\P_{j+1}(0,n_1(k+1)+b)$ which vanishes for $b\ge1$ since $j<n=n_1(k+1)$.
Therefore, the sum of (\AW) and (\AX) equals
$$\align 
&\sum _{b=1} ^{k}
(-1)^b\(x^{b}+y^{b}\)\P_{j+1}(0,b)\\
&\kern.3cm
+\sum _{h=0} ^{n_1-2} \sum _{b=1} ^{k}
(-1)^{(h+1)(k+1)+b}
\(x^{(h+1)(k+1)+b}+y^{(h+1)(k+1)+b}\)\P_{j+1}(0,(h+1)(k+1)+b)\\
&=\sum _{h=0} ^{n_1-1} \sum _{b=1} ^{k}
(-1)^{h(k+1)+b}
\(x^{h(k+1)+b}+y^{h(k+1)+b}\)\P_{j+1}(0,h(k+1)+b).
\tag\AY
\endalign$$
In total, by taking the sum of (\AV) and (\AY), we see that the
$(0,j)$-entry in our new matrix, given in (\ASAS), is equal to
$$\multline 
-\P_{j+1}(0,0)+\sum _{m=0} ^{n_1(k+1)-1} (-1)^{m}
\(x^{m}+y^{m}\)\P_{j+1}(0,m)\\
-(-1)^{n_1(k+1)}\sum _{s=1} ^{n_1-1} 
t^{\min\{s,n_1-s\}}x^{s(k+1)}y^{(n_1-s)(k+1)}
\delta_{j,n-1}.
\endmultline\tag\AZa$$

We now claim that 
$$
-\P_{j+1}(0,0)+\sum _{m=0} ^{j+1} (-1)^{m}
\(x^{m}+y^{m}\)\P_{j+1}(0,m)=0.
\tag\AZ
$$
(The reader should keep in mind that $n=n_1(k+1)$.)
In order to see this, we appeal to (\AGa). Thereby, the
left-hand side of (\AZ) becomes
$$\align 
&\coef{z^0}\Bigg(
-z^{j+1}\(1+\frac {x} {z}\)^{j+1}\(1+\frac {y} {z}\)^{j+1}\\
&\kern3cm
+\sum _{m=0} ^{j+1} (-1)^{m}
\(x^{m}+y^{m}\)z^{j+1-m}\(1+\frac {x} {z}\)^{j+1}\(1+\frac {y} {z}\)^{j+1}
\Bigg)\\
&=\coef{z^0}
z^{j+1}\(1+\frac {x} {z}\)^{j+1}\(1+\frac {y} {z}\)^{j+1}
\Bigg(-1+\frac {1-\(-\frac {x} {z}\)^{j+2}} {1+\frac {x} {z}}
+\frac {1-\(-\frac {y} {z}\)^{j+2}} {1+\frac {y} {z}}
\Bigg)\\
&=\coef{z^0}
z^{j+1}\(1+\frac {x} {z}\)^{j}\(1+\frac {y} {z}\)^{j}
\(1-\frac {xy} {z^2}
\)\\
&\kern4cm
-\coef{z^0}
z^{-1}\(1+\frac {x} {z}\)^{j}\(1+\frac {y} {z}\)^{j}
\((-x)^{j+2}+(-y)^{j+2} 
\)\\
&=\coef{z^{-1}}
\frac {1} {(j+1)}\frac {d} {dz}
\bigg(z^{j+1}\(1+\frac {x} {z}\)^{j+1}\(1+\frac {y} {z}\)^{j+1}
\bigg)=0,
\endalign$$
establishing the claim.

We are now in the position to conclude the proof of (\AJ) in the case
where $n=n_1(k+1)$. By using (\AZ) in (\AZa), we 
see that the $(0,j)$-entry in the new matrix is given by
$$\multline
-(-1)^{n_1(k+1)}\(x^{n_1(k+1)}+y^{n_1(k+1)}\)\delta_{j,n-1}\\
-(-1)^{n_1(k+1)}\sum _{s=1} ^{n_1-1} 
t^{\min\{s,n_1-s\}}x^{s(k+1)}y^{(n_1-s)(k+1)}
\delta_{j,n-1}\\=
-(-1)^{n_1(k+1)}\sum _{s=0} ^{n_1} 
t^{\min\{s,n_1-s\}}x^{s(k+1)}y^{(n_1-s)(k+1)}
\delta_{j,n-1}.
\endmultline$$
In particular, this means that all the entries in row $0$, except for the last,
vanish. It is therefore now easy to compute the determinant of
the matrix we have obtained: as we observed in the paragraph above
(\AO), this matrix is an ``almost" upper triangular matrix, meaning
that it is a matrix $(\hat M_{i,j})$ for which $\hat M_{i,j}=0$ if $j<i-1$.
Furthermore (cf.\ the remark after (\AQ)), we have $\hat
M_{i,i-1}=1$ for $i\ge1$.
Now, in addition, all entries in row $0$ vanish, except for the last,
which is equal to
$$
-(-1)^{n_1(k+1)}\sum _{s=0} ^{n_1} 
t^{\min\{s,n_1-s\}}x^{s(k+1)}y^{(n_1-s)(k+1)}.
$$
Hence, the determinant of the matrix we have obtained equals
the above expression times the product of the entries $\hat
M_{i,i-1}$, $i\ge1$ (which is equal to $1$), times the sign
$(-1)^{n-1}=(-1)^{n_1(k+1)-1}$, that is, it is equal to
$$
\sum _{s=0} ^{n_1} 
t^{\min\{s,n_1-s\}}x^{s(k+1)}y^{(n_1-s)(k+1)}.
\tag\AZAZ$$

We can now sum up. The row operations at the very beginning and the
factorisation of powers of $t$ and $xy$ from the rows resulted in a factor
(see (\AO) and (\AP))
$$(-1)^{n_1\binom {k+1}2}t^{\fl{\frac {n_1} {2}}}
(xy)^{(k+1)^2\binom {n_1}2}.$$
The determinant of the matrix we had obtained after these operations
turned out to be equal to (\AZAZ). The product of these two
expressions is indeed equal to the right-hand side of (\AJ) in the
case where $n=n_1(k+1)$.

\medskip
Finally, we treat the case where $n=n_1(k+1)+k$. In fact, this case can be
reduced to the previous one. Namely, after one has performed the
manipulations described at the beginning, after which the new matrix
entries are given by (\AK)--(\ALa), one is faced with a block matrix
$$\pmatrix A&B\\0&C\endpmatrix,$$
where $A$ is exactly the $(n_1(k+1))\times (n_1(k+1))$ matrix
that is obtained by applying these manipulations in the previous case
(that is, in the case where $n=n_1(k+1)$), and where $C$ is a
$k\times k$ ``reflected upper triangular" matrix. 
(By ``reflected upper triangular", we mean a matrix where all entries above 
the {\it anti}-diagonal of the matrix are equal to $0$.)
Moreover, all entries on the anti-diagonal of the matrix $C$ are equal
to $t^{\chi(n_1\text{ odd})}(xy)^{n_1(k+1)}$. Hence, our determinant
is equal to the result of the previous case 
multiplied by
$$(-1)^{\binom k2}t^{k\cdot\chi(n_1\text{ odd})}(xy)^{n_1k(k+1)}.$$
This is exactly in agreement with the right-hand side of (\AJ) in the
case where $n=n_1(k+1)+k$.\quad \quad \qed
\enddemo

\subhead 6. Main theorems, II\endsubhead
In this section we present two further families of determinant
evaluations, where the entries of the matrices of which the
determinant is taken are built out of path generating functions.
By (\SL), and (\SM), 
the special case $t=1$ of Theorems~\TG\
and \TH\ below implies Theorems~\TAa\ and \TBa, respectively.

The reader should compare these theorems 
with Theorems~\TE\ and \TF. Evidently, there are strong
similarities, with the only essential difference being located
in the first argument of the path generating
function in the second term of the matrix entries. Why we
have chosen to present the matrix entries in Theorems~\TG\ and \TH\ as
sums rather than as differences (as opposed to the presentation of
Theorems~\TE\ and \TF), and why we have chosen a slightly different
exponent of $xy$, will become clear in Section~8. Clearly, one
could transform one presentation into the other by replacement of $t$ by
$-txy$.

\proclaim{Theorem \TG}
For all positive integers $n$ and $k$, we have
$$\multline
\det_{0\le i,j\le n-1}\(\P_{j}(i,k)+t(xy)^{i}\P_{j}(-i,k)\)\\=
\cases (-1)^{kn_1+\binom k2}(1+t)t^{k(n_1-1)}(xy)^{k(n_1-1)(2kn_1-k+1)}&n=2kn_1-k+1,\\
(-1)^{kn_1}(1+t)t^{kn_1-1}(xy)^{kn_1(2kn_1-k-1)}&n=2kn_1,\\
0&n\not\equiv 0,k+1~(\text {\rm mod }2k),
\endcases
\endmultline\tag\BA$$
while for $k=0$ we have
$$
\det_{0\le i,j\le n-1}\(\P_{j}(i,0)+t(xy)^{i}\P_{j}(-i,0)\)=
(1+t)^n(xy)^{\binom n2}.
\tag\BAa$$
\endproclaim

\remark{Remarks}
(1) 
If $k=1$, the first two cases on the right-hand side of (\BA)
coincide, so that we have
$$
\det_{0\le i,j\le n-1}\(\P_{j}(i,1)+t(xy)^{i}\P_{j}(-i,1)\)=
\cases 
(-1)^{n_1}(1+t)t^{n_1-1}(xy)^{2n_1(n_1-1)}&n=2n_1,\\
0&n\text { odd}.
\endcases
$$

\smallskip
(2) By Formula~(\AF), the determinant evaluation in Theorem~\TG\
also implies a formula for negative $k$. More precisely, 
using also (\AE), we have
$$\align 
\P_{j}(i,-k)+t(xy)^{i}\P_{j}(&-i,-k)
=\P_{j}(0,-i-k)+t(xy)^{i}\P_{j}(0,i-k)\\
&=(xy)^{i+k}\P_{j}(0,i+k)+t(xy)^{k}\P_{j}(0,k-i)\\
&=t(xy)^{k}\left(\P_{j}(i,k)
+t^{-1}(xy)^{i}\P_{j}(-i,k)\right).
\tag\BAb
\endalign$$
Aside from some trivial factors, the expression in the last line is
again in the form as the expression for the matrix element of the
determinant on the left-hand side of (\BA).
We omit the explicit statement of the resulting formula.
\endremark

\demo{Proof}
If $k=0$, then the matrix of which we want to compute the determinant 
is upper triangular, and the entry on the
main diagonal in the $i$-th row is
$$
\P_{i}(i,0)+t(xy)^i\P_{i}(-i,0)=(1+t)(xy)^i,
$$
$i=0,1,\dots,n-1$.
The assertion in (\BAa) follows immediately.

The rest of the
proof is analogous to the one of Theorem~\TE. We content ourselves
in outlining the key steps, leaving details to the reader.

In the matrix 
$$
\(\P_{j}(i,k)+t(xy)^{i}\P_{j}(-i,k)\)_{0\le i,j\le n-1}\ ,
$$
of which we want to compute the determinant, we replace row
$(2kh+b)$ by
$$\multline 
\sum _{\ell=0} ^{h}(-1)^{\ell-h}t^{\ell-h}(xy)^{(h-\ell)k}\cdot
\big(\text {row $(2\ell k+b)$}\big)\\
+\sum _{\ell=1} ^{h}(-1)^{\ell-h}t^{\ell-h-1}(xy)^{(h-\ell)k+b}\cdot
\big(\text {row $(2\ell k-b)$}\big)
\endmultline\tag\BC$$
if $0\le b\le k-1$, and by
$$\multline 
\sum _{\ell=0} ^{h}(-1)^{\ell-h}t^{\ell-h}(xy)^{(h-\ell)k}\cdot
\big(\text {row $(2\ell k+b)$}\big)\\
+\sum _{\ell=1} ^{h+1}(-1)^{\ell-h}t^{\ell-h-1}(xy)^{(h-\ell)k+b}\cdot
\big(\text {row $(2\ell k-b)$}\big)
\endmultline\tag\BD$$
if $k+1\le b\le 2k-1$. Again, it is easy to see that this can be achieved by
elementary row manipulations.
(We remind the reader that the rows are labelled by $0,1,2,\dots$.) 
We must pay attention to the fact that, since in the case where $b=0$
and $h\ge1$
the coefficient of row $2hk$ in (\BC) is $1+t^{-1}$, these
manipulations change the value of the determinant. To be precise, they
create a factor of
$$\(1+t^{-1}\)^{\fl{(n-1)/2k}},
\tag\BE$$
by which we must divide the result in the end.

The $(i,j)$-entry of the new matrix, with $i=2hk+b$, is given by
$$
(xy)^{2hk}\P_j(0,(2h+1)k-b)
+t(xy)^{2hk+b}\P_j(0,(2h+1)k+b)
$$
if $0\le b\le k$, and by
$$
-t(xy)^{(2h+1)k}\P_j(0,(2h+3)k-b)
+t(xy)^{2hk+b}\P_j(0,(2h+1)k+b)
$$
if $k+1\le b\le 2k-1$.

In particular, this means that the $(i,j)$-entry, with $i=2hk+b$,
vanishes in the case where $0\le b\le k$ whenever
$j<(2h+1)k-b$, and that it vanishes in the case where $k+1\le b\le 2k-1$
whenever $j<(2h+3)k-b$. 
Hence, if
$n=2hk+b$ with $1\le b\le k$, then row $2hk$ consists entirely of
zeroes since $n-1<2hk+k$. Similarly, if
$n=2hk+b$ with $k+2\le b\le 2k-1$, then row $2hk+k+1$ consists entirely of
zeroes since $n-1<2hk+2k-1$. 
Consequently, the determinant equals zero
in the case where $n\hbox{${}\not\equiv{}0,k+1$}$~(mod~$2k$).

If $n=2n_1k$, then one can transform the matrix which we have
obtained by the above manipulations into an upper triangular matrix,
using the permutation of the rows given by
$$i=2hk+b\mapsto
\cases 
(2h+1)k-b&0\le b\le k,\\
(2h+3)k-b&k+1\le b\le 2k-1,
\endcases
\tag\BF$$
$0\le i\le n-1$. Reading the entries 
along the main diagonal of this upper triangular matrix, we find
$$\align 
&1+t;1,\dots,1;-t(xy)^k,\dots,-t(xy)^k,\\
&(1+t)(xy)^{2k};(xy)^{2k}\dots,(xy)^{2k};-t(xy)^{3k},\dots,-t(xy)^{3k},\\
&(1+t)(xy)^{4k};(xy)^{4k},\dots,(xy)^{4k};-t(xy)^{5k},\dots,-t(xy)^{5k},\\
&\hbox to 9cm{\leaders \hbox{.}\hfil}\\
&(1+t)(xy)^{(2n_1-2)k};(xy)^{(2n_1-2)k},
\dots,(xy)^{(2n_1-2)k};-t(xy)^{(2n_1-1)k},\dots,-t(xy)^{(2n_1-1)k},
\tag\BG\endalign$$
where, when arranged as above, there are exactly $2k$ entries in each
line, the first always containing a factor $1+t$, followed by $k$
equal entries, which are in their turn followed by $k-1$ equal
entries. The product of these entries is
$(1+t)^{n_1}t^{(k-1)n_1}(xy)^{kn_1(2kn_1-k-1)}$. In order to arrive at
the final result, this expression has to be divided by (\BE),
by the sign of the permutation (\BF), and by the signs arising in
(\BG). If everything is put together, we obtain the right-hand side in
(\BA) for $n=2kn_1$.

The case of $n=2kn_1-k+1$ can be treated in the same manner. We leave
the details to the reader.\quad \quad \qed
\enddemo

\proclaim{Theorem \TH}
For all positive integers $n$ and $k$, we have
$$\multline
\det_{0\le i,j\le n-1}\(\P_{j+1}(i,k)+t(xy)^{i}\P_{j+1}(-i,k)\)\\=
\cases 
(-1)^{k(n_1-1)-1}(1+t)t^{kn_1-2}(xy)^{kn_1(2kn_1-k-3)+k}P_{n-k+2,k}(x,y,t)
\hskip-1cm\\
&n=2kn_1-1,\\
(-1)^{kn_1+\binom k2}(1+t)t^{k(n_1-1)}(xy)^{k(n_1-1)(2kn_1-k+1)}P_{n,k}(x,y,t)
\hskip-4cm\\
&n=2kn_1-k+1,\\
(-1)^{kn_1+\binom {k+1}2}(1+t)t^{k(n_1-1)}(xy)^{k(n_1-1)(2kn_1-k-1)}P_{n-k,k}(x,y,t)
\hskip-4cm\\
&n=2kn_1-k,\\
(-1)^{kn_1}(1+t)t^{kn_1-1}(xy)^{kn_1(2kn_1-k-1)}P_{n,k}(x,y,t)&n=2kn_1,\\
0\hskip7cm 
n\not\equiv 0,k,k+1,2k-1~(\text {\rm mod }2k),
\hskip-4cm
\endcases\\
\endmultline\tag\BB$$
where
$$\multline 
P_{m,k}(x,y,t)\\
=\cases 
\sum _{s=0} ^{m/k}
(-1)^s t^{\min\{s,\frac {m} {k}-s\}}
x^{s k}y^{m-s k}&\text{if }m\equiv0~(\text{\rm mod }k),\\
\sum _{s=0} ^{\fl{m/k}}
(-1)^s t^{\min\{s,\cl{m/k}-s\}}
\(x^{s k}y^{m-s k}+
        x^{m-s k}y^{s k}\)&\text{if }m\not\equiv0~(\text{\rm mod }k),
\endcases
\endmultline$$
while for $k=1$ we have
$$\multline
\det_{0\le i,j\le n-1}\(\P_{j+1}(i,1)+t(xy)^{i}\P_{j+1}(-i,1)\)\\=
\cases 
(-1)^{n_1+1}(1+t)t^{n_1-1}(xy)^{2(n_1-1)^2}P_{n-1,1}(x,y,t)
&n=2n_1-1,\\
(-1)^{n_1}(1+t)t^{n_1-1}(xy)^{2n_1(n_1-1)}P_{n,1}(x,y,t)&n=2n_1,\\
\endcases\\
\endmultline\tag\BBa$$
and for $k=0$ we have
$$
\det_{0\le i,j\le n-1}\(\P_{j+1}(i,0)+t(xy)^{i}\P_{j+1}(-i,0)\)=
(1+t)^n(xy)^{\binom n2}(x^n+y^n).
\tag\BBb$$
\endproclaim

\remark{Remark}
Again, by Formula~(\AF), the determinant evaluation in
Theorem~\TH\
also implies formulae for negative $k$. To see this, all one has
to do is to replace $j$ by $j+1$ in the calculation (\BAb).
We omit the explicit statement of the resulting formula.
\endremark

\demo{Proof}
As was the case also earlier, we have to treat the case $k=0$
separately. Using (\AE) and (\AF), we see that, in this case, we have
$$
\P_{j+1}(i,0)+t(xy)^{i}\P_{j+1}(-i,0)=(1+t)(xy)^i\P_{j+1}(0,i).
$$
Hence, the determinant that we want to compute equals
$$(1+t)^n(xy)^{\binom n2}
\det_{0\le i,j\le n-1}\big(\P_{j+1}(0,i)\big).
\tag\BBc$$
We note that, for $i=1,2,\dots,n-1$, 
the first $i-1$ entries in row $i$ of the matrix
$\(\P_{j+1}(0,i)\)_{0\le i,j\le n-1}$ vanish, while
the entry in column $i-1$ equals $\P_i(0,i)=1$.
We now replace row~$0$ by
$$
-(\text {row $0$})+\sum _{i=0} ^{n-1} (-1)^{i}
\(x^{i}+y^{i}\)\cdot(\text {row $i$}).
$$
Since the coefficient of row~$0$ in the above linear combination of
rows equals $1$, this operation does not change the value of the
determinant. Using (\AZ), we see that, in the new $0$-th row, all
entries vanish except for the last one in column $n-1$, which equals
$-(-1)^n(x^n+y^n)$. It is now easy to compute the determinant now
obtained: it is equal to 
$$-(-1)^{n-1}(-1)^n(x^n+y^n)\cdot 1^{n-1}=x^n+y^n.$$
Substituting this for the determinant in (\BBc) leads to the
right-hand side of (\BBb), as required.

From now on let $k\ge1$.
Also here, we have done a similar proof already when establishing
Theorem~\TF. Therefore, again, we shall be brief here.

We start by applying the row operations described in (\BC) and (\BD). 
We obtain a new matrix, where the $(i,j)$-entry, with $i=2hk+b$,
of the new matrix is given by
$$
(xy)^{2hk}\P_{j+1}(0,(2h+1)k-b)
+t(xy)^{2hk+b}\P_{j+1}(0,(2h+1)k+b)
$$
if $0\le b\le k$, and by
$$
-t(xy)^{(2h+1)k}\P_{j+1}(0,(2h+3)k-b)
+t(xy)^{2hk+b}\P_{j+1}(0,(2h+1)k+b)
$$
if $k+1\le b\le 2k-1$.
As earlier, at this point we can already read off that the determinant
vanishes if $n\not\equiv 0,k,k+1,2k-1~(\text {\rm mod }2k)$.

We concentrate now on the case where $n=2kn_1$.
We reorder the rows according to the permutation (\BF). Subsequently,
we divide each entry in row $i$, $i\ge1$, by the first non-zero entry
in its row. Clearly, since this changes the determinant, the
corresponding factor has to be taken into account in the end.

The resulting matrix is again an ``almost" upper triangular matrix,
that is, a matrix $(M_{i,j})$ for which $M_{i,j}=0$ if $j<i-1$.  
Furthermore, the first non-zero entry in row $i$, the entry
$M_{i,i-1}$, equals $1$ for all $i\ge1$.

In this matrix, we replace the $0$-th row by
$$\sum _{h=0} ^{2n_1-1}\sum _{b=0} ^{k-1} (-1)^{hk+b}
\sum _{s=0} ^{h}d(h,b,s)\,x^{sk}y^{(h-s)k}
\cdot
\big(\text {row $(hk+b)$}\big),\tag\BH$$
where the coefficients $d(h,b,s)$ are given by
$$
d(h,b,s)=\cases 
(-1)^{h-s}t^{\min\{s+\chi(h\text{ odd}),h-s\}}x^b+
(-1)^{s}t^{\min\{s,h-s+\chi(h\text{ odd})\}}y^b,\hskip-2cm\\
&\text{if }b\ne0,\\
(-1)^st^{\min\{s,h-s\}},
&\text{if $b=0$ and $h$ is even},
\\
\chi(s=0\text{ or }s=h),
&\text{if $b=0$ and $h$ is odd}.
\endcases
$$
Again, since the coefficient of the $0$-th row in the linear combination
(\BH) is $1$, this does not change the value of the determinant.

We claim that, in the new matrix, all entries in row $0$ are zero,
except for the last one. To compute the last one, one can proceed as
in the analogous situation in the proof of Theorem~\TF. Since there are
no new aspects which arise here, we omit the details, leaving them to
the reader.

\medskip
The remaining three cases can be treated similarly.\quad \quad \qed
\enddemo

\subhead 7. Specialisations\endsubhead
In this section we list specialisations of our results obtained in the
previous sections. The special values of $x$ and $y$ that we choose 
are those that we discussed at the end of Section~2.
We state all of our results for non-negative values of $k$ only.
However, we wish to point out that, for  most of them, 
our results in the previous
sections also imply corresponding results for negative values of $k$,
cf\. the remarks after Theorems~\TAa--\TBa, \TE--\TH. We omit their
explicit statement however for the sake of brevity.

We begin with Theorem~\TA. If we set $x=-y=\sqrt{-1}$ there, then,
using (\SB) and Lemma~\TC, we obtain the following two results.

\proclaim{Corollary \TI}
For all positive integers $n$ and non-negative integers $k$, we have
$$\det_{0\le i,j\le n-1}\(\frac {2k+1} {i+j+k+1}\binom
{2i+2j}{i+j+k}\)=
\cases (-1)^{k n_1+\binom k2}&n=n_1(2k+1)-k,\\
(-1)^{k n_1}&n=n_1(2k+1),\\
0&n\not\equiv 0,k+1~(\text {\rm mod }2k+1).
\endcases
$$
\endproclaim

\proclaim{Corollary \TJ}
For all positive integers $n$ and non-negative integers $k$, we have
$$\det_{0\le i,j\le n-1}\(\frac {2k+1} {i+j+k+2}\binom
{2i+2j+2}{i+j+k+1}\)=
\cases (-1)^{k n_1+\binom {k+1}2}&n=n_1(2k+1)-k-1,\\
(-1)^{k n_1}&n=n_1(2k+1),\\
0&n\not\equiv 0,k~(\text {\rm mod }2k+1).
\endcases
$$
\endproclaim

On the other hand, under the same specialisation, Lemma~\TD\ suggests
that also the determinant
$$\det_{0\le i,j\le n-1}\(\frac {2k} {i+j+k+1}\binom
{2i+2j+1}{i+j+k}\)
$$
is always $0$, $1$, or $-1$. However, since on the right-hand side of
(\SH) we encounter the square of the above determinant, we do not know
whether we get $+1$ or $-1$ for the cases where the determinant is
non-zero. Fortunately, there is a different specialisation which
disposes of this problem, see Corollary~\TL.

Next we specialise $x=\frac {1} {2}(1+\sqrt{-3})$ and $y=\frac {1}
{2}(1-\sqrt{-3})$ in Theorem~\TA. By (\SD), we obtain the following
result.

\proclaim{Corollary \TK}
For all positive integers $n$ and non-negative integers $k$, we have
$$\det_{0\le i,j\le n-1}\(
\sum _{\ell\ge0} ^{}\(\binom {i+j}{\ell,\ell+k}
-\binom {i+j}{\ell,\ell+k+2}\)\)=
\cases 
(-1)^{n_1\binom {k+1}2}&n=n_1(k+1),\\
0&n\not\equiv 0~(\text {\rm mod }k+1).
\endcases
$$
\endproclaim

Finally, the specialisation $x=y=1$ in Theorem~\TA\ yields the following
determinant identity upon appealing to (\SF) and replacing $k$ by
$k-1$.

\proclaim{Corollary \TL}
For all positive integers $n$ and $k$, we have
$$\det_{0\le i,j\le n-1}\(\frac {2k} {i+j+k+1}\binom
{2i+2j+1}{i+j+k}\)=
\cases 
(-1)^{n_1\binom k2}&n=n_1k,\\
0&n\not\equiv 0~(\text {\rm mod }k).
\endcases
$$
\endproclaim

\medskip
Now we turn our attention to Theorem~\TB.
If we set $x=-y=\sqrt{-1}$ there, then,
using (\SB) and Lemma~\TC, in addition to obtaining Corollary~\TJ\
again, we obtain the following result.

\proclaim{Corollary \TM}
For all positive integers $n$ and $k\ge2$, we have
$$\multline
\det_{0\le i,j\le n-1}\(\frac {2k+1} {i+j+k+3}\binom
{2i+2j+4}{i+j+k+2}\)\\=
\cases 
(-1)^{k n_1+\binom {k}2+1}&n=n_1(2k+1)-k-2,\\
(-1)^{k n_1+\binom {k+1}2}(n+k+1)&n=n_1(2k+1)-k-1,\\
(-1)^{k (n_1+1)}(n+1)&n=n_1(2k+1)-1,\\
(-1)^{k n_1}&n=n_1(2k+1),\\
0&n\not\equiv 0,k-1,k,2k~(\text {\rm mod }2k+1),
\endcases
\endmultline
\tag\BHd$$
while for $k=0$ we have
$$
\det_{0\le i,j\le n-1}\(\frac {1} {i+j+3}\binom
{2i+2j+4}{i+j+2}\)=n+1.
$$
\endproclaim

\remark{Remark}
The formula in (\BHd) is also valid for $k=1$. Explicitly, we have
$$
\det_{0\le i,j\le n-1}\(\frac {3} {i+j+4}\binom
{2i+2j+4}{i+j+3}\)=
\cases 
(-1)^{n_1+1}3n_1&n=3n_1-2,\\
(-1)^{n_1+1}3n_1&n=3n_1-1,\\
(-1)^{n_1}&n=3n_1.
\endcases
$$
\endremark

Also here, Lemma~\TD\ suggests a further determinant evaluation, which
turns out to be obtainable by another specialisation,
see Corollary~\TO.

Next we specialise $x=\frac {1} {2}(1+\sqrt{-3})$ and $y=\frac {1}
{2}(1-\sqrt{-3})$ in Theorem~\TB. By (\SD), we obtain the following
result.

\proclaim{Corollary \TN}
For all positive integers $n$ and $k$, we have
$$\multline
\det_{0\le i,j\le n-1}\(
\sum _{\ell\ge0} ^{}\(\binom {i+j+1}{\ell,\ell+k}
-\binom {i+j+1}{\ell,\ell+k+2}\)\)
\\=
\cases 
-(-1)^{n_1\binom {k+2}2}&n=n_1(3k+3)-2k-3\text{ and }k\not\equiv2~(\text{\rm mod }3),\\
-(-1)^{n_1\binom {k+2}2}&n=n_1(3k+3)-2k-2\text{ and }k\not\equiv2~(\text{\rm mod }3),\\
(-1)^{(n_1-1)\binom {k+2}2}&n=n_1(3k+3)-k-2\text{ and }k\not\equiv2~(\text{\rm mod }3),\\
(-1)^{n_1\binom {k+2}2}&n=n_1(3k+3)\text{ and }k\not\equiv2~(\text{\rm mod }3),\\
(-1)^{n_1\binom {k+1}2+\frac {1} {3}n_1(k+1)}(n_1+1)&
n=n_1(k+1)\text{ and }k\equiv2~(\text{\rm mod }3),\\
(-1)^{n_1\binom {k+1}2+\binom k2+\frac {1} {3}n_1(k+1)}(n_1+1)&
n=n_1(k+1)+k\text{ and }k\equiv2~(\text{\rm mod }3),\\
0&\text{otherwise},
\endcases
\endmultline
$$
while for $k=0$ we have
$$
\det_{0\le i,j\le n-1}\(
\sum _{\ell\ge0} ^{}\(\binom {i+j+1}{\ell,\ell+k}
-\binom {i+j+1}{\ell,\ell+k+2}\)\)
=
\cases 
1&n\equiv0,1~(\text{\rm mod }6),\\
-1&n\equiv3,4~(\text{\rm mod }6),\\
0&\text{otherwise}.
\endcases
$$
\endproclaim

Finally, setting $x=y=1$ in Theorem~\TB\ yields the following
determinant identity upon appealing to (\SF) and replacing $k$ by
$k-1$.

\proclaim{Corollary \TO}
For all positive integers $n$ and $k$, we have
$$\multline
\det_{0\le i,j\le n-1}\(\frac {2k} {i+j+k+2}\binom
{2i+2j+3}{i+j+k+1}\)\\
=\cases (-1)^{n_1\binom {k+1}2}(n_1+1)&n=n_1(k+1),\\
(-1)^{n_1\binom {k+1}2+\binom k2}(n_1+1)&n=n_1(k+1)+k,\\
0&n\not\equiv 0,k~(\text {\rm mod }k+1).
\endcases
\endmultline
$$
\endproclaim

\medskip
Next we consider the corresponding specialisations of Theorem~\TAa.
If we set $x=-y=\sqrt{-1}$ there, then,
using (\SA) and Lemma~\TC, we obtain the following two results.

\proclaim{Corollary \TP}
For all positive integers $n$ and $k$, we have
$$\det_{0\le i,j\le n-1}\(\binom {2i+2j} {i+j+k}\)=
\cases 
(-1)^{n_1k}&n=2n_1k,\\
(-1)^{n_1k+\binom {k}2}&n=2n_1k-k+1,\\
0&n\not\equiv 0,k+1~(\text {\rm mod }2k),
\endcases$$
while for $k=0$ we have
$$\det_{0\le i,j\le n-1}\(\binom {2i+2j} {i+j}\)=
2^{n-1}.
\tag\BHdd
$$
\endproclaim

\proclaim{Corollary \TQ}
For all positive integers $n$ and $k$, we have
$$\det_{0\le i,j\le n-1}\(\binom {2i+2j+2} {i+j+k+1}\)=
\cases 
(-1)^{\binom {n_1}2k+n_1\binom k2}&n=n_1k,\\
0&n\not\equiv 0~(\text {\rm mod }k),
\endcases$$
while for $k=0$ we have
$$\det_{0\le i,j\le n-1}\(\binom {2i+2j+2} {i+j+1}\)=
2^{n}.
$$
\endproclaim

Still considering the specialisation $x=-y=\sqrt{-1}$, 
Lemma~\TD\ hints at a further determinant evaluation, which
is given in the theorem below. 

\proclaim{Theorem \TQa}
For all positive integers $n$ and $k$, we have
$$\det_{0\le i,j\le n-1}\(\binom {2i+2j+1} {i+j+k}\)=
\cases 
1&n=(2k-1)n_1,\\
(-1)^{n_1\binom {k}2}&n=(2k-1)n_1-k+1,\\
0&n\not\equiv 0,k~(\text {\rm mod }2k-1),
\endcases\tag\BHa$$
while for $k=0$ we have
$$\det_{0\le i,j\le n-1}\(\binom {2i+2j+1} {i+j}\)=1.
$$
\endproclaim

\remark{Remark}
If $k=1$, the first two cases on the right-hand side of (\BHa)
coincide, so that we have
$$\det_{0\le i,j\le n-1}\(\binom {2i+2j+1} {i+j+1}\)=1.
$$
\endremark

\demo{Sketch of proof of Theorem~\TQa}
Lemma~\TD\ is not sufficient to prove the assertion
 because it only
yields a formula for the square of the determinant in (\BHa). So, we have
to find a direct proof. 

By using the path decomposition argument in the paragraph below (\SK),
or by Chu--Vandermonde convolution, we have
$$\binom {2i+2j+1}{i+j+k}=
\sum _{\ell=-i} ^{i}\binom {2i}{i+\ell}\binom {2j+1}{j+k-\ell}.$$
Then, in the same style as in the paragraph above (\SL), one can do
row manipulations to see that
$$\det_{0\le i,j\le n-1}\(\binom {2i+2j+1} {i+j+k}\)=
\frac {1} {2}\det_{0\le i,j\le n-1}\(\binom {2j+1} {j+k-i}+
\binom {2j+1} {j+k+i}\).$$
In order to evaluate the latter determinant, we can proceed as in the
proof of Theorem~\TG.\quad \quad \qed
\enddemo

If we specialise $x=\frac {1} {2}(1+\sqrt{-3})$ and $y=\frac {1}
{2}(1-\sqrt{-3})$ in Theorem~\TAa, then, by (\SC), we obtain the following
result.

\proclaim{Corollary \TR}
For all positive integers $n$ and $k$, we have
$$
\det_{0\le i,j\le n-1}\(
\sum _{\ell\ge0} ^{}\binom {i+j}{\ell,\ell+k}\)
=
\cases (-1)^{kn_1+\binom k2}&n=2kn_1-k+1,\\
(-1)^{kn_1}&n=2kn_1,\\
0&n\not\equiv 0,k+1~(\text {\rm mod }2k),
\endcases
$$
while for $k=0$ we have
$$
\det_{0\le i,j\le n-1}\(
\sum _{\ell\ge0} ^{}\binom {i+j}{\ell,\ell}\)
=2^{n-1}.
$$
\endproclaim

The specialisation $x=y=1$ in Theorem~\TAa\ yields
the result in Corollary~\TP\ a second time.

\medskip
At last, we consider the corresponding specialisations of Theorem~\TBa.
If we set $x=-y=\sqrt{-1}$ there, then,
using (\SA) and Lemma~\TC, we obtain Corollary~\TQ\ again, but also
the additional determinant evaluation below.

\proclaim{Corollary \TS}
For all positive integers $n$ and $k\ge1$, we have
$$\multline
\det_{0\le i,j\le n-1}\(\binom {2i+2j+4} {i+j+k+2}\)\\=
\cases 
(-1)^{n_1k}&n=2n_1k,\\
(-1)^{n_1k+\binom {k+2}2}&n=2n_1k-k-1,\\
2(-1)^{n_1k+\binom {k+1}2}(n+k)&n=2n_1k-k,\\
2(-1)^{n_1k+k}(n+1)&n=2n_1k-1,\\
0&n\not\equiv 0,k-1,k,2k-1~(\text {\rm mod }2k),
\endcases
\endmultline
\tag\BHe
$$
while for $k=0$ we have
$$
\det_{0\le i,j\le n-1}\(\binom {2i+2j+4} {i+j+2}\)=
2^n(2n+1).
$$
\endproclaim

\remark{Remark}
The formula in (\BHe) is also valid for $k=1$. Explicitly, we have
$$
\det_{0\le i,j\le n-1}\(\binom {2i+2j+4} {i+j+3}\)=
\cases 
(-1)^{n_1}&n=2n_1,\\
(-1)^{n_1+1}4n_1&n=2n_1-1.
\endcases
$$
\endremark

Also here, Lemma~\TD\ hints at a further determinant evaluation.
We formulate it in the conjecture below. 

\proclaim{Conjecture \TSa}
For all positive integers $n$ and $k\ge2$, we have
$$\det_{0\le i,j\le n-1}\(\binom {2i+2j+3} {i+j+k+1}\)=
\cases 
2n_1+1&n=(2k-1)n_1,\\
(-1)^{k+1}(4n_1)&n=(2k-1)n_1-1,\\
(-1)^{\binom {k}2}(4n_1)&n=(2k-1)n_1-k+1,\\
(-1)^{\binom {k-1}2}(2n_1-1)&n=(2k-1)n_1-k,\\
0&\text{otherwise},
\endcases\tag\BHb$$
while for $k=0,1$ we have
$$\det_{0\le i,j\le n-1}\(\binom {2i+2j+3} {i+j+k+1}\)=2n+1.
$$
\endproclaim

\remark{Remark}
We believe that this conjecture can be
proved in a similar way as Theorem~\TQa\ above; that is, one would
first transform the determinant via
$$\det_{0\le i,j\le n-1}\(\binom {2i+2j+3} {i+j+k+1}\)=
\frac {1} {2}\det_{0\le i,j\le n-1}\(\binom {2j+3} {j+k-i+1}+
\binom {2j+3} {j+k+i+1}\),$$
and then proceed in the spirit of the proof of Theorem~\TH. However,
we did not try to work this out. We should point out that, by
Lemma~\TD, the only unproven part in (\BHb) concerns the signs. 
\endremark

If we specialise $x=\frac {1} {2}(1+\sqrt{-3})$ and $y=\frac {1}
{2}(1-\sqrt{-3})$ in Theorem~\TBa, then, by (\SC), we obtain the following
result.

\proclaim{Corollary \TT}
For all positive integers $n$ and $k\ge2$, we have
$$
\det_{0\le i,j\le n-1}\(
\sum _{\ell\ge0} ^{}\binom {i+j+1}{\ell,\ell+k}\)
=
\cases 
(-1)^{kn_1/2}&n=kn_1\text{ and }k\equiv0~(\text{\rm mod }6),\\
(-1)^{\binom {n_1+1}2}&n=kn_1\text{ and }k\equiv3~(\text{\rm mod }12),\\
(-1)^{\binom {n_1}2}&n=kn_1\text{ and }k\equiv9~(\text{\rm mod }12),\\
(-1)^{kn_1+\binom {k+1}2}&n=6kn_1-5k\text{ and }3\nmid k,\\
3(-1)^{k(n_1+1)+\fl{(k+1)/6}}\kern-2pt&n=6kn_1-5k+1\text{ and }3\nmid k,\\
3(-1)^{k(n_1+1)+\fl{k/3}}&n=6kn_1-4k-1\text{ and }3\nmid k,\\
2(-1)^{kn_1+1}&n=6kn_1-4k\text{ and }3\nmid k,\\
2(-1)^{kn_1+\binom {k}2+1}&n=6kn_1-3k\text{ and }3\nmid k,\\
3(-1)^{k(n_1+1)+\fl{(k+4)/6}}\kern-2pt&n=6kn_1-3k+1\text{ and }3\nmid k,\\
3(-1)^{kn_1+\fl{k/3}+1}&n=6kn_1-2k-1\text{ and }3\nmid k,\\
(-1)^{k(n_1+1)}&n=6kn_1-2k\text{ and }3\nmid k,\\
(-1)^{kn_1+\binom {k+1}2}&n=6kn_1-k\text{ and }3\nmid k,\\
(-1)^{kn_1}&n=6kn_1\text{ and }3\nmid k,\\
0&\text{otherwise},
\endcases
$$
while for $k=1$ we have
$$
\det_{0\le i,j\le n-1}\(
\sum _{\ell\ge0} ^{}\binom {i+j+1}{\ell,\ell+1}\)
=
\cases 
1&n\equiv0,1,4,5~(\text{\rm mod }12),\\
2&n\equiv2,3~(\text{\rm mod }12),\\
-1&n\equiv6,7,10,11~(\text{\rm mod }12),\\
-2&n\equiv8,9~(\text{\rm mod }12),
\endcases
$$
and for $k=0$ we have
$$
\det_{0\le i,j\le n-1}\(
\sum _{\ell\ge0} ^{}\binom {i+j+1}{\ell,\ell}\)
=
\cases 
(-8)^{n_1-1}&n=3n_1-2,\\
2(-8)^{n_1-1}&n=3n_1-1,\\
(-8)^{n_1}&n=3n_1.
\endcases
$$
\endproclaim

Finally, the specialisation $x=y=1$ in Theorem~\TBa\ yields
the result in Corollary~\TQ\ a second time.

\medskip
Obviously, we could have also performed analogous specialisations in
Theorems~\TE--\TH, thus obtaining even more determinant evaluations. 
For the sake of brevity, we mention just one of these. It is the
special case $x=y=1$ of Theorem~\TE.

\proclaim{Corollary \TU}
For all positive integers $n$ and non-negative integers $k$, we have
$$\multline
\det_{0\le i,j\le n-1}\(\binom {2j} {j+k-i}-t\binom {2j} {j+k+i+2}\)\\=
\cases (-1)^{n_1\binom
{k+1}2}t^{k\fl{\frac {n_1}2}}&n=n_1(k+1),\\
0&n\not\equiv 0~(\text {\rm mod }k+1).
\endcases
\endmultline$$
\endproclaim

\subhead 8. Concluding remarks and questions\endsubhead
We close our article by some comments on the results that we have
obtained, and by posing some open questions.

\medskip
{\smc 8.1. Is there a connection to symplectic and orthogonal
characters?}
Given a partition 
$\la=(\la_1,\la_2,\dots,\la_n)$ 
(i.e., a non-increasing sequence of non-negative integers), 
the (irreducible) {\it symplectic character}
$\sp_\la(x_1,x_2,\dots,x_n)$ can be
defined by (see \cite{\FuHaAA, Prop.~24.22})
$$\multline
\sp_\la(x_1,x_2,\dots,x_n)\\=
\frac {1} {2}\det_{1\le i,j\le n}
\(h_{\la_i-i+j}(x_1^{\pm1},x_2^{\pm1},\dots,x_n^{\pm1})
+h_{\la_i-i-j+2}(x_1^{\pm1},x_2^{\pm1},\dots,x_n^{\pm1})\),
\endmultline
\tag\BI
$$
where, for $m\ge1$, $h_m(y_1,y_2,\dots,y_N):=
\sum _{1\le i_1\le \dots\le i_m\le N} ^{}y_{i_1}\cdots y_{i_m}$ is the
$m$-th {\it  complete homogeneous symmetric function} in the variables
$y_1,y_2,\dots,y_N$, $h_0(y_1,y_2,\dots,y_N):=1$, and 
$h_m(x_1^{\pm1}, x_2^{\pm1},\dots,x_n^{\pm1})$ is short for
$h_m(x_1,x_1^{-1},x_2,x_2^{-1},\dots,x_n,x_n^{-1})$. It can as well be
written alternatively in the form (see \cite{\FuHaAA, Cor.~24.24})
$$\multline
\sp_\la(x_1,x_2,\dots,x_n)\\=\det_{1\le i,j\le \la_1}
\(e_{\la'_i-i+j}(x_1^{\pm1},x_2^{\pm1},\dots,x_n^{\pm1})
-e_{\la'_i-i-j}(x_1^{\pm1},x_2^{\pm1},\dots,x_n^{\pm1})\),
\endmultline
\tag\BJ
$$
where, for $m\ge1$, $e_m(y_1,y_2,\dots,y_N):=
\sum _{1\le i_1< \dots< i_m\le N} ^{}y_{i_1}\cdots y_{i_m}$ is the
$m$-th {\it elementary symmetric function} in the variables
$y_1,y_2,\dots,y_N$, $e_0(y_1,y_2,\dots,y_N):=1$,
$\la'$ denotes the partition conjugate to $\la$, 
and where we use an analogous convention for
the short form $e_m(x_1^{\pm1},x_2^{\pm1},\dots,x_n^{\pm1})$.
On the other hand, the (irreducible) odd special {\it orthogonal character}
$\so_\la(x_1,x_2,\dots,x_n)$ can be
defined by (see \cite{\FuHaAA, Prop.~24.46})
$$\multline
\so_\la(x_1,x_2,\dots,x_n)\\=
\det_{1\le i,j\le n}
\(h_{\la_i-i+j}(x_1^{\pm1},x_2^{\pm1},\dots,x_n^{\pm1},1)
-h_{\la_i-i-j}(x_1^{\pm1},x_2^{\pm1},\dots,x_n^{\pm1},1)\),
\endmultline
\tag\BK$$
where 
$h_m(x_1^{\pm1},x_2^{\pm1},\dots,x_n^{\pm1},1)$ is short for
$h_m(x_1,x_1^{-1},x_2,x_2^{-1},\dots,x_n,x_n^{-1},1)$. It can as well be
written alternatively in the form (see \cite{\FuHaAA, Cor.~24.35})
$$\multline
\so_\la(x_1,x_2,\dots,x_n)\\=\frac {1} {2}\det_{1\le i,j\le \la_1}
\(e_{\la'_i-i+j}(x_1^{\pm1},x_2^{\pm1},\dots,x_n^{\pm1},1)
+e_{\la'_i-i-j+2}(x_1^{\pm1},x_2^{\pm1},\dots,x_n^{\pm1},1)\),
\endmultline
\tag\BL$$
with an analogous convention how to read
$e_m(x_1^{\pm1},x_2^{\pm1},\dots,x_n^{\pm1},1)$.

If one now compares the right-hand sides of (\BJ) and (\BK) with the determinants in 
Theorems~\TE\ and \TF, respectively the right-hand sides of (\BI) and
(\BL) with the determinants in Theorems~\TG\ and \TH, and if one
recalls that (\BI)--(\BL) have been interpreted in \cite{\FuKrAA} in terms of
generating functions for certain families of non-intersecting lattice
paths, then one observes striking similarities. 
One is led to think that one should be able
to specialise the partition $\la$ and the variables
$x_1,x_2,\dots,x_n$ appropriately so that the determinants in
Theorems~\TE--\TH\  are obtained (at least for $t=1$). However, we
were not able to make this speculation concrete.

\medskip
{\smc 8.2. A combinatorial derivation of (\SL) and (\SM)?}
The determinantal relations (\SJ) and (\SK) were derived by
combinatorial means, making appeal to the Lindstr\"om--Gessel--Viennot
theorem presented here in Theorem~\GV. We could also have derived
these relations by some row manipulations, but we believe that the
combinatorial argument is much more illuminating. On the other hand,
the determinantal relations (\SL) and (\SM) were derived by row
manipulations. This leads naturally to the question whether there are
also combinatorial explanations for (\SL) and (\SM)? Indeed, the
determinants on the left-hand side can be combinatorially
interpreted as generating functions for families of non-intersecting
lattice paths by using again Theorem~\GV. Moreover, there is also a
combinatorial model available for the right-hand side determinants,
which would interpret them as generating functions for families of
non-intersecting paths where, in addition, reflections of paths do
also not intersect other paths (cf\. \cite{\FuKrAA, Sec.~7} for more detailed
explanations on this model). However, we were not able to use these
combinatorial interpretations to develop a combinatorial
understanding of (\SL) and (\SM).

\medskip
{\smc 8.3. Determinant evaluations of E\u gecio\u glu, Redmond and Ryavec.}
In \cite{\EgRRAB, \EgRRAD}, E\u gecio\u glu, Redmond and Ryavec go in
a direction somewhat ``orthogonal" to ours, in that they consider the
Hankel determinants 
$$\det_{0\le i,j\le n-1}\(\binom {2i+2j+k} {i+j}\)$$
(among others). This should be compared to the determinants in
Corollaries~\TP, \TQ\ and \TS. 
E\u gecio\u glu, Redmond and Ryavec develop a complex method based on differential
equations for polynomial generalisations of such determinants, which
enables them to prove closed form evaluations in the cases
$k=0,1,\dots,4$. However, if $k>4$, there do not seem to be ``nice"
formulae for these determinants, as opposed to our families of
determinants. On the other hand, E\u gecio\u glu, Redmond and Ryavec
conjecture (see \cite{\EgRRAD, Sec.~11}) that also their determinants
follow a modular pattern, depending on $n$ and $k$, in general.
Although there is only marginal overlap between their results and
ours, it is still possible that there is a unifying picture for both
sets of results (and conjectures) lurking behind.

\end